\documentclass {article}
\usepackage {amsfonts, amssymb, exscale}
\usepackage [all, 2cell] {xy}
\catcode `\@ 11
\if@compatibility
\else
  \@settopoint\textwidth
\fi
\if@compatibility
\else
  \@settopoint\textwidth
\fi
\if@compatibility
\else
  \ifcase\@ptsize
  \or
  \or
  \fi
\fi
\advance\textheight by \topskip
\let\@oddhead\@empty
\let\@evenhead\@empty
\def\@oddfoot{\normalfont\rightmark \hfil
              \pagename{} \thepage}
\def\@evenfoot{\@oddfoot}
\newcommand\pagename{Page}
\def\@begintheorem#1#2{\trivlist
   \item[\hskip \labelsep{\bfseries #1\ #2:}]\slshape}
\renewenvironment{abstract}{%
    \if@twocolumn
      \section*{\abstractname}%
    \else
      \small
      \begin{center}%
        {\bfseries \abstractname\vspace{-.5em}\vspace{\z@}}%
      \end{center}%
      \quotation
      \noindent \ignorespaces
    \fi}
    {\if@twocolumn\else\endquotation\fi}
\def\today{\number\year\space \ifcase\month\or
  January\or February\or March\or April\or May\or June\or
  July\or August\or September\or October\or November\or December\fi
  \space\number\day}
\def\pr@m@s{%
  \ifx'\@let@token
    \expandafter\pr@@@s
  \else
    \ifcat^\@let@token
      \expandafter\expandafter\expandafter\pr@@@t
    \else
      \ifcat\space\@let@token
        \def~ {\futurelet\@let@token\pr@m@s}%
        \expandafter\expandafter\expandafter\expandafter%
        \expandafter\expandafter\expandafter~%
      \else
        \egroup
      \fi
    \fi
  \fi}
\def\thefootnote{[\@arabic\c@footnote]}
\renewcommand \thepart {\@arabic\c@part}
\renewcommand \thesection {\ifnum \c@part > 0 \thepart.\fi \@arabic\c@section}
\renewcommand\part{%
   \if@noskipsec \leavevmode \fi
   \par
 \c@section \z@
   \addvspace{4ex}%
   \@afterindentfalse
   \secdef\@part\@spart}
\catcode `\@ 12
\newcommand \prefix [1] {$#1$-}
\newcommand \pprefix [1] {($#1$)-}
\newcommand \zero {\prefix 0}
\newcommand \two {\prefix 2}
\newcommand \Two [1] {\two \uppercase {#1}}
\newcommand \too {\pprefix 2}
\renewcommand \deg [1] {degree-$#1$}
\newcommand \D {\prefix D}
\newcommand \G {\prefix G}
\newcommand \gee {\pprefix G}
\newcommand \GGG {\prefix \GG }
\newcommand \Gee {\pprefix \GG }
\renewcommand \H {\prefix H}
\renewcommand \a {a\space }
\newcommand \A {A\space }
\newcommand \BB {{\mathcal B}}
\newcommand \bb {{\mathfrak b}}
\newcommand \CC {{\mathcal C}}
\newcommand \C {{\mathbf C}}
\DeclareMathAlphabet \mathbfrak U {euf} b n
\newcommand \CCC {{\mathbfrak C}}
\newcommand \EE {{\mathcal E}}
\newcommand \ee {{\mathfrak e}}
\newcommand \FF {{\mathcal F}}
\newcommand \ff {{\mathfrak f}}
\newcommand \GG {{\mathcal G}}
\renewcommand \gg {{\mathfrak g}}
\newcommand \I {I\space }
\newcommand \ii {{\mathfrak i}}
\newcommand \jj {{\mathfrak j}}
\newcommand \mm {{\mathfrak m}}
\newcommand \PP {{\mathcal P}}
\newcommand \pp {{\mathfrak p}}
\newcommand \RR {{\mathcal R}}
\newcommand \rr {{\mathfrak r}}
\renewcommand \tt {{\mathfrak t}}
\newcommand \UU {{\mathcal U}}
\newcommand \uu {{\mathfrak u}}
\newcommand \ww {{\mathfrak w}}
\newcommand \XX {{\mathcal X}}
\newcommand \xx {{\mathfrak x}}
\newcommand \YY {{\mathcal Y}}
\newcommand \yy {{\mathfrak y}}
\newcommand \ZZ {{\mathcal Z}}
\newcommand \zz {{\mathfrak z}}
\newcommand \figureplus {+} 
\newcommand \q [1] {\lq {#1}\rq }
\newcommand \strong [1] {\textbf {#1}}
\newenvironment {closeitemize}
 {\itemize \itemsep - \itemsep \relax } {\enditemize }
\catcode `\@ 11
\newcommand \unbig {\futurelet \@let@token \unbigpal }
\newcommand \unbigpal [1]
 {{\ifx \@let@token .\kern \nulldelimiterspace \relax \else
  \ifx \@let@token <\langle \else
  \ifx \@let@token >\rangle \else
  #1\fi \fi \fi }}
\catcode `\@ 12
\newcommand \lr [3] {\mathopen {\unbig #1}#2\mathclose {\unbig #3}}
\newcommand \tuple [1] {\lr <{#1}>}
\newcommand \pair [2] {\tuple {{#1},{#2}}}
\newcommand \triple [3] {\tuple {{#1},{#2},{#3}}}
\newcommand \quadruple [4] {\tuple {{#1},{#2},{#3},{#4}}}
\newtheorem {prop} {Proposition}
\newtheorem {thm} {Theorem}
\newenvironment {prf}
 {\trivlist \item [\hskip \labelsep {\bfseries Proof:}]}
 {\hfil \strut \hfill $\blacksquare $\endtrivlist }
\newcommand \id [1] {\underline {#1}}
\newcommand \iid [1] {\id {#1}}
\newcommand \idid [1] {\iid {\id #1}}
\newcommand \src [1] {\lr |{#1}|}
\newcommand \longto \longrightarrow
\newcommand \longTo \Longrightarrow
\newcommand \maplongto [1] {\mathrel {\mathop {\longrightarrow }\limits ^{#1}}}
\newcommand \maplongTo [1] {\mathrel {\mathop {\Longrightarrow }\limits ^{#1}}}
\newcommand \matrixxy [1] {\matrix {\xymatrix @!0{#1}}}
\newcommand \dricell [1] {\drtwocell \omit {\omit \textstyle #1}}
\newcommand \dreqcell {\drtwocell \omit {=}}
\newcommand \driicell [1] {\drtwocell \omit {#1}}
\newcommand \drloweriicell [2] {\drlowertwocell _{{#1}}{#2}}
\newcommand \drupperiicell [2] {\druppertwocell ^{{#1}}{#2}}
\newcommand \drfulliicell [3] {\drtwocell ^{{#1}}_{{#2}}{#3}}
\newcommand \ddricell [1] {\ddrtwocell \omit {\omit \textstyle #1}}
\newcommand \ddreqcell {\ddrtwocell \omit {=}}
\newcommand \ddriicell [1] {\ddrtwocell \omit {#1}}
\newcommand \rrloweriicell [2] {\rrlowertwocell _{{#1}}{#2}}
\newcommand \rrupperiicell [2] {\rruppertwocell ^{{#1}}{#2}}
\newcommand \rrfulliicell [3] {\rrtwocell ^{{#1}}_{{#2}}{#3}}
\newcommand \drricell [1] {\drrtwocell \omit {\omit \textstyle #1}}
\newcommand \drreqcell {\drrtwocell \omit {=}}
\newcommand \drriicell [1] {\drrtwocell \omit {#1}}
\newcommand \ddrricell [1] {\ddrrtwocell \omit {\omit \textstyle #1}}
\newcommand \ddrreqcell {\ddrrtwocell \omit {=}}
\newcommand \ddrriicell [1] {\ddrrtwocell \omit {#1}}
\newcommand \dstring [1] {\ar @{-}[dd]^{\textstyle \mathstrut #1}}
\newcommand \ddstring [1] {\ar @{-}[dddd]^{\textstyle \mathstrut #1}}
\newcommand \dddstring [1] {\ar @{-}[dddddd]^{\textstyle \mathstrut #1}}
\newcommand \ddddstring [1] {\ar @{-}[dddddddd]^{\textstyle \mathstrut #1}}
\newcommand \ustring [1] {\dstring {\isub {#1}}}
\newcommand \uustring [1] {\ddstring {\isub {#1}}}
\newcommand \uuustring [1] {\dddstring {\isub {#1}}}
\newcommand \uuuustring [1] {\ddddstring {\isub {#1}}}
\newcommand \estring {\ar @{.}[dd]}
\newcommand \dStringl [1] {\ar @{=}[dd]_{\textstyle #1}}
\newcommand \ddStringl [1] {\ar @{=}[dddd]_{\textstyle #1}}
\newcommand \dddStringl [1] {\ar @{=}[dddddd]_{\textstyle #1}}
\newcommand \ddddStringl [1] {\ar @{=}[dddddddd]_{\textstyle #1}}
\newcommand \dString {\ar @{=}[dd]}
\newcommand \ddString {\ar @{=}[dddd]}
\newcommand \dddString {\ar @{=}[dddddd]}
\newcommand \rstring [1] {\ar @{--}[rr]|{\textstyle #1}}
\newcommand \rrstring [1] {\ar @{--}[rrrr]|{\textstyle #1}}
\newcommand \rrrstring [1] {\ar @{--}[rrrrrr]|{\textstyle #1}}
\newcommand \rrrrstring [1] {\ar @{--}[rrrrrrrr]|{\textstyle #1}}
\newcommand \rString {\ar @{--}[rr]}
\newcommand \rrString {\ar @{--}[rrrr]}
\newcommand \rrrString {\ar @{--}[rrrrrr]}
\newcommand \pstring [1] {{\textstyle #1}} 
\newcommand \st {\;{\mathstrut \smash {\lower 0.5 ex \hbox {\vdots }}}\;}
\newcommand \famb [2] {\lr ({{#1}\st {#2}})}
\newcommand \setb [2] {\lr \{{{#1}\st {#2}}\}}
\newcommand \faml [1] {\lr ({#1})}
\newcommand \nerve [2] {#1^{\lr [{#2}]}}
\newcommand \cat [2] {#1^{#2}}
\newcommand \nervecat [3] {#1^{#2\lr [{#3}]}}
\newcommand \del [2] {#1_{\lr [{#2}]}}
\newcommand \dom [2] {#1_{#2}}
\newcommand \deldom [3] {#1_{\lr [{#2}]#3}}
\newcommand \domdel [3] {#1_{#2\lr [{#3}]}}
\newcommand \lnr [2] {\mathord {}{#1}\mathclose {#2}}
\newcommand \restrict [2] {\lnr {#1}|_{#2}}
\newcommand \lnmrn [3] {\mathord {}{#1}\mathord {#2}{#3}\mathord {}}
\newcommand \fract [2] {\lnmrn {#1}/{#2}}
\newcommand \slice [2] {\fract {#1}{#2}}
\newcommand \Bun [4] {\mathbf {Bun}_{#1}\lr ({#2,#3,#4})}
\newcommand \Grp {{\mathbf {Grp}}}
\newcommand \bsub [2] {b_{{#1}{#2}}}
\newcommand \bbsub [2] {\bb _{{#1}{#2}}}
\newcommand \bpsub [2] {b'_{{#1}{#2}}}
\newcommand \bbpsub [2] {\bb '_{{#1}{#2}}}
\newcommand \bcbpsub [2] {\lr ({b;b'})_{{#1}{#2}}}
\newcommand \bbcbbpsub [2] {\lr ({\bb ;\bb '})_{{#1}{#2}}}
\newcommand \esub 1
\newcommand \eesub 1 
\newcommand \gsub [2] {g_{{#1}{#2}}}
\newcommand \ggsub [2] {\gg _{{#1}{#2}}}
\newcommand \gpsub [2] {g'_{{#1}{#2}}}
\newcommand \ggpsub [2] {\gg '_{{#1}{#2}}}
\newcommand \isub [1] {#1^{-1}}
\newcommand \iisub [1] {#1^{-1}} 
\newcommand \msub [2] {\lr ({{#1}{#2}})}
\newcommand \mmsub [2] {\lr ({{#1}{#2}})} 
\newcommand \jtildesub [2] {{#1}/_{\!#2}}
\newcommand \jjtildesub [2] {{#1}/_{\!#2}} 
\newcommand \jptildesub [2] {{#1}/_{\!#2}} 
\newcommand \jjptildesub [2] {{#1}/_{\!#2}} 
\newcommand \jpptildesub [2] {{#1}/_{\!#2}} 
\newcommand \rsub [2] {\lr ({{#1}{#2}})} 
\newcommand \rrsub [2] {\lr ({{#1}{#2}})} 
\newcommand \alphasub [3] {\alpha \lr ({{#1}{#2}{#3}})}
\newcommand \betasub [3] {\beta _{{#1}{#2}{#3}}}
\newcommand \gammasub [3] {\gamma _{{#1}{#2}{#3}}}
\newcommand \gammapsub [3] {\gamma '_{{#1}{#2}{#3}}}
\newcommand \deltasub [3] {\delta _{{#1}{#2}{#3}}}
\newcommand \ndeltasub [3] {\bar \deltasub {#1}{#2}{#3}}
\newcommand \epsilonsub [1] {\epsilon \lr ({#1})}
\newcommand \zetasub [3] {\zeta {#1}_{{#2}{#3}}}
\newcommand \zetatildesub [3] {\tilde \zeta {#1}_{{#2}{#3}}}
\newcommand \etasub [1] {\eta _{#1}}
\newcommand \etapsub [1] {\eta '_{#1}}
\newcommand \thetasub [3] {\theta {#1}_{{#2}{#3}}}
\newcommand \iotasub [1] {\iota \lr ({#1})}
\newcommand \lambdasub [1] {\lambda \lr ({#1})}
\newcommand \musub [3] {\mu \lr ({{#1}{#2}{#3}})}
\newcommand \xisub [2] {\xi _{{#1}{#2}}}
\newcommand \rhosub [1] {\rho \lr ({#1})}
\newcommand \sigmasub [3] {\sigma _{{#1}{#2}{#3}}}
\newcommand \upsilonsub [1] {\upsilon \lr ({#1})}
\newcommand \Aut [1] {\mathop {\underline {\mathit {Aut}}}\lr ({#1})}
\newcommand \semitimes {\rtimes }
\newcommand \imestimes {\ltimes }
\begin {document}
\UseAllTwocells
\title {Higher gauge theory I:\\
\Two bundles}
\author {Toby Bartels\footnote
{I reserve no copyright or patent rights to this work;
see \texttt {http://toby.bartels.name/copyright/}.}\\
\texttt {toby\figureplus 2431488868@math.ucr.edu}\\
Department of Mathematics\\
University of California\\
Riverside CA 92521\\
USA}
\date {2006 June 24}
\maketitle
\begin {abstract}
\I categorify the definition of fibre bundle,
replacing smooth manifolds with differentiable categories,
Lie groups with coherent Lie \two groups,
and bundles with \a suitable notion of \two bundle.
To link this with previous work,
\I show that certain \two categories of principal \two bundles
are equivalent to certain \two categories of (nonabelian) gerbes.
This relationship can be (and has been) extended
to connections on \two bundles and gerbes.
\par The main theorem, from \a perspective internal to this paper,
is that the \two category of \two bundles
over \a given \two space under \a given \two group
is (up to equivalence) independent of the fibre
and can be expressed in terms of cohomological data (called \two transitions).
From the perspective of linking to previous work on gerbes,
the main theorem is
that when the \two space is the \two space corresponding to \a given space
and the \two group is the automorphism \two group of \a given group,
then this \two category
is equivalent to the \two category of gerbes over that space under that group
(being described by the same cohomological data).
\end {abstract}
\tableofcontents
\part * {Introduction}
Gauge theory, built on the theory of fibre bundles,
has underlain all verified theories of basic physics for half \a century.
Geometrically, there is \a direct relationship
between \a connection on \a principal bundle under \a group
---which (at least relative to local coordinates)
assigns group elements to paths in spacetime---
and the effects of \a physical field on particles
that (in the Feynman interpretation) travel paths in spacetime.
On \a more purely mathematical side,
studying bundles over \a given space \( B \) under \a given group \( G \)
allows one to approach the \G valued \deg 1 cohomology of \( B \);
this can be also be expressed sheaf-theoretically,
which can be generalised, ultimately,
to \deg 1 cohomology valued in sheaves of groups.
Indeed, everything done with bundles can be done with sheaves.
\par More recently, much of this has been categorified (see below) with gerbes.
This began on the mathematical side,
with \deg 2 cohomology valued in (nonabelian) groups (or sheaves thereof);
but there is also hope for physical applications.
\A categorified theory of bundles
could describe physical theories of strings,
which are (in spacetime) \two dimensional.
Before this work began,
several authors \cite {BM, Attal, HP}
have begun discussing connections on gerbes.
However, gerbes are special stacks,
which in turn are \a categorification of sheaves;
they are not directly \a categorification of bundles,
and bundles are more familiar to working physicists.
The goal of this paper is to develop,
analogously to \a direct development of fibre bundles,
\a theory of categorified fibre bundles.
Rather than invent \a new term, \I call these \emph {\two bundles}.
\par As sheaves are more general than bundles,
so gerbes are more general than \two bundles.
However, \I also have the opportunity here
to generalise from groups to \two groups.
\Two groups~\cite {HDA5} are \a categorification of groups,
although any group may be made into \a \two group,
the \emph {automorphism \two group} of the group.
Previous work on gerbes implicitly uses \two groups,
but usually only automorphism \two groups of groups.
Abelian gerbes use \a different method (that applies only to abelian groups)
of turning \a group into \a \two group.
Lawrence Breen \cite {crosmodgerb}
and Branislav Jur\v co \cite {crosmodbungerb}
have discussed using crossed modules,
which include the above as special cases.
In \a theory of \two bundles, however,
it is natural to use general \two groups from the beginning;
indeed, \I will use
the most general (properly weakened) \emph {coherent} \two groups.
(Crossed modules are equivalent to \q {strict} \two groups.)
Thus, this paper is at once
both more general and less general than work on gerbes;
however, it is to be hoped that gerbes themselves
will also be generalised to (stacks of) coherent \two groups in the future.
\par Categorifying the first paragraph of this Introduction,
one expects to say that higher gauge theory,
built on \a theory of \two bundles,
can underlie stringy theories of physics:
\begin {quote}
Geometrically, there is \a direct relationship
between \a connection on \a principal \two bundle under \a \two group
---which (at least relative to local coordinates)
assigns \two group object-elements to paths in spacetime
and \two group morphism-elements to surfaces in spacetime---
and the effects of \a physical field on (open) strings
that travel surfaces between paths in spacetime.
\end {quote}
Although this paper does not go into connections on \two bundles,
John Baez and Urs Schreiber have begun that theory
(based on earlier versions of this paper) in \cite {HGT2} and \cite {HGT}.
Schreiber has even applied this to physics
in his doctoral dissertation~\cite {Urs}.
\par There is another approach to gerbes
that, like this one, comes closer to the traditional view of bundles:
the concept of bundle gerbe
(introduced in \cite {bungerb},
generalised to arbitrary groups in \cite {nonabbungerb},
and generalised to crossed modules in \cite {crosmodbungerb}).
This approach removes all (explicit) category theory,
both the category theory involved in extending bundles to sheaves
and the categorification involved in raising to a higher dimension.
Like them, \I also downplay the first aspect,
but the second is my prime motivator,
a concept of central importance to this exposition.
\par In general, \emph {categorification}~\cite {CF, categorification}
is \a project in which mathematical structures are enriched
by replacing equations between elements of \a set
with specific isomorphisms between objects in \a category.
The resulting categorified concept is generally richer,
because \a given pair of isomorphic objects
might have many isomorphisms between them
(whose differences were previously ignored),
and other interesting morphisms (previously invisible)
might not be invertible at all.
The classic historical example
is the revolution in homology theory,
begun by Leopold Vietoris and Emmy Noether
at the beginning of the 20th century
(see \cite {Noether} for \a brief history),
in which they replaced Betti numbers (elements of \a set)
with homology groups (objects of \a category)
and studied induced group homomorphisms.
\par Previous work using only \emph {strict} Lie \two groups
(as by Breen~\cite {crosmodgerb} in the guise of crossed modules
and by Baez~\cite {gauge.tex} explicitly)
violates the spirit of categorification,
since the axioms for \a strict \two group
require \emph {equations} between certain objects;
\a more natural definition from the perspective of categorification
should replace these with \emph {isomorphisms},
and then enrich the theory
by positing coherence laws between those isomorphisms.
The \emph {coherent} \two groups do precisely this; \I use them throughout.
Although not directly relevant to this paper,
Baez and Alissa Crans~\cite {HDA6}
have developed \a corresponding theory of Lie \two algebras;
here one finds
\a canonical \prefix 1parameter family
of categorifications of \a Lie algebra---
but these categorifications must be so called \q {weak} Lie \two algebras
(which correspond to coherent Lie \two groups).
\part {Uncategorified bundles} \label 1
To begin with, \I review the theory of bundles with \a structure group,
which is what \I intend to categorify.
Some of this exposition covers rather elementary material;
my purpose is to summarise that material in purely arrow-theoretic terms,
to make clear how it is to be categorified in part~\ref 2.
\section {Categorical preliminaries} \label {1prelim}
Here, \I review the necessary category theory,
which will become \two category theory in section~\ref {2prelim}.
\par In theory, \I work with \a specific category \( \C \),
to be described in section~\ref {1C}.
(It's no secret;
it's just the category of smooth manifolds and smooth functions.
But \I don't want to focus attention on that fact yet.)
However, \I will not normally make reference
to the specific features of this category;
\I will just refer to it as the category \( \C \).
Thus the ideas in this paper
generalise immediately to other categories,
such as categories of topological spaces, algebraic varieties, and the like.
All that matters
is that \( \C \) support the categorical operations described below.
\subsection {Notation and terminology} \label {1term}
\A \strong {space} is an object in \( \C \),
and \a \strong {map} is \a morphism in \( \C \).
Uppercase italic letters like \q {\( X \)} denote spaces,
and lowercase italic letters like \q {\( x \)} denote maps.
However, the identity morphism on \( X \) will be denoted \q {\( \id X \)}.
\par In general, the names of maps will label arrows directly,
perhaps in \a small inline diagram like \( X \maplongto x Y \maplongto y Z \),
which denotes \a composition of maps,
or perhaps in \a huge commutative diagram.
\I will endeavour to make such diagrams easy to read
by orienting maps to the right when convenient, or if not then downwards.
Sometimes \I will place \a lowercase Greek letter like \q {\( \psi \)}
inside \a commutative diagram,
which at this stage is merely \a convenient label for the diagram.
(When \a small commutative diagram, previously discussed,
appears inside of \a large, complicated commutative diagram later on,
these labels can help you keep track of what's what.)
When \I categorify the situation in part~\ref 2, however,
the meaning of these labels will be upgraded to natural transformations.
\par Except in section~\ref {1C} (at least within part~\ref 1),
my discussion is in purely arrow-theoretic terms,
making no direct reference to points, open subsets, and the like.
However, to make the connection to more set-theoretic presentations,
there is \a useful generalised arrow-theoretic concept of element.
Specifically, an \strong {element} \( x \) of \a space \( Y \)
is simply \a space \( X \) together with \a map \( X \maplongto x Y \).
In certain contexts, the map \( x \) really has the same information in it
as an ordinary (set-theoretic) element of (the underlying set of) \( Y \).
On the other hand, if \( X \) is chosen differenlty,
then the map \( x \) may described more complicated features, like curves.
The most general element is in fact the identity map
\( X \maplongto { \id X } X \).
You (as reader) may imagine the elements as set-theoretic points if that helps,
but their power in proofs lies in their complete generality.
\subsection {Products} \label {1prod}
Given spaces \( X \) and \( Y \),
there must be \a space \( X \times Y \),
the \strong {Cartesian product} of \( X \) and \( Y \).
One can also form more general Cartesian products,
like \( X \times Y \times Z \) and so on.
In particular, there is \a \strong {singleton space} \( 1 \),
which is the Cartesian product of no spaces.
\par The generalised elements of \a Cartesian product
may be taken to be ordered pairs;
that is, given maps \( X \maplongto x Y \) and \( X \maplongto { x ' } Y ' \),
there is \a \strong {pairing map}
\( X \maplongto { \pair x { x ' } } Y \times Y ' \);
conversely, given \a map \( X \maplongto y Y \times Y ' \),
there is \a unique pair of maps
\( X \maplongto x Y \) and \( X \maplongto { x ' } Y ' \)
such that \( y = \pair x { x ' } \).
Similarly, given maps
\( X \maplongto x Y \), \( X \maplongto { x ' } Y ' \),
and \( X \maplongto { x ' ' } Y ' ' \),
there is \a map
\( X \maplongto { \triple x { x ' } { x ' ' } } Y \times Y ' \times Y ' ' \);
conversely, given \a map \( X \maplongto y Y \times Y ' \times Y ' ' \),
there is \a unique triple of maps
\( X \maplongto x Y \), \( X \maplongto { x ' } Y \),
and \( X \maplongto { x ' ' } Y ' ' \)
such that \( y = \triple x { x ' } { x ' ' } \).
And so on.
Additionally, given any space \( X \),
there is \a \strong {trivial map} \( X \maplongto { \hat X } 1 \);
conversely, every map \( X \maplongto y 1 \) must be equal to \( \hat X \).
(That is, the singleton space has \a unique generalised element,
once you fix the source \( X \).)
Inverting this notation,
it's sometimes convenient
to abbreviate the map
\( X \maplongto { \pair { \id X } { \id X } } X \times X \)
as \( X \maplongto { \check X } X \times X \).
\par \I will only occasionally want to use the pairing maps directly;
more often, I'll be working with product maps.
Specifically, given maps \( X \maplongto x Y \)
and \( X ' \maplongto { x ' } Y ' \),
there is \a \strong {product map}
\( X \times X ' \maplongto { x \times x ' } Y \times Y ' \).
This respects composition;
given \( X \maplongto x Y \maplongto y Z \)
and \( X ' \maplongto { x ' } Y ' \maplongto { y ' } Z \),
the product map \( X \times X ' \longto Z \times Z ' \)
is \( X \times X ' \maplongto { x \times x ' } Y \times Y '
\maplongto { y \times y ' } Z \times Z ' \).
This also respects identities;
the product map of \( X \maplongto { \id X } X \)
and \( X ' \maplongto { \id { X ' } } X \)
is the identity map
\( X \times X ' \maplongto { \id { X \times X ' } } X \times X ' \).
Similarly, given maps
\( X \maplongto x Y \), \( X ' \maplongto { x ' } Y ' \),
and \( X ' ' \maplongto { x ' ' } Y ' ' \),
there is \a product map
\( X \times X ' \times X ' '
\maplongto { x \times x ' \times x ' ' } Y \times Y ' \times Y ' ' \)
with similar nice properties;
and so on.
The product maps and pairing maps fit together nicely:
given \( X \maplongto { \pair x { x ' } } Y \times Y ' \),
\( x \) and \( x ' \) may be reconstructed
as (respectively) \( X \maplongto { \pair x { x ' } } Y \times Y '
\maplongto { \id Y \times \hat { Y ' } } Y \)
and \( X \maplongto { \pair x { x ' } } Y \times Y '
\maplongto { \hat Y \times \id { Y ' } } Y ' \);
and the product map \( X \times X ' \maplongto { x \times x ' } Y \times Y ' \)
is equal to \( \pair { x \times \hat { X ' } } { \hat X \times x ' } \).
\par \( \C \) is \a monoidal category under \( \times \),
so the Mac Lane Coherence Theorem~\cite [VII.2] {CWM} applies,
allowing me to treat \( X \times 1 \) and \( X \) as the same space,
\( X \maplongto { \check X } X \times X
\maplongto { X \times \check X } X \times X \times X \)
and \( X \maplongto { \check X } X \times X
\maplongto { \check X \times X } X \times X \times X \)
as the same well defined map,
and so on.
By Mac Lane's theorem,
there will always be \a canonical functor that links any identified pair;
and these will all fit together properly.
\I will make such identifications in the future without further comment.
(In fact, \I already did so towards the end of the previous paragraph.)
\subsection {Pullbacks} \label {1pull}
Unlike products, pullbacks~\cite [III.4] {CWM}
need not always exist in \( \C \).
Thus, \I will have to treat this in more detail,
so that \I can discuss exactly what it means for \a pullback to exist.
\par \A \strong {pullback diagram}
consists of spaces \( X \), \( Y \), and \( Z \),
and maps \( X \maplongto x Z \) and \( Y \maplongto y Z \):
\begin {equation} \label {pullback diagram}
\matrixxy { && X \ar [dd] ^ x \\\\ Y \ar [rr] _ y && Z }
\end {equation}
Given \a pullback diagram, \a \strong {pullback cone}
is \a space \( C \)
together with maps \( C \maplongto z X \)
and \( C \maplongto w Y \)
that make the following diagram commute:
\begin {equation} \label {pullback cone}
\matrixxy
{ C \ar [dd] _ w \ar [rr] ^ z \ddrricell \omega && X \ar [dd] ^ x \\\\
Y \ar [rr] _ y && Z }
\end {equation}
(Recall that here \( \omega \) is just \a name for the commutative diagram,
to help keep track of it.)
Given pullback cones \( C \) and \( C ' \),
\a \strong {pullback cone morphism} from \( C \) to \( C ' \)
is \a map \( C \maplongto u C ' \)
that makes the following diagrams commute:
\begin {equation} \label {pullback cone morphism}
\matrixxy
{ C \ar `d ^dr [dddr] _ w [dddr] \ar [dr] ^ u \\
\relax \dricell \psi & C ' \ar [dd] ^ { w ' } \\ & \\ & Y }
\qquad
\matrixxy
{ \\ C \ar [dr] _ u \ar `r _dr [drrr] ^ z [drrr] & \relax \dricell \chi \\
& C ' \ar [rr] _ { z ' } && X }
\end {equation}
\par \A \strong {pullback} of the given pullback diagram
is \a pullback cone \( P \)
that is \emph {universal}
in the sense that, given any other pullback cone \( C \),
there is \a unique pullback cone morphism
from \( C \) to \( P \).
Pullbacks always exist in, for example, the category of sets
(as the subset of \( X \times Y \)
consisting of those pairs
whose left and right components are mapped to the same element of \( Z \)).
However, they don't always exist for the category \( \C \),
so \I will need to make sure that they exist in the relevant cases.
\par In the rest of this paper,
\I will want to define certain spaces as pullbacks, when they exist.
If such definitions are to be sensible,
then \I must show that it doesn't matter which pullback one uses.
\begin {prop} \label {pullback isomorphism}
Given the pullback diagram (\ref {pullback diagram})
and pullbacks \( P \) and \( P ' \),
the spaces \( P \) and \( P ' \) are isomorphic in \( \C \).
\end {prop}
\begin {prf}
Since \( P \) is universal,
there is \a pullback cone morphism \( P ' \maplongto u P \):
\begin {equation} \label {pullback morphism}
\matrixxy
{ P '
\ar `d ^dr [dddr] _ { w ' } [dddr] \ar [dr] ^ u
\ar `r _dr [drrr] ^ { z ' } [drrr] &
\relax \dricell \chi \\
\relax \dricell \psi & P \ar [dd] ^ w \ar [rr] _ z \ddrricell \omega &&
X \ar [dd] ^ x \\ & \\
& Y \ar [rr] _ y && Z }
\end {equation}
Since \( P ' \) is universal,
there is also \a pullback cone morphism \( P \maplongto { \bar u } P ' \):
\begin {equation} \label {competing pullback morphism}
\matrixxy
{ P
\ar `d ^dr [dddr] _ w [dddr] \ar [dr] ^ { \bar u }
\ar `r _dr [drrr] ^ z [drrr] &
\relax \dricell { \bar \chi } \\
\relax \dricell { \bar \psi } &
P ' \ar [dd] ^ { w ' } \ar [rr] _ { z ' } \ddrricell { \omega ' } &&
X \ar [dd] ^ x \\ & \\
& Y \ar [rr] _ y && Z }
\end {equation}
Composing these one way,
\I get \a pullback cone morphism from \( P \) to itself:
\begin {equation} \label {composed pullback automorphism}
\matrixxy
{ P
\ar `d ^r [ddddrr] _ w [ddddrr] \ar [dr] ^ { \bar u }
\ar `r _d [ddrrrr] ^ z [ddrrrr] &&
\relax \dricell { \bar \chi } \\
& P '
\ar `d ^dr [dddr] _ { w ' } [dddr] \ar [dr] ^ u
\ar `r _dr [drrr] ^ { z ' } [drrr] &
\relax \dricell \chi & \\
\relax \dricell { \bar \psi } & \relax \dricell \psi &
P \ar [dd] ^ w \ar [rr] _ z \ddrricell \omega && X \ar [dd] ^ x \\ && \\
&& Y \ar [rr] _ y && Z }
\end {equation}
But since the identity on \( P \) is also \a pullback cone morphism,
the universal property
shows that \( P \maplongto { \bar u } P ' \maplongto u P \)
is the identity on \( P \).
Similarly, \( P ' \maplongto u P \maplongto { \bar u } P ' \)
is the identity on \( P ' \).
Therefore, \( P \) and \( P ' \) are isomorphic.
\end {prf}
\subsection {Equivalence relations} \label {1eqrel}
Given \a space \( U \),
\a \strong {binary relation} on \( U \)
is \a space \( \nerve R 2 \) equipped with maps
\( \nerve R 2 \maplongto { \del j 0 } U \)
and \( \nerve R 2 \maplongto { \del j 1 } U \)
that are \emph {jointly monic};
this means that given generalised elements
\( X \maplongto x \nerve R 2 \) and \( X \maplongto y \nerve R 2 \)
of \( \nerve R 2 \),
if these diagrams commute:
\begin {equation} \label {relation}
\matrixxy
{ X \ar [dd] _ y \ar [drr] ^ x \\
\relax \drricell { \del \chi 0 } &&
\relax \nerve R 2 \ar [dd] ^ { \del j 0 } \\
\relax \nerve R 2 \ar [drr] _ { \del j 0 } && \\ && U }
\qquad
\matrixxy
{ X \ar [ddr] _ y \ar [rr] ^ x & \relax \ddricell { \del \chi 1 } &
\relax \nerve R 2 \ar [ddr] ^ { \del j 1 } \\\\
& \relax \nerve R 2 \ar [rr] _ { \del j 1 } & & U }
\end {equation}
then \( x = y \).
The upshot of this is that an element of \( \nerve R 2 \)
is determined by two elements of \( U \),
or equivalently by \a single element of \( U \times U \);
in other words, the set of elements (with given domain) of \( \nerve R 2 \)
is \a subset of the set of elements of \( U \times U \).
(Thus the link with the set-theoretic notion of binary relation.)
\par \A binary relation is \strong {reflexive}
if it is equipped with \a \emph {reflexivity map}
\( U \maplongto { \del j { 0 0 } } \nerve R 2 \)
that is \a section of both \( \del j 0 \) and \( \del j 1 \);
that is, the following diagrams commute:
\begin {equation} \label {reflexivity}
\matrixxy
{ U \ar `d ^dr [dddr] _ { \id U } [dddr] \ar [dr] ^ { \del j { 0 0 } } \\
\relax \dricell { \del \omega 0 } &
\relax \nerve R 2 \ar [dd] ^ { \del j 0 } \\ & \\
& U }
\qquad
\matrixxy
{ U \ar [dr] _ { \del j { 0 0 } } \ar `r _dr [drrr] ^ { \id U } [drrr] &
\relax \dricell { \del \omega 1 } \\
& \relax \nerve R 2 \ar [rr] _ { \del j 1 } && U }
\end {equation}
\I say \q {equipped with} rather than \q {such that there exists},
but \a reflexivity map (if extant) is unique by the relation's joint monicity.
In terms of generalised elements,
given any element \( X \maplongto x U \) of \( U \),
composing with \( \del j { 0 0 } \)
gives an element of \( \nerve R 2 \);
the commutative diagrams state
that this corresponds to the element \( \pair x x \) of \( U \times U \).
(Thus the link with the set-theoretic notion of reflexivity.)
\par Assume that the kernel pair of \( \del j 0 \) exists,
and let it be the space \( \nerve R 3 \):
\begin {equation} \label {equivalence relation triple}
\matrixxy
{ \relax \nerve R 3
\ar [dd] _ { \del j { 0 1 } } \ar [drr] ^ { \del j { 0 2 } } \\
\relax \drricell { \del \omega { 0 0 } } &&
\relax \nerve R 2 \ar [dd] ^ { \del j 0 } \\
\relax \nerve R 2 \ar [drr] _ { \del j 0 } && \\ && U }
\end {equation}
\A binary relation is \strong {right Euclidean}
if it is equipped with \a (right) \emph {Euclideanness map}
\( \nerve R 3 \maplongto { \del j { 1 2 } } \nerve R 2 \)
making these diagrams commute:
\begin {equation} \label {Euclideanness}
\matrixxy
{ \relax \nerve R 3
\ar [dd] _ { \del j { 0 1 } } \ar [rr] ^ { \del j { 1 2 } }
\ddrricell { \del \omega { 0 1 } } &&
\relax \nerve R 2 \ar [dd] ^ { \del j 1 } \\\\
\relax \nerve R 2 \ar [rr] _ { \del j 0 } && U }
\qquad
\matrixxy
{ \relax \nerve R 3
\ar [ddr] _ { \del j { 0 2 } } \ar [rr] ^ { \del j { 1 2 } } &
\relax \ddricell { \del \omega { 1 1 } } &
\relax \nerve R 2 \ar [ddr] ^ { \del j 1 } \\\\
& \relax \nerve R 2 \ar [rr] _ { \del j 1 } & & U }
\end {equation}
Like the reflexivity map, any Euclideanness map is unique.
In terms of elements,
each element of \( \nerve R 3 \)
is uniquely defined by elements \( x \), \( y \), and \( z \) of \( U \)
such that \( \pair x y \) and \( \pair x z \)
give elements of \( \nerve R 2 \).
The Euclideanness condition states precisely
that under these circumstances,
\( \pair y z \) also gives an element of \( \nerve R 2 \).
As Euclid \cite [CN1] {Euclid} put it (in Heath's translation),
referring to magnitudes of geometric figures,
\q {Things which are equal to the same thing are also equal to one another.};
hence the name.
\par Putting these together,
an \strong {equivalence relation}
is \a binary relation that is both reflexive and Euclidean.
(The usual definition uses symmetry and transitivity instead of Euclideanness;
in the presence of reflexivity, however, these are equivalent.
Using Euclideanness
---more to the point, avoiding explicit mention of symmetry---
helps me to maintain \a consistent orientation on all diagrams;
it also privleges \( \del j 0 \),
which seems bad at first
but is actually desirable
in light of the axioms for covers in section~\ref {1cover}.)
\par Given an equivalence relation as above,
\a \strong {quotient} of the equivalence relation
is \a space \( \nerve R 0 \) and \a map \( U \maplongto j \nerve R 0 \)
that \emph {coequalises} \( \del j 0 \) and \( \del j 1 \).
This means, first, that this diagram commutes:
\begin {equation} \label {quotient}
\matrixxy
{ \relax \nerve R 2
\ar [dd] _ { \del j 0 } \ar [rr] ^ { \del j 1 } \ddrricell \omega &&
U \ar [dd] ^ j \\\\
U \ar [rr] _ j && \relax \nerve R 0 }
\end {equation}
and also that, given any map \( U \maplongto x \nerve R 0 \)
such that this diagram commutes:
\begin {equation} \label {quotient cone}
\matrixxy
{ \relax \nerve U 2
\ar [dd] _ { \del j 0 } \ar [rr] ^ { \del j 1 } \ddrricell { \omega _ x } &&
U \ar [dd] ^ x \\\\
U \ar [rr] _ x && X }
\end {equation}
there is \a \emph {unique} map \( \nerve R 0 \maplongto { \tilde x } X \)
such that this diagram commutes:
\begin {equation} \label {quotient universality}
\matrixxy
{ U \ar [dr] _ j \ar `r _dr [drrr] ^ x [drrr] &
\relax \dricell { \tilde \omega _ x } \\
& \relax \nerve R 0 \ar [rr] _ { \tilde x } && X }
\end {equation}
Thus if an equivalence relation has \a quotient,
then \I can define \a map out of this quotient
by defining \a map (with certain properties) out of the relation's base space.
\subsection {Covers} \label {1cover}
\I will require \a notion of (open) cover
expressed in \a purely arrow-theoretic way.
For me, \strong {covers}
will be certain maps \( U \maplongto j B \) in the category \( \C \).
I'll describe which maps these are in section~\ref {1opencover}
(although again it's no secret:
they're the surjective local diffeomorphisms);
for now, \I list and examine
the category-theoretic properties that covers must have.
\par Here are the properties of covers that \I will need:
\begin {closeitemize}
\item All isomorphisms are covers;
\item \A composite of covers is \a cover;
\item The pullback of \a cover along any map exists and is \a cover;
\item The quotient of every equivalence relation involving \a cover
exists and is \a cover;
and \item Every cover is the quotient of its kernel pair.
\end {closeitemize}
(Compare \cite [B.1.5.7] {Elephant}.)
\par \I will be particularly interested
in the pullback of \( U \) along itself,
that is the pullback of this diagram:
\begin {equation} \label {cover kernel pair diagram}
\matrixxy { & & U \ar [dd] ^ j \\ \\ U \ar [rr] _ j & & B }
\end {equation}
This pullback, \( \nerve U 2 \), is the \strong {kernel pair} of \( U \);
it amounts to an open cover of \( B \)
by the pairwise intersections of the original open subsets.
In the pullback diagram:
\begin {equation} \label {cover kernel pair}
\matrixxy
{ \relax \nerve U 2
\ar [dd] _ { \del j 1 } \ar [rr] ^ { \del j 0 } \ddrricell \omega & &
U \ar [dd] ^ j \\ \\
U \ar [rr] _ j & & B }
\end {equation}
\( \del j 0 \) may be interpreted as the (collective) inclusion
of \( U _ k \cap U _ { k ' } \) into \( U _ k \),
while \( \del j 1 \) is the inclusion
of \( U _ k \cap U _ { k ' } \) into \( U _ { k ' } \).
But \I will not refer to these individual intersections,
only to the kernel pair \( \nerve U 2 \) as \a whole.
\par Notice that this commutative diagram:
\begin {equation} \label {cover kernel pair diagonal cone}
\matrixxy
{ U \ar [dd] _ { \id U } \ar [rr] ^ { \id U } \ddrreqcell &&
U \ar [dd] ^ j \\\\
U \ar [rr] _ j && B }
\end {equation}
is \a pullback cone for~(\ref {cover kernel pair}),
so the universal property of that pullback
defines \a map \( U \maplongto { \del j { 0 0 } } \nerve U 2 \).
This map may be thought of as
the inclusion of \( U _ k \) into \( U _ k \cap U _ k \).
(But don't let this trick you into thinking
that \( \del j { 0 0 } \) is an isomorphism;
in general, most elements of \( \nerve U 2 \) are not of that form.)
This universal map comes with the following commutative diagrams:
\begin {equation} \label {cover kernel pair diagonal}
\matrixxy
{ U \ar `d ^dr [dddr] _ { \id U } [dddr] \ar [dr] ^ { \del j { 0 0 } } \\
\relax \dricell { \del \omega 1 } &
\relax \nerve U 2 \ar [dd] ^ { \del j 1 } \\ & \\
& U }
\qquad
\matrixxy
{ U \ar [dr] _ { \del j { 0 0 } } \ar `r _dr [drrr] ^ { \id U } [drrr] &
\relax \dricell { \del \omega 0 } \\
& \relax \nerve U 2 \ar [rr] _ { \del j 0 } & & U }
\end {equation}
\par In one instance, \I will even need the space \( \nerve U 3 \),
built out of the triple intersections
\( U _ k \cap U _ { k ' } \cap U _ { k ' ' } \).
This space may also be defined as \a pullback:
\begin {equation} \label {cover kernel triple}
\matrixxy
{ \relax \nerve U 3
\ar [dd] _ { \del j { 1 2 } } \ar [rr] ^ { \del j { 0 1 } }
\ddrricell { \del \omega { 0 1 } } & &
\relax \nerve U 2 \ar [dd] ^ { \del j 1 } \\ \\
\relax \nerve U 2 \ar [rr] _ { \del j 0 } & & U }
\end {equation}
\par Notice that this commutative diagram:
\begin {equation} \label {cover kernel triple outer map cone}
\matrixxy
{ \relax \nerve U 3
\ar [dd] _ { \del j { 1 2 } } \ar [rr] ^ { \del j { 0 1 } }
\ddrricell { \del \omega { 0 1 } } & &
\relax \nerve U 2 \ar [dd] _ { \del j 1 } \ar [ddrr] ^ { \del j 0 } \\
& & \relax \ddrricell \omega \\
\relax \nerve U 2 \ar [ddrr] _ { \del j 1 } \ar [rr] ^ { \del j 0 } &
\relax \ddrricell \omega &
U \ar [ddrr] ^ j & & U \ar [dd] ^ j \\ & & & & \\
& & U \ar [rr] _ j & & B }
\end {equation}
is \a pullback cone for~(\ref {cover kernel pair}),
so the universal property of that pullback
defines \a map \( \nerve U 3 \maplongto { \del j { 0 2 } } \nerve U 2 \).
This map may be thought of
as the inclusion of \( U _ k \cap U _ { k ' } \cap U _ { k ' ' } \)
into \( U _ k \cap U _ { k ' ' } \).
This universal map comes with the following commutative diagrams:
\begin {equation} \label {cover kernel pair triple outer map}
\matrixxy
{ \relax \nerve U 3
\ar [dd] _ { \del j { 1 2 } } \ar [rr] ^ { \del j { 0 2 } }
\relax \ddrricell { \del \omega { 1 1 } } & &
\relax \nerve U 2 \ar [dd] ^ { \del j 1 } \\ \\
\relax \nerve U 2 \ar [rr] _ { \del j 1 } & & U }
\qquad
\matrixxy
{ \relax \nerve U 3
\ar [dd] _ { \del j { 0 2 } } \ar [rr] ^ { \del j { 0 1 } }
\relax \ddrricell { \del \omega { 0 0 } } & &
\relax \nerve U 2 \ar [dd] ^ { \del j 0 } \\ \\
\relax \nerve U 2 \ar [rr] _ { \del j 0 } & & U }
\end {equation}
These squares are in fact pullbacks themselves.
\section {The category of spaces} \label {1C}
Here \I want to explain how the category of smooth manifolds
has the properties needed from section~\ref {1prelim}.
This section (and this section alone)
is not analogous to the corresponding section in part~\ref 2.
\par \A \strong {space} (object) is \a smooth manifold,
and \a \strong {map} (morphism in \( \C \)) is \a smooth function.
\A composite of smooth functions is smooth,
and the identity function on any smooth manifold is smooth,
so this forms \a category \( \C \).
An isomorphism in this category is \a diffeomorphism.
\par If the space \( X \) consists of \a single point,
then the map \( X \maplongto x Y \) has the same information in it
as an ordinary (set-theoretic) element of \( Y \).
Furthermore, if two maps from \( Y \) to \( Z \)
agree on each of these elements, then they are equal.
Thus when \I speak of generalised elements,
you (as reader) may think of set-theoretic elements;
everything will make perfect sense.
\subsection {Products and pullbacks} \label {1Clims}
The singleton space \( 1 \)
is the connected \zero dimensional manifold,
which consists of \a single point.
Every function to this space is smooth,
and there is always \a unique function from any space to \( 1 \).
Similarly, the product of smooth manifolds \( X \) and \( Y \) exists,
although \I won't spell it in so much detail.
\par Pullbacks are trickier, because the don't always exist.
For example, in the pullback diagram~(\ref {pullback diagram}),
take \( X \) to be the singleton space,
\( Z \) to be the real line, and \( Y \) to be the plane \( Z \times Z \);
let \( x \) be the constant map to the real number \( 0 \),
and let \( y \) be multiplication.
Then the pullback of this diagram, if it existed,
would the union of the two axes in \( Y \)
---but that is not \a manifold.
(It is possible to generalise the notion of manifold;
Baez \& Schreiber try this in \cite [Definition~4ff] {HGT}.
However, it's enough here that certain pullbacks exist.)
\par When pullbacks do exist, they are the spaces of solutions to equations.
That is why the counterexample above
should have been the set of solutions to the equation \( \msub x y = 0 \)
(where \( x \) and \( y \) are elements of \( Z \),
so that \( \pair x y \) is an element of \( Y \cong Y \times 1 \)).
\A useful example is where one of the maps in the pullback
is the inclusion of \a subset;
then the pullback is the preimage of that subset along the other map.
If the subset is open, then this always exists (is \a manifold)
and always is itself an open subset.
This is the first part of the idea of cover;
but we need to modify this so that covers are regular epimorphisms
(which inclusion of \a proper open subsets never is).
\subsection {Open covers} \label {1opencover}
The usual notion of open cover of the space \( B \)
is an index set \( K \)
and \a family \( \famb { U _ k } { k \in K } \) of open subsets of \( B \),
with union \( \bigcup _ k U _ k = B \).
This family can be combined
into the disjoint union
\( U : = \setb { \lr ( { k , x } ) \in K \times B } { x \in U _ k } \)
and the single map \( U \maplongto j B \)
given by \( \lr ( { k , x } ) \mapsto x \).
Notice that this map \( j \) is \a surjective local diffeomorphism.
\par Conversely, consider
any surjective local diffeomorphism \( U \maplongto j B \).
Since \( j \) is \a \emph {local} diffeomorphism,
there exists an open cover \( \famb { U _ k } { k \in K } \) of \( U \)
such that \( j \) becomes an open embedding when restricted to any \( U _ k \).
Since \( j \) is surjective,
the family \( \famb { j \lr ( { U _ k } ) } { k \in K } \)
is therefore an open cover of \( B \).
Also note that if \I start with an open cover of \( B \),
turn it into \a surjective local diffeomorphism,
and turn that back into an open cover,
then \I get back the same open cover that \I started with.
\par Thus, \I will define \a \strong {cover}
to be any surjective local diffeomorphism \( U \maplongto j B \).
It seems that \cite [7.C] {JST}
was the first to use this generality as \a specification for \( j \).
(Note that their \q {\( E \)} and \q {\( p \)}
refer to my \( U \) and \( j \);
while they call my \( E \) and \( p \)
---to be introduced in section~\ref {1bun} below---
\q {\( A \)} and \q {\( \alpha \)}.)
\par You could reinterpret this paper with slightly different notions of cover.
For example, you could insist that covers always be the form given above,
\a map from (a space diffeomorphic to) \a disjoint union of open subsets;
but that wouldn't make any significant difference.
\A more interesting possibility might be to take all surjective submersions
(without requiring them to be immersions as well).
Compare the uses of these in~\cite {surjsubm}.
\section {Bundles} \label {1bun}
The most general notion of bundle is quite simple;
if \( B \) is \a space,
then \a \strong {bundle} over \( B \)
is simply \a space \( E \) together with \a map \( E \maplongto p B \).
Before long, of course,
\I will restrict attention to bundles with more structure than this---
eventually, to locally trivial bundles with \a structure group.
\par \I will generally denote the bundle by the same name as the space \( E \),
and the default name for the map is always \q {\( p \)}---
with subscripts, primes, or other decorations as needed.
In diagrams, \I will try to orient \( p \) at least partially downward,
to appeal to the common intuition of \( E \) as lying \q {over} \( B \).
\subsection {The category of bundles} \label {1C/B}
For \a proper notion of equivalence of bundles,
\I should define the category \( \slice \C B \) of bundles over \( B \).
\par Given bundles \( E \) and \( E ' \) over \( B \),
\a \strong {bundle morphism} from \( E \) to \( E ' \)
is \a map \( E \maplongto f E ' \)
making the following diagram commute:
\begin {equation} \label {bundle morphism}
\matrixxy
{ E \ar [dddr] _ p \ar [rr] ^ f \ddrricell \pi & &
E ' \ar [dddl] ^ { p ' } \\ \\
& & \\ & B }
\end {equation}
Given bundle morphisms
\( f \) from \( E \) to \( E ' \) and \( f ' \) from \( E ' \) to \( E ' ' \),
the composite bundle morphism
is simply the composite map \( E \maplongto f E ' \maplongto { f ' } E ' ' \);
notice that this diagram commutes:
\begin {equation} \label {composite bundle morphism}
\matrixxy
{ E \ar [dddrr] _ p \ar [rr] ^ f & \relax \ddricell \pi &
E '\ar [ddd] ^ { p ' } \ar [rr] ^ { f ' } \ddricell { \pi ' } & &
E ' ' \ar [dddll] ^ { p ' ' } \\ \\ & & & \\
& & B }
\end {equation}
Also, the identity bundle morphism on \( E \) is the identity map on \( E \).
\begin {prop} \label {bundle cat}
Given \a space \( B \),
bundles over \( B \) and their bundle morphisms form \a category.
\end {prop}
\begin {prf}
\( \slice \C B \) is just
\a slice category~\cite [II.6] {CWM} of the category \( \C \) of spaces.
\end {prf}
\par \I now know what it means
for bundles \( E \) and \( E ' \) to be \strong {equivalent bundles}:
isomorphic objects in the category \( \slice \C B \).
There must be \a bundle morphism from \( E \) to \( E ' \)
and \a bundle morphism from \( E ' \) to \( E \)
whose composite bundle morphisms, in either order,
are identity bundle morphisms.
What this amounts to, then,
are maps \( E \maplongto f E ' \) and \( E ' \maplongto { \bar f } E \)
making the following diagrams commute:
\begin {equation} \label {bundle equivalence}
\matrixxy
{ E \ar [dddr] _ p \ar [rr] ^ f \ddrricell \pi & &
E ' \ar [dddl] ^ { p ' } \\ \\ & & \\
& B }
\qquad
\matrixxy
{ E ' \ar [dddr] _ { p ' } \ar [rr] ^ { \bar f } \ddrricell { \bar \pi } & &
E \ar [dddl] ^ p \\ \\ & & \\
& B }
\qquad
\matrixxy
{ E \ar [rrrr] ^ { \id E } \ar [drr] _ f & \relax \drricell \kappa & & & E \\
& & E ' \ar [urr] _ { \bar f } & }
\qquad
\matrixxy
{ E ' \ar [rrrr] ^ { \id { E ' } } \ar [drr] _ { \bar f } &
\relax \drricell { \bar \kappa } & & & E ' \\
& & E \ar [urr] _ f & }
\end {equation}
In particular, \( E \) and \( E ' \)
are equivalent as spaces (diffeomorphic).
\subsection {Trivial bundles} \label {1trivbun}
If \( B \) and \( F \) are spaces,
then the Cartesian product \( F \times B \)
is automatically \a bundle over \( B \).
Just let the map be \( F \times B \maplongto { \hat F \times \id B } B \),
the projection that forgets about \( F \).
\par Of this bundle, it can rightly be said that \( F \) is its fibre.
I'll want to define \a more general notion, however,
of bundle over \( B \) with fibre \( F \).
To begin, I'll generalise this example to the notion of \a trival bundle.
Specifically, \a \strong {trivial bundle} over \( B \) with fibre \( F \)
is simply \a bundle over \( B \) equipped with
\a bundle equivalence to it from \( F \times B \).
\par In more detail,
this is \a space \( E \)
equipped with maps \( E \maplongto p B \),
\( F \times B \maplongto f E \),
and \( E \maplongto { \bar f } F \times B \)
making the following diagrams commute:
\begin {equation} \label {trivial bundle}
\matrixxy
{ F \times B
\ar [dddr] _ { \hat F \times \id B } \ar [rr] ^ f \ddrricell \pi & &
E \ar [dddl] ^ p \\ \\ & & \\
& B }
\qquad
\matrixxy
{ E \ar [dddr] _ p \ar [rr] ^ { \bar f } \ddrricell { \bar \pi } & &
F \times B \ar [dddl] ^ { \hat F \times \id B } \\ \\ & & \\
& B }
\qquad
\matrixxy
{ F \times B \ar [rrrr] ^ { \id F \times \id B } \ar [drr] _ f &
\relax \drricell \kappa & & & F \times B \\
& & E \ar [urr] _ { \bar f } & }
\qquad
\matrixxy
{ E \ar [rrrr] ^ { \id E } \ar [drr] _ { \bar f } &
\relax \drricell { \bar \kappa } & & & E \\
& & F \times B \ar [urr] _ f & }
\end {equation}
\subsection {Restrictions of bundles} \label {1restrict}
Ultimately, I'll want to deal with \emph {locally} trivial bundles,
so \I need \a notion of restriction to an open cover.
In this section, I'll work quite generally,
where the open cover is modelled by \emph {any} map \( U \maplongto j B \).
\I will call such \a map an (unqualified) \strong {subspace} of \( B \);
this is just \a radical generalisation
of terminology like \q {immersed submanifold},
which can be seen sometimes in differential geometry~\cite {submanifold}.
(Thus logically, there is no difference between
\a subspace of \( B \) and \a bundle over \( B \);
but one does different things with them,
and they will be refined in different ways.)
\par Given \a bundle \( E \)
and \a subspace \( U \) (equipped with the map \( j \) as above),
\I get \a pullback diagram:
\begin {equation} \label {restriction diagram}
\matrixxy { & & E \ar [dd] ^ p \\ \\ U \ar [rr] _ j & & B }
\end {equation}
Then the \strong {restriction} \( \restrict E U \) of \( E \) to \( U \)
is \emph {any} pullback of this diagram, if any exists.
\I name the associated maps as in this commutative diagram:
\begin {equation} \label {restriction}
\matrixxy
{ \relax \restrict E U
\ar [dd] _ { \restrict p U } \ar [rr] ^ { \tilde j }
\ddrricell { \tilde \omega } & &
E \ar [dd] ^ p \\ \\
U \ar [rr] _ j & & B }
\end {equation}
Notice that \( \restrict E U \) becomes
both \a subspace of \( E \) and \a bundle over \( U \).
\par By the unicity of pullbacks,
the restriction \( \restrict E U \), if it exists,
is well defined up to diffeomorphism;
but is it well defined up to equivalence of bundles over \( U \)?
The answer is yes,
because the pullback cone morphisms
in diagram (\ref {composed pullback automorphism})
become bundle morphisms when applied to this situation.
\subsection {Locally trivial bundles} \label {1loctrivbun}
\I can now combine the preceding ideas
to define \a locally trivial bundle with \a given fibre,
which is the general notion of fibre bundle without \a fixed structure group
(compare \cite [I.1.1, I.2.3] {Steenrod}).
\par Given spaces \( B \) and \( F \),
suppose that \( B \) has been supplied with \a subspace \( U \maplongto j B \),
and let this subspace be \a \emph {cover}
in the sense of section~\ref {1cover}.
In this case, the restriction \( \restrict E U \)
(like any pullback along \a cover map)
must exist.
\A \strong {locally trivial bundle}
over \( B \) with fibre \( F \) subordinate to \( U \)
is \a bundle \( E \) over \( B \)
such that the trivial bundle \( F \times U \) over \( U \)
is isomorphic to \( \restrict E U \).
In fact, since the restriction is defined to be \emph {any} pullback,
\I can simply take \( F \times U \) to \emph {be} \( \restrict E U \),
by specifying the appropriate \( \tilde j \).
\par So to sum up, \a locally trivial bundle consists of the following items:
\begin {closeitemize}
\item \a \strong {base space} \( B \);
\item \a \strong {cover space} \( U \);
\item \a \strong {fibre space} \( F \);
\item \a \strong {total space} \( E \);
\item \a \strong {cover map} \( U \maplongto j B \);
\item \a \strong {projection map} \( E \maplongto p B \);
and \item \a \strong {pulled-back map}
\( F \times U \maplongto { \tilde j } E \);
\end {closeitemize}
such that this diagram is \a pullback:
\begin {equation} \label {summary pullback}
\matrixxy
{ F \times U
\ar [dd] _ { \hat F \times \id U } \ar [rr] ^ { \tilde j }
\ddrricell { \tilde \omega } & &
E \ar [dd] ^ p \\ \\
U \ar [rr] _ j & & B }
\end {equation}
\par In studying fibre bundles,
\I regard \( B \) as \a fixed structure on which the bundle is defined.
The fibre \( F \) (like the group \( G \) in the next section)
indicates the type of problem that one is considering;
while the cover \( U \) (together with \( j \))
is subsidiary structure that is not preserved by bundle morphisms.
Accordingly, \( E \) (together with \( p \)) \q {is} the bundle;
\I may use superscripts or primes on \( p \) (as well as \( U \) and \( j \))
if I'm studying more than one bundle.
\par If \( X \) is any space,
\( X \maplongto w F \) is \a generalised element of \( F \),
and \( X \maplongto x U \) is \a generalised element of \( U \),
then let \( \jtildesub w x \)
be the composite
\( X \maplongto { \pair w x } F \times U \maplongto { \tilde j } E \).
(The mnemonic is that \( \jtildesub w x \)
is the interpretation of the point \( w \) in the fibre
\emph {over} the point given in local coordinates as \( x \).)
This notation will be convenient in the sequel.
\section {\G spaces} \label {1Gspace}
Although the theory of groups acting on sets is well known~\cite [I.5] {Lang},
\I will review it in strictly arrow-theoretic terms
to make the categorification clearer.
\subsection {Groups} \label {1group}
The arrow-theoretic definition of group
may be found for example in \cite [page~2f] {CWM}.
\par First, \a \strong {monoid} is \a space \( G \)
together with maps
\( 1 \maplongto e G \) (the identity, or unit, element)
and \( G \times G \maplongto m G \) (multiplication),
making the following diagrams commute:
\begin {equation} \label {monoid}
\matrixxy
{ G \times G \times G
\ar [dd] _ { \id G \times m } \ar [rrr] ^ { m \times \id G } &
\relax \ddricell \alpha && G \times G \ar [dd] ^ m \\\\
G \times G \ar [rrr] _ m &&& G }
\qquad
\matrixxy
{ G \ar [ddrrr] _ { \id G } \ar [rrr] ^ { e \times \id G } &&
\relax \ddricell \lambda & G \times G \ar [dd] ^ m \\\\
&&& G }
\qquad
\matrixxy
{ G \ar [ddrrr] _ { \id G } \ar [rrr] ^ { \id G \times e } &&
\relax \ddricell \rho & G \times G \ar [dd] ^ m \\\\
&&& G }
\end {equation}
These are (respectively) the \emph {associative law}
and the \emph {left and right unit laws}
(hence the letters chosen to stand for the diagrams).
\par In terms of generalised elements,
if \( x \) and \( y \) are elements of \( G \)
(that is, \( X \maplongto x G \) and \( X \maplongto y G \)
are maps for some arbitrary space \( X \)),
then \I will write \( \msub x y \)
for the composite \( X \maplongto { \pair x y } G \times G \maplongto m G \).
Also, I'll write \( \esub \)
for \( X \maplongto { \hat X } 1 \maplongto e G \)
(suppressing \( X \)).
Using this notation,
if \( x \), \( y \), and \( z \) are elements of \( G \),
then the associative law
says that \( \msub { \msub x y } z = \msub x { \msub y z } \),
and the unit laws
say (respectively) that \( \msub \esub x = x \)
and \( \msub x \esub = x \).
\par Then \a \strong {group} is \a monoid \( G \)
together with \a map \( G \maplongto i G \) (the inverse operator)
such that the following diagrams commute:
\begin {equation} \label {inverse}
\matrixxy
{ G \ar [dd] _ { \hat G } \ar [rrr] ^ { \pair i { \id G } } &
\relax \ddricell \epsilon &&
G \times G \ar [dd] ^ m \\ \\
1 \ar [rrr] _ e &&& G }
\qquad
\matrixxy
{ G \ar [dd] _ { \pair { \id G } i } \ar [rrr] ^ { \hat G } &
\relax \ddricell \iota &&
1 \ar [dd] ^ e \\ \\
G \times G \ar [rrr] _ m &&& G }
\end {equation}
These are (respectively) the \emph {left and right inverse laws}.
\par In terms of elements,
if \( X \maplongto x G \) is an element of \( G \),
then I'll write \( X \maplongto x G \maplongto i G \) as \( \isub x \).
Then the inverse laws
say (respectively) that \( \msub { \isub x } x = \esub \)
and \( \esub = \msub x { \isub x } \).
\par Note that since maps in the category \( \C \) are smooth functions,
\a group in that category is automatically \a \emph {Lie group}.
At this point, \I could define (internal) homomorphisms
to get the category \( \C ^ \Grp \) of internal groups.
But in the following, I'm only interested in \a single fixed group.
For the rest of part~\ref 1, let \( G \) be this fixed group.
\subsection {String diagrams for groups} \label {1string}
In part~\ref 2,
I'll make extensive use of string diagrams (\cite {string})
to perform calculations in \two groups.
But in fact, these can already be used in groups,
and (although they're not essential)
I'll make reference to them throughout the rest of part~\ref 1 for practice.
\par If \( X \maplongto x G \)
is \a generalised element of the group (or monoid) \( G \),
then this map is drawn as in the diagram on the left below.
Multiplication is shown by juxtaposition;
given elements \( x \) and \( y \),
\( \msub x y \) is in the middle below.
The identity element \( \esub \) is generally invisible,
but it may be shown as on the right below to stress its presence.
\begin {equation} \label {monoid maps string}
\matrixxy { \dstring x \\\\ & }
\qquad \qquad \qquad
\matrixxy { \dstring x & \dstring y \\\\ & }
\qquad \qquad \qquad
\matrixxy { & \estring \\\\ & }
\end {equation}
String diagrams reflect the associativity of monoids,
so there is no distinction between
\( \msub { \msub x y } z \) and \( \msub x { \msub y z } \),
both of which are shown on the left below.
Similarly, the picture on the right below
may be interpreted as \( x \),
\( \msub \esub x \), or \( \msub x \esub \);
\( \esub \) is properly invivisible here.
You can also consider these maps
to be pictures of the associative and unit laws;
these laws are invisible to string diagrams.
\begin {equation} \label {monoid string}
\matrixxy { \dstring x & \dstring y & \dstring z \\\\ & & }
\qquad \qquad \qquad
\matrixxy { & \dstring x \\\\ & }
\end {equation}
\par String diagrams are boring for monoids,
but the inverse operation is more interesting.
\I will denote \( \isub x \) simply by the diagram at the left,
but the inverse laws are now interesting:
\begin {equation} \label {inverse string}
\matrixxy { \ustring x \\\\ & }
\qquad \qquad \qquad
\matrixxy { & \estring & \\\\ \dstring x \rstring \iota && \ustring x \\\\ && }
\qquad \qquad \qquad
\matrixxy { \ustring x && \dstring x \\\\
\rstring \epsilon & \estring & \\\\ & }
\end {equation}
The point of such diagrams
is to say that the expression at the top of the diagram
is equal to the expression at the bottom,
for the reasons given in the middle.
(Of course, since these equations go both ways,
they can also appear upside down.)
\par \I can use string diagrams
to prove equations in elementary string theory.
For example, to prove
that \( \msub x y = \msub x z \) implies that \( y = z \),
you could use this string diagram:
\begin {equation} \label {stringexample}
\matrixxy { && && \ddstring y \\\\
\uustring x \rstring \iota && \ddstring x \\\\
&& \rstring \chi && \ddstring z \\\\
\rstring \iota && \\\\
&& && && }
\end {equation}
where \( \chi \) is the given fact \( \msub x y = \msub x z \).
This argument is essentially this:
\( y = \msub \esub y = \msub { \msub { \isub x } x } y
= \msub { \isub x } { \msub x y } = \msub { \isub x } { x z }
= \msub { \msub { \isub x } x } z = \msub \esub z = z \),
which uses (in turn) the left unit law~\( \lambda \),
the right inverse law~\( \iota \), the associative law~\( \alpha \),
the given fact \( \chi \), the associative law again,
the right inverse law again, and finally the left unit law again.
Notice that the associative and unit laws are invisible in the string diagram,
while the others appear explicitly;
this is typical.
\par The rest of part~\ref 1 will have several more examples,
with each string diagram also interpreted as an equation.
This is perhaps overkill;
but in part~\ref 2, these equations will become specific isomorphisms,
and the string diagrams will be quite helpful in keeping track of them.
\subsection {Action on \a space} \label {1act}
\I want \( G \) to act on the fibre of \a bundle.
The inverse map \( i \) will play no role in this section,
which could be applied when \( G \) is just \a monoid
(as was explicitly recognised at least as early as \cite [III.A.1] {Weyl}).
\par \A \strong {right \G space}
is \a space \( F \) equipped with
\a map \( F \times G \maplongto r F \)
that makes the following diagrams commute:
\begin {equation} \label {Gspace}
\matrixxy
{ F \times G \times G
\ar [dd] _ { \id F \times m } \ar [rrr] ^ { r \times \id G } &
\relax \ddricell \mu && F \times G \ar [dd] ^ r \\\\
F \times G \ar [rrr] _ r &&& F }
\qquad
\matrixxy
{ F \ar [ddrr] _ { \id F } \ar [rr] ^ { \id F \times e } &
\relax \dricell \upsilon &
F \times G \ar [dd] ^ r \\ && \\
&& F }
\end {equation}
If more than one \G space is around at \a time,
then I'll use subscripts or primes on \( r \) to keep things straight.
\par In terms of generalised elements,
if \( X \maplongto w F \) is \a generalised element of \( F \)
and \( X \maplongto x G \) is \a generalised element of \( G \),
then \I will write \( \rsub x z \)
for the composite \( X \maplongto { \pair w x } F \times G \maplongto r F \),
overloading the notation that introduced in section~\ref {1group} for \( m \).
Then if \( w \) is an element of \( F \)
and \( x \) and \( y \) are elements of \( G \),
then the laws (\ref {Gspace})
say (respectively) that \( \rsub { \rsub w x } y = \rsub w { \msub x y } \)
and that \( \rsub w \esub = w \).
\par The axioms \( \mu \) and \( \upsilon \)
are the associative and right unit laws of the right action.
These fit in with the laws for the group itself:
\begin {prop} \label {identity Gspace}
The group \( G \) is itself \a \G space.
\end {prop}
\begin {prf}
Set \( r _ G : = m \).
The requirements (\ref {Gspace})
then reduce to (some of) the requirements for \a group.
\end {prf}
I'll use this to define \G torsors and principal \G bundles.
\subsection {The category of \G spaces} \label {1C^G}
Just as spaces form \a category \( \C \),
so \G spaces form \a category \( \C ^ G \).
\par Given right \G spaces \( F \) and \( F ' \),
\a \strong {\G map} from \( F \) to \( F ' \)
is \a map \( F \maplongto t F ' \)
making the following diagram commute:
\begin {equation} \label {Gmap}
\matrixxy
{ F \times G \ar [dd] _ r \ar [rrr] ^ { t \times \id G } &
\relax \ddricell \phi && F ' \times G \ar [dd] ^ { r ' } \\ \\
F \ar [rrr] _ t &&& F ' }
\end {equation}
Note that \a composition of \G maps is \a \G map,
since both squares below commute:
\begin {equation} \label {composite Gmap}
\matrixxy
{ F \times G \ar [dd] _ r \ar [rrr] ^ { t \times \id G } &
\relax \ddricell \phi &&
F ' \times G \ar [dd] _ { r ' } \ar [rrr] ^ { t ' \times \id G } &
\relax \ddricell { \phi ' } & & F ' ' \times G \ar [dd] ^ { r ' ' } \\ \\
F \ar [rrr] _ t &&& F ' \ar [rrr] _ { t ' } & & & F ' ' }
\end {equation}
The identity map on \( F \) is also \a \G map.
Indeed, \G spaces and \G maps form \a category \( \C ^ G \).
\par Now \I know what it means
for \G spaces \( F \) and \( F ' \) to be \strong {equivalent}:
There must be maps
\( F \maplongto t F ' \) and \( F ' \maplongto { \bar t } F \)
such that the following diagrams commute:
\begin {eqnarray} \label {Gspace equivalence}
&
\matrixxy
{ F \times G \ar [dd] _ r \ar [rrr] ^ { t \times \id G } &
\relax \ddricell \phi && F ' \times G \ar [dd] ^ { r ' } \\ \\
F \ar [rrr] _ t &&& F ' }
\qquad
\matrixxy
{ F ' \times G \ar [dd] _ { r ' } \ar [rrr] ^ { \bar t \times \id G } &
\relax \ddricell { \bar \phi } && F \times G \ar [dd] ^ r \\ \\
F ' \ar [rrr] _ { \bar t } &&& F }
& \nonumber \\ \\ &
\matrixxy
{ F \ar [drr] _ t \ar [rrrr] ^ { \id F } & \relax \drricell \tau & && F \\
&& F ' \ar [urr] _ { \bar t } & }
\qquad
\matrixxy
{ F ' \ar [drr] _ { \bar t } \ar [rrrr] ^ { \id { F ' } } &
\relax \drricell { \bar \tau } & && F ' \\
&& F \ar [urr] _ t & }
& \nonumber
\end {eqnarray}
\subsection {\G torsors} \label {1torsor}
When \I define principal bundles,
\I will want the fibre of the bundle to be \( G \).
However, it should be good enough
if the fibre is only \emph {equivalent} to \( G \).
Such \a \G space is sometimes called \q {homogeneous},
but \I prefer the noun \q {torsor}.
Thus, \a \strong {right \G torsor} is any right \G space
that is equivalent, as \a \G space, to \( G \) itself.
\par In more detail, this is \a space \( F \)
equipped with maps \( F \times G \maplongto r F \),
\( F \maplongto t G \), and \( G \maplongto { \bar t } F \),
such that the following diagrams commute:
\begin {eqnarray} \label {Gtorsor}
&
\matrixxy
{ F \times G \times G
\ar [dd] _ { \id F \times m } \ar [rrr] ^ { r \times \id G } &
\relax \ddricell \mu && F \times G \ar [dd] ^ r \\\\
F \times G \ar [rrr] _ r &&& F }
\qquad
\matrixxy
{ F \ar [ddrr] _ { \id F } \ar [rr] ^ { \id F \times e } &
\relax \dricell \upsilon &
F \times G \ar [dd] ^ r \\ && \\
&& F }
\qquad
\matrixxy
{ F \times G \ar [dd] _ r \ar [rrr] ^ { t \times \id G } &
\relax \ddricell \phi && G \times G \ar [dd] ^ m \\ \\
F \ar [rrr] _ t &&& G }
& \nonumber \\ \\ &
\matrixxy
{ G \times G \ar [dd] _ m \ar [rrr] ^ { \bar t \times \id G } &
\relax \ddricell { \bar \phi } && F \times G \ar [dd] ^ r \\ \\
G \ar [rrr] _ { \bar t } &&& F }
\qquad
\matrixxy
{ F \ar [drr] _ t \ar [rrrr] ^ { \id F } & \relax \drricell \tau & && F \\
&& G \ar [urr] _ { \bar t } & }
\qquad
\matrixxy
{ G \ar [drr] _ { \bar t } \ar [rrrr] ^ { \id G } &
\relax \drricell { \bar \tau } & && G \\
&& F \ar [urr] _ t & }
& \nonumber
\end {eqnarray}
\section {\G bundles} \label {1Gbun}
\I can now put the above ideas together to get the concept of \G bundle.
\subsection {Definition of \G bundle} \label {1def}
Suppose that \I have \a space \( B \) with \a cover \( U \maplongto j B \),
as well as \a group \( G \).
In these circumstances,
\a \strong {\G transition} of the cover \( U \)
is \a map \( \nerve U 2 \maplongto g G \)
such that the following diagrams commute:
\begin {equation} \label {transition}
\matrixxy
{ \relax \nerve U 3
\ar [ddrrrr] _ { \del j { 0 2 } }
\ar [rrrr] ^ { \pair { \del j { 0 1 } } { \del j { 1 2 } } } &&&&
\relax \nerve U 2 \times \nerve U 2 \ar [rrr] ^ { g \times g }
\relax \ddricell \gamma &&&
G \times G \ar [dd] ^ m \\\\
&&&& \relax \nerve U 2 \ar [rrr] _ g &&& G }
\qquad
\matrixxy
{ U \ar [dd] _ { \del j { 0 0 } } \ar [rr] ^ { \hat U } \ddrricell \eta &&
1 \ar [dd] ^ e \\ \\
\relax \nerve U 2 \ar [rr] _ g & & G }
\end {equation}
Notice that if \( U \) appears directly as \a disjoint union,
then \( g \) breaks down into \a family of maps on each disjunct;
thus one normally speaks of \emph {transition maps} in the plural.
\par In terms of elements,
if \( \pair x y \) defines an element of \( \nerve U 2 \),
then let \( \gsub x y \) be the composite
\( X \maplongto { \pair x y } \nerve U 2 \maplongto g G \).
Then the law \( \gamma \)
says that \( \msub { \gsub x y } { \gsub y z } = \gsub x z \)
(for an element \( \triple x y z \) of \( \nerve U 3 \)),
while the law \( \eta \) says that \( \esub = \gsub x x \)
(for an element \( X \maplongto x U \) of \( U \)).
The string diagrams for \( \gsub x y \),
for \( \gammasub x y z \), and for \( \etasub x \)
are (respectively):
\begin {equation} \label {transition string}
\matrixxy { & \dstring g \\ x & & y \\ & }
\qquad \qquad
\matrixxy
{ \dstring g && && \dstring g \\ && y \\
\rrstring \gamma && \dstring g && \\ x && && z \\ && }
\qquad \qquad
\matrixxy { & \estring \\ x \\ & \pstring \eta \dstring g \\\\ & }
\end {equation}
(Here I've drawn in the invisible identity string diagram
to clarify the height of the diagram for \( \eta \).)
\I should also point out that the string diagrams of this form
cover only some of the meaningful expressions that one can build using \( g \);
for example, there is no diagram for \( \msub { \gsub w x } { \gsub y z } \).
However, \I will never need to use this expression;
the string diagrams can handle those expressions that \I want.
\par As you will see, there is an analogy to be made
between \a \G transition and \a group.
In this analogy, \( g \) corresponds to the underlying set of the group,
the commutative diagram \( \gamma \)
corresponds to the operation of multiplication,
and the commutative diagram \( \eta \) corresponds to the identity element.
(This analogy will be even more fruitful when categorified.)
Another possible analogy is with \a group \emph {homomorphism};
in fact, the following proposition (and its proof)
are analogous to the result that \a function between groups,
if it preserves multiplication, must also preserve the identity.
(The \emph {truth} is that \a \G transition
is \a \emph {groupoid} homomorphism;
see the discussion in section~\ref {disc2space},
where groupoid homomorphisms appear in the guise of \two maps.)
\begin {prop} \label {gamma implies eta}
If \a map \( \nerve U 2 \maplongto g G \)
satisfies the law \( \gamma \), then it also satisfies \( \eta \).
\end {prop}
\begin {prf}
\( \eta \) may be given in terms of generalised elements:
\( \esub = \msub { \gsub x x } { \isub { \gsub x x } }
= \msub { \msub { \gsub x x } { \gsub x x } } { \isub { \gsub x x } }
= \msub { \gsub x x } { \msub { \gsub x x } { \isub { \gsub x x } } }
= \msub { \gsub x x } \esub = \gsub x x \);
using (in order) the inverse law \( \iota \),
the assumed law \( \gamma \), the associative law \( \alpha \),
\( \iota \) again, and the unit law \( \rho \).
\end {prf}
\I think that these proofs are easier to read in terms of elements,
so I'll normally give them that way, without drawing \a commutative diagram.
But it's important to know that such \a diagram can be drawn;
the reasoning does not depend on
\( \C \)'s being \a concrete category with literal points as elements.
This proof can also be written as \a string diagram:
\begin {equation} \label {gamma implies eta string}
\matrixxy
{ && \dstring g \rrstring \iota && && \uustring g \\
&& && & x \\ \ddstring g \rrstring \gamma && && \dstring g \\\\
&& x && \rstring \iota && \\\\ & }
\end {equation}
\par If \( G \) acts on some (right) \G space \( F \),
and if \( E \) is \a locally trivial bundle over \( B \) with fibre \( F \),
then \I can form the following diagram:
\begin {equation} \label {Gbundle}
\matrixxy
{ F \times \nerve U 2
\ar [ddrrrr] _ { \id F \times \del j 0 }
\ar [rrrr] ^ { \id F \times \pair g { \del j 1 } } &&&&
F \times G \times U \ar [rrr] ^ { r \times \id U } \ddricell \theta &&&
F \times U \ar [dd] ^ { \tilde j } \\\\
&&&& F \times U \ar [rrr] _ { \tilde j } &&& E }
\end {equation}
By definition, \I have \a \strong {\G bundle} if this diagram commutes.
\I say that the \G bundle \( E \)
is \strong {associated} with the \G transition \( g \).
This isn't quite the normal usage of the term \q {associated};
see section~\ref {1assocbun} for the connection.
\par In terms of elements,
if \( X \maplongto w F \) is an element of \( F \)
and \( \pair x y \) defines an element of \( \nerve U 2 \),
then this law states
that \( \jtildesub { \rsub w { \gsub x y } } y = \jtildesub w x \).
As \a string diagram:
\begin {equation} \label {Gbundle string}
\matrixxy
{ \ddStringl w && \dstring g && \dString \\ && & y \\
&& \rString && \theta \dString \\ & x \\ && && }
\end {equation}
The bottom part of this diagram represents \( \jtildesub w x \),
\a general notation that covers all meanings of \( \tilde j \);
while the top part,
which represents \( \jtildesub { \rsub w { \gsub x y } } y \),
cannot be generalised to, say,
\( \jtildesub { \rsub w { \gsub x y } } z \),
even though that expression makes sense for any element \( z \) of \( \UU \).
(But \I will have no need to refer to such expressions.)
\par The law \( \theta \) has certain analogies with \a group action;
these will show up in part~\ref 2.
\par \A \strong {principal \G bundle}
is simply \a \G bundle whose fibre \( F \) is \( G \) itself.
\subsection {The category of \G transitions} \label {1G^B}
To define \a morphism of \G bundles properly,
\I first need the notion of morphism of \G transitions.
\par Given covers \( U \) and \( U ' \) of \( B \),
\a \G transition \( g \) on \( U \),
and \a \G transition \( g ' \) on \( U ' \),
\a \strong {\G transition morphism} from \( g \) to \( g ' \)
is \a map \( U \cap U ' \maplongto b G \)
such that these diagrams commute:
\begin {equation} \label {transition morphism left}
\matrixxy
{ \relax \nerve U 2 \cap U '
\ar [ddrrrr] _ { \del j { 0 2 } }
\ar [rrrr] ^ { \pair { \del j { 0 1 } } { \del j { 1 2 } } } &&&&
\relax \nerve U 2 \times U \cap U ' \ar [rrrr] ^ { g \times b }
\relax \ddricell \sigma &&&&
G \times G \ar [dd] ^ m \\\\
&&&& U \cap U ' \ar [rrrr] _ b &&&& G }
\end {equation}
and:
\begin {equation} \label {transition morphism right}
\matrixxy
{ U \cap \nerve { U ' } 2
\ar [ddrrrr] _ { \del j { 0 2 } }
\ar [rrrr] ^ { \pair { \del j { 0 1 } } { \del j { 1 2 } } } &&&&
U \cap U ' \times \nerve { U ' } 2 \ar [rrrr] ^ { b \times g ' }
\relax \ddricell \delta &&&&
G \times G \ar [dd] ^ m \\\\
&&&& U \cap U ' \ar [rrrr] _ b &&&& G }
\end {equation}
In terms of generalised elements,
if \( \pair x { x ' } \) defines an element of \( U \cap U ' \),
then let \( \bsub x { x ' } \)
be \( X \maplongto { \pair x { x ' } } U \cap U ' \maplongto b G \).
Then if \( \triple x y { x ' } \)
defines an element of \( \nerve U 2 \cap U ' \),
then the law \( \sigma \)
says that \( \msub { \gsub x y } { \bsub y { x ' } } = \bsub x { x ' } \);
and if \( \triple x { x ' } { y ' } \)
defines an element of \( U \cap \nerve { U ' } 2 \),
then the law \( \delta \)
says that \( \msub { \bsub x { x ' } } { \gpsub { x ' } { y ' } }
= \bsub x { y ' } \).
Here are the string diagrams for \( \bsub x { x ' } \),
for \( \sigmasub x y { x ' } \), and for \( \deltasub x { x ' } { y ' } \):
\begin {equation} \label {transition morphism string}
\matrixxy { & \dstring b \\ x & & x ' \\ & }
\qquad \qquad
\matrixxy
{ \dstring g && && \dstring b \\ && y \\
\rrstring \sigma && \dstring b && \\ x && && x ' \\ && }
\qquad \qquad
\matrixxy
{ \dstring b && && \dstring { g ' } \\ && x ' \\
\rrstring \delta && \dstring b && \\ x && && y ' \\ && }
\end {equation}
\par In the analogy
where the \G transitions \( g , g ' \) are analogous to groups,
the transition morphism \( b \)
is analogous to \a space
that is acted on by \( g \) on the left and by \( g ' \) on the right.
(This is sometimes called
\a \emph {bihomogeneous space} or \a \emph {bitorsor},
hence the letter \q {\( b \)} in my notation.)
If \( \gamma \) and \( \eta \) in~(\ref {transition})
are analogous to the multiplication and identity in the groups,
then \( \sigma \) and \( \delta \) here
are analogous (respectively) to the left and right actions.
(The notation recalls the Latin words \q {sinister} and \q {dexter}.)
\par As \a simple example, take \( U ' \) be the same cover as \( U \),
and take \( g \), \( g ' \), and \( b \) to be all the same map \( g \);
then diagrams (\ref {transition morphism left})
and~(\ref {transition morphism right})
both reduce to the first diagram in~(\ref {transition}).
In this way,
the \G transition \( g \) serves as its own identity \G transition morphism.
\begin {prop} \label {Gtransition morphism composite}
Suppose that \( g , g ' , g ' ' \) are all \G transitions
on (respectively) covers \( U , U ' , U ' ' \),
and suppose that \( b \) is \a \G transition morphism from \( g \) to \( g ' \)
while \( b ' \) is \a \G transition morphism from \( g ' \) to \( g ' ' \).
Then \I can define \a \G transition morphism
\( b ; b ' \) from \( g \) to \( g ' ' \)
such that this diagram commutes:
\begin {equation} \label {composite transition morphism universal property}
\matrixxy
{ U \cap U ' \cap U ' '
\ar [ddrrrrr] _ { \del j { 0 2 } }
\ar [rrrrr] ^ { \pair { \del j { 0 1 } } { \del j { 1 2 } } } &&&&&
U \cap U ' \times U ' \cap U ' '
\ar [rrrr] ^ { b \times b ' } \ddricell \beta &&&&
G \times G \ar [dd] ^ m \\\\
&&&&& U \cap U ' ' \ar [rrrr] _ { b ; b ' } &&&& G }
\end {equation}
\end {prop}
Given an element \( \pair x { x ' ' } \) of \( U \cap U ' ' \),
\I will denote its composite with \( b ; b ' \) as \( \bcbpsub x { x ' ' } \).
Then given an element \( \triple x { x ' } { x ' ' } \)
of \( U \cap U ' \cap U ' ' \),
the property (\ref {composite transition morphism universal property})
of \( b ; b ' \)
states that \( \msub { \bsub x { x ' } } { \bpsub { x ' } { x ' ' } }
= \bcbpsub x { x ' ' } \).
\begin {prf}
This is probably most clearly expressed
in terms of an element \( \quadruple x { x ' } { y ' } { y ' ' } \)
of \( U \cap \nerve { U ' } 2 \cap U ' ' \).
Applying the laws \( \sigma ' \), \( \alpha \), and \( \delta \) in turn,
\I find that \( \msub { \bsub x { x ' } } { \bpsub { x ' } { y ' ' } }
= \msub { \bsub x { x ' } }
{ \msub { \gpsub { x ' } { y ' } } { \bpsub { y ' } { y ' ' } } }
= \msub { \msub { \bsub x { x ' } } { \gpsub { x ' } { y ' } } }
{ \bpsub { y ' } { y ' ' } }
= \msub { \bsub x { y ' } } { \bpsub { y ' } { y ' ' } } \).
As \a string diagram:
\begin {equation} \label {pulledback composite transition morphism string}
\matrixxy
{ & \ddstring b && && && \dstring { b ' } \\ & && && && && & y ' ' \\
& && x ' && \dstring { g ' } \rrstring { \sigma ' } && &&
\ddstring { b ' } \\\\
& \rrstring \delta && \dstring b && && y ' \\ x \\ & && && && && }
\end {equation}
Thus, this diagram commutes:
\begin {equation} \label {pulledback composite transition morphism}
\matrixxy
{ U \cap \nerve { U ' } 2 \cap U ' '
\ar [dddddd] _ { \del j { 0 2 3 } }
\ar [rrrrrrrrrrr] ^ { \del j { 0 1 3 } } &&&&& &&&& &&
U \cap U ' \cap U ' '
\ar [dd] ^ { \pair { \del j { 0 1 } } { \del j { 1 2 } } } \\\\
&&&&& &&&& && U \cap U ' \times U ' \cap U ' ' \ar [dd] ^ { b \times b ' } \\\\
&&&&& &&&& && G \times G \ar [dd] ^ m \\\\
U \cap U ' \cap U ' '
\ar [rrrrr] _ { \pair { \del j { 0 1 } } { \del j { 1 2 } } } &&&&&
U \cap U ' \times U ' \cap U ' ' \ar [rrrr] _ { b \times b ' } &&&&
G \times G \ar [rr] _ m && G }
\end {equation}
This diagram is an example of~(\ref {quotient cone})
for the cover \( U \cap U ' \cap U ' '
\maplongto { \del j { 0 2 } } U \cap U ' ' \),
so (\ref {quotient universality}) defines
\a map \( U \cap U ' ' \maplongto { b ; b ' } G \)
satisfying~(\ref {composite transition morphism universal property}).
\par Now given an element \( \quadruple x y { y ' } { y ' ' } \)
of \( \nerve U 2 \cap U ' \cap U ' ' \),
\I use (in turn) the laws
\( \beta \), \( \alpha \), \( \sigma \), and \( \beta \) again
to deduce that \( \msub { \gsub x y } { \bcbpsub y { y ' ' } }
= \msub { \gsub x y }
{ \msub { \bsub y { y ' } } { \bpsub { y ' } { y ' ' } } }
= \msub { \msub { \gsub x y } { \bsub y { y ' } } }
{ \bpsub { y ' } { y ' ' } }
= \msub { \bsub x { y ' } } { \bpsub { y ' } { y ' ' } }
= \bcbpsub x { y ' ' } \).
As \a string diagram:
\begin {equation}
\label {composite transition morphism left action homomorph string}
\matrixxy
{ \ddstring g && && && \dstring { b ; b ' } \\ && & y \\
&& && \dstring b \rrstring \beta && && \ddstring { b ' } \\\\
\rrstring \sigma && \dstring b && && y ' \\ && &&& &&& & y ' ' \\
x && \rrrstring \beta &&& \dstring { b ; b ' } &&& \\\\ && &&& }
\end {equation}
Thus, this diagram commutes:
\begin {equation}
\label {composite transition morphism left action homomorph}
\matrixxy
{ \relax \nerve U 2 \cap U ' \cap U ' '
\ar [dddddd] _ { \del j { 0 1 3 } } \ar [rrrrr] ^ { \del j { 0 1 3 } } &&& &&
\relax \nerve U 2 \cap U ' '
\ar [dd] ^ { \pair { \del j { 0 1 } } { \del j { 1 2 } } } \\\\
&&& && \relax \nerve U 2 \times U \cap U ' '
\ar [dd] ^ { g \times b ; b ' } \\\\
&&& && G \times G \ar [dd] ^ m \\\\
\relax \nerve U 2 \cap U ' ' \ar [rrr] _ { \del j { 0 2 } } &&&
U \cap U ' ' \ar [rr] _ { b ; b ' } && G }
\end {equation}
Because the cover \( \nerve U 2 \cap U ' \cap U ' '
\maplongto { \del j { 0 1 3 } } U \cap U ' \cap U ' ' \)
is an epimorphism, it follows that this diagram commutes:
\begin {equation} \label {composite transition morphism left action}
\matrixxy
{ \relax \nerve U 2 \cap U ' '
\ar [ddrrrr] _ { \del j { 0 2 } }
\ar [rrrr] ^ { \pair { \del j { 0 1 } } { \del j { 1 2 } } } &&&&
\relax \nerve U 2 \times U \cap U ' ' \ar [rrrr] ^ { g ' \times b ; b ' } &&&&
G \times G \ar [dd] ^ m \\\\
&&&& U \cap U ' ' \ar [rrrr] _ { b ; b ' } &&&& G }
\end {equation}
This is precisely the law \( \sigma \) for \( b ; b ' \).
(In terms of an element \( \triple x y { y ' ' } \)
of \( \nerve U 2 \cap U ' ' \),
this is the same equation
\( \msub { \gsub x y } { \bcbpsub y { y ' ' } } = \bcbpsub x { y ' ' } \)
as above;
it is because \( \del j { 0 1 3 } \) is an epimorphism
that \I don't need to posit \( y ' \) anymore.
In more set-theoretic terms, since \( \del j { 0 1 3 } \) is onto,
some \( y ' \) must always exist.)
\par The law \( \delta \) for \( b ; b ' \) follows by an analogous argument;
thus, \( b ; b ' \) is \a \G transition morphism, as desired.
\end {prf}
This \( b ; b ' \)
is the composite of \( b \) and \( b ' \)
in the category \( G ^ B \) of \G transitions in \( B \).
\begin {prop} \label {Gtransition cat}
\( G ^ B \) is \a category.
\end {prop}
\begin {prf}
\I must check that composition of \G transitions
is associative and has identities.
The key to this
is that the diagram (\ref {pulledback composite transition morphism})
not only shows (as in the previous proof)
that that \a map \( b ; b ' \)
satisfying~(\ref {composite transition morphism universal property})
exists,
but that \( b ; b ' \)
is the \emph {unique} map
satisfying~(\ref {composite transition morphism universal property}).
Associativity follows.
\end {prf}
\I now know what it means
to say that the \G transitions \( g \) and \( g ' \) are equivalent:
there are \G transition morphisms
\( b \) from \( g \) to \( g ' \)
and \( \bar b \) from \( g ' \) to \( g \),
such that \( b ; \bar b \) equals \( g \)
(which is the identity bundle morphism on \( g \))
and \( \bar b ; b \) equals \( g ' \) (the identity on~\( g ' \)).
\subsection {The category of \G bundles} \label {1Bun}
To classify \G \two bundles,
\I need \a proper notion of equivalence of \G bundles.
For this, \I should define the category
\( \Bun \C G B F \) of \G bundles over \( B \) with fibre \( F \).
\par Assume \G bundles \( E \) and \( E ' \) over \( B \),
both with the given fibre \( F \),
and associated with the \G transitions \( g \) and \( g ' \) (respectively).
Then \a \strong {\G bundle morphism} from \( E \) to \( E ' \)
is \a bundle morphism \( f \) from \( E \) to \( E ' \)
together with \a \G transition morphism \( b \) from \( g \) to \( g ' \),
such that this diagram commutes:
\begin {equation} \label {Gbundle morphism}
\matrixxy
{ F \times U \cap U '
\ar [dd] _ { \id F \times \del j 0 }
\ar [rrrr] ^ { \id F \times \pair b { \del j 1 } } &&&&
F \times G \times U ' \ar [dd] ^ { r \times \id { U ' } } \\
& \relax \ddrricell \zeta \\
F \times U \ar [ddrr] _ { \tilde j } && &&
F \times U ' \ar [dd] ^ { \tilde j } \\
&& & \\ && E \ar [rr] _ f && E ' }
\end {equation}
In terms of generalised elements,
if \( w \) is an element of \( F \)
and \( \pair x { x ' } \) defines an element of \( U \cap U ' \),
then \( \jptildesub { \rsub w { \bsub x { x ' } } } { x ' } \)
is equal to the composite of \( \jtildesub w x \) with \( f \);
this also can be written as \a string diagram:
\begin {equation} \label {Gbundle morphism string}
\matrixxy { \ddStringl w && \dstring b && \dString \\ && & x ' \\
&& \rString && \zeta \dString \\ & x \\ && && }
\end {equation}
(In this string diagram,
the transition labelled \( \zeta \)
does not ---like usual--- indicate equality of the top and bottom sides,
but instead that the top is the composite of the bottom with \( f \).
In other words,
\( f \) is applied to the portion in \( E \)
so that the whole diagram may be interpreted in \( E ' \).)
\I think of \( f \) as \emph {being} the \G bundle morphism
and say that \( f \)
is \strong {associated} with the transition morphism \( b \).
\par If \( E \) and \( E ' \) are the same \G bundle
(so associated with the same \G transition map)
and \( f \) is the identity map \( \id E \),
then these diagrams (\ref {Gbundle morphism}, \ref {Gbundle morphism string})
are \a special case of the diagrams (\ref {Gbundle}, \ref {Gbundle string}),
with \( b \) taken to be the identity \G transition morphism on \( g \),
which is \( g \) itself.
In this way, every \G bundle has an identity \G bundle automorphism.
\par Given \G bundle morphisms
\( f \) from \( E \) to \( E ' \) and \( f ' \) from \( E ' \) to \( E ' ' \),
the composite \G bundle morphism \( f ; f ' \)
is simply the composite bundle morphism
\( E \maplongto f E ' \maplongto { f ' } E ' ' \)
together with the composite \G transition \( b ; b ' \)
as described in section~\ref {1G^B}.
\begin {prop} \label {Gbundle morphism composite}
This \( f ; f ' \) really is \a \G bundle morphism.
\end {prop}
\begin {prf}
Let \( w \) be an element of \( F \),
and let \( \triple x { x ' } { x ' ' } \)
define an element of \( U \cap U ' \cap U ' ' \).
Then using (in turn)
the laws \( \zeta \), \( \zeta ' \), \( \mu \), and \( \beta \),
\I see that \( \jpptildesub { \rsub w { \bcbpsub { x ' ' } x } } { x ' ' }
= \jpptildesub { \rsub w
{ \msub { \bsub x { x ' } } { \bpsub { x ' } { x ' ' } } } }
{ x ' ' }
= \jpptildesub { \rsub { \rsub w { \bsub x { x ' } } }
{ \bpsub { x ' } { x ' ' } } }
{ x ' ' } \),
which is the composite of \( f ' \)
with \( \jptildesub { \rsub w { \bsub x { x ' } } } { x ' } \),
which is the composite of \( f \) with \( \jtildesub w x \);
or as \a string diagram:
\begin {equation} \label {Gbundle morphism composite string}
\matrixxy
{ \ddddStringl w && && \dstring { b ; b ' } && && \ddString \\
&& && && & x ' ' \\
&& \ddstring b \rrstring \beta && && \dstring { b ' } && \\\\
&& && x ' && \rString && \zeta ' \dString \\\\
&& \rrrString && && && \zeta \dString \\ & x \\ && && && && }
\end {equation}
Because \( \nerve U 3 \maplongto { \del j { 0 2 } } \nerve U 2 \) is epic,
\I can ignore \( x ' \)
(just as \I ignored \( y ' \)
in the proof of Proposition~\ref {Gtransition morphism composite});
thus this diagram commutes:
\begin {equation} \label {composite Gbundle morphism}
\matrixxy
{ F \times U \cap U ' '
\ar [dd] _ { \id F \times \del j 0 }
\ar [rrrrrr] ^ { \id F \times \pair { b ; b ' } { \del j 1 } } &&&&&&
F \times G \times U ' ' \ar [dd] ^ { r \times \id { U ' ' } } \\
&& \relax \ddrricell { \zeta ; \zeta ' } \\
F \times U \ar [ddrr] _ { \tilde j } && && &&
F \times U ' ' \ar [dd] ^ { \tilde j ' ' } \\
&& && \\ && E \ar [rr] _ f && E ' \ar [rr] _ { f ' } && E ' ' }
\end {equation}
This is simply the diagram (\ref {Gbundle morphism})
for the \G bundle morphism \( f ; f ' \).
Therefore, \( f ; f ' \) really is \a \G bundle morphism.
\end {prf}
\begin {prop} \label {Gbundle cat}
Given \a space \( B \),
bundles over \( B \) and their bundle morphisms form \a category.
\end {prop}
\begin {prf}
Composition of bundle morphisms is just composition of maps,
which \I know to be associative;
and Proposition~\ref {Gtransition cat}
proves that composition of the associated \G transition morphisms
is associative;
thus, composition of \G bundle morphisms is associative.
Similarly, this composition has identities;
therefore, \G bundles form \a category \( \Bun \C G B F \).
\end {prf}
\par \I now know what it means
for \G bundles \( E \) and \( E ' \) to be \strong {equivalent \G bundles}:
isomorphic objects in the category \( \Bun \C G B F \).
There must be \a \G bundle morphism from \( E \) to \( E ' \)
and \a \G bundle morphism from \( E ' \) to \( E \)
whose composite \G bundle morphisms, in either order,
are identity \G bundle morphisms.
In particular, \( E \) and \( E ' \) are equivalent as bundles.
\par \I should make \a remark
about when the action of \( G \) on \( F \) is unfaithful.
Many references will define \a bundle morphism
to be \a bundle morphism \( f \)
\emph {such that there exists} \a transition morphism \( b \),
while I've defined it to be \( f \) \emph {equipped with} \( b \).
The difference affects the notion of equality of bundle morphisms:
whether \( f = f ' \) as maps is enough,
or if \( b = b ' \) (as maps) must also hold.
This is relevant to the next section~\ref {1assocbun} as well;
the propositions there are fine as far as they go,
but the final theorem requires \a functor from \( \Bun \C G B F \)
to the category \( B ^ G \) of \G transitions,
so each \G bundle morphism
must be associated with \a unique \G transition morphism.
You can fix this either by passing from \( G \) to \( \fract G N \),
where \( N \) is the kernel of the action of \( G \) on \( F \);
or even by requiring this action to be faithful
in the definition of \G bundle
(as is done, for example, in \cite [2.3] {Steenrod}).
\I find the theory cleaner without these restrictions;
in any case, there is no problem for principal bundles,
since \( G \) acts on itself faithfully.
\subsection {Associated bundles} \label {1assocbun}
The local data given in terms of \G transitions
is in fact sufficient to recreate the associated bundle.
\begin {prop} \label {associated bundle}
Given \a cover \( U \maplongto j B \),
\a \G transition \( \nerve U 2 \maplongto g G \),
and \a \G space \( F \),
there is \a \G bundle \( E \) over \( B \) with fibre \( F \)
associated with the transition \( g \).
\end {prop}
\begin {prf}
\I will construct \( E \) as the quotient of an equivalence relation
from \( F \times \nerve U 2 \) to \( F \times U \).
One of the maps in the equivalence relation is \( \id F \times \del j 0 \);
the other
is \( F \times \nerve U 2
\maplongto { \id F \times \pair g { \del j 1 } } F \times G \times U
\maplongto { r \times \id U } F \times U \).
(You can see these maps in diagram~(\ref {Gbundle}),
where they appear with the quotient map
\( F \times U \maplongto { \tilde j } E \)
that this proof will construct.)
Since \( \id F \times \del j 0 \) is \a cover
and every equivalence relation involving \a cover has \a quotient,
the desired quotient does exist, satisfying (\ref {Gbundle})---
at least, if this really is an equivalence relation!
\par To begin with, it is \a relation; that is, the two maps are jointly monic.
Given two elements
\( \pair w { \pair x y } \) and \( \pair { w ' } { \pair { x ' } { y ' } } \)
of \( F \times \nerve U 2 \),
the monicity diagrams (\ref {relation})
say that \( \pair w x = \pair { w ' } { x ' } \)
and \( \pair { \rsub w { \gsub x y } } y
= \pair { \rsub { w ' } { \gsub { x ' } { y ' } } } { y ' } \);
if these hold, then certainly \( w = w ' \), \( x = x ' \), and \( y = y ' \),
which is what \I need.
\par The reflexivity map of the equivalence relation
is \( F \times U
\maplongto { \id F \times \del j { 0 0 } } F \times \nerve U 2 \);
this is obviously \a section of \( \id F \times \del j 0 \).
It's \a section of \( F \times \nerve U 2
\maplongto { \id F \times \pair g { \del j 1 } } F \times G \times U
\maplongto { r \times \id U } F \times U \)
using the laws \( \upsilon \) and \( \eta \);
in terms of elements,
\( \pair w x = \pair { \rsub w \esub } x
= \pair { \rsub w { \gsub x x } } x \)
for an element \( \pair w x \) of \( F \times U \);
or in \a string diagram:
\begin {equation} \label {bundle construction eqrel reflexivity}
\matrixxy { \ddStringl w \\ && & x \\ && \pstring \eta \dstring g \\\\ && }
\end {equation}
\par The kernel pair of \( \id F \times \del j 0 \)
is \( F \times \nerve U 3 \):
\begin {equation} \label {bundle construction eqrel kernel pair}
\matrixxy
{ F \times \nerve U 3
\ar [dd] _ { \id F \times \del j { 0 1 } }
\ar [rrr] ^ { \id F \times \del j { 0 2 } } &&&
F \times \nerve U 2 \ar [dd] ^ { \id F \times \del j 0 } \\\\
F \times \nerve U 2 \ar [rrr] _ { \id F \times \del j 0 } &&&
F \times U }
\end {equation}
Let the Euclideanness map
be \( F \times \nerve U 3
\maplongto { \id F \times \pair { \del j { 0 1 } } { \del j { 1 2 } } }
F \times \nerve U 2 \times \nerve U 2
\maplongto { \id F \times g \times \id { \nerve U 2 } }
F \times G \times \nerve U 2
\maplongto { r \times \id { \nerve U 2 } } F \times \nerve U 2 \).
Thus, to prove that \I have an equivalence relation,
\I need to show that these diagrams commute:
\begin {equation} \label {bundle construction eqrel royal}
\matrixxy
{ F \times \nerve U 3
\ar [dd] _ { \id F \times \del j { 0 1 } }
\ar [rrrrr] ^
{ \id F \times \pair { \del j { 0 1 } } { \del j { 1 2 } } } &&&&&
F \times \nerve U 2 \times \nerve U 2
\ar [rrrrr] ^ { \id F \times g \times \id { \nerve U 2 } } &&&&&
F \times G \times \nerve U 2 \ar [rrrr] ^ { r \times \id { \nerve U 2 } } &&&&
F \times \nerve U 2 \ar [dd] ^ { \id F \times \del j 0 } \\\\
F \times \nerve U 2
\ar [rrrrrrrrrr] _ { \id F \times \pair g { \del j 1 } } &&&&& &&&&&
F \times G \times U \ar [rrrr] _ { r \times \id U } &&&& F \times U }
\end {equation}
and:
\begin {equation} \label {bundle construction eqrel Euclidean}
\matrixxy
{ F \times \nerve U 3
\ar [dd] _ { \id F \times \pair { \del j { 0 1 } } { \del j { 1 2 } } }
\ar [rrrrrrr] ^ { \id F \times \del j { 0 2 } } &&&& &&&
F \times \nerve U 2 \ar [dddd] ^ { \id F \times \pair g { \del j 1 } } \\\\
F \times \nerve U 2 \times \nerve U 2
\ar [dd] _ { \id F \times g \times \id { \nerve U 2 } } \\\\
F \times G \times \nerve U 2
\ar [dd] _ { r \times \id { \nerve U 2 } } &&&& &&&
F \times G \times U \ar [dd] ^ { r \times \id U } \\\\
F \times \nerve U 2 \ar [rrrr] _ { \id F \times \pair g { \del j 1 } } &&&&
F \times G \times U \ar [rrr] _ { r \times \id U } &&& F \times U }
\end {equation}
The first of these diagrams
is essentially one of the structural diagrams of the cover \( U \)
(all that business with \( g \) and \( r \) is just along for the ride);
in terms of elements,
both sides describe \( \pair { \rsub w { \gsub x y } } y \)
for an element \( \pair w { \triple x y z } \) of \( F \times \nerve U 3 \).
The second diagram holds using the laws \( \mu \) and \( \gamma \);
in terms of elements,
\( \pair { \rsub w { \gsub x z } } z
= \pair { \rsub w { \msub { \gsub x y } { \gsub y z } } } z
= \pair { \rsub { \rsub w { \gsub x y } } { \gsub y z } } z \);
in \a string diagram:
\begin {equation} \label {bundle construction eqrel Euclidean string}
\matrixxy
{ \ddStringl w && && \dstring g \\ & x & && && & z \\
&& \dstring g \rrstring \gamma && && \dstring g \\ && && y \\ && && && }
\end {equation}
\par Therefore, \I really do have an equivalence relation,
so some quotient \( F \times U \maplongto { \tilde j } E \) must exist
satisfying~(\ref {Gbundle}).
To define the bundle map \( E \maplongto p B \),
first note that this diagram commutes:
\begin {equation} \label {bundle construction projection cone}
\matrixxy
{ F \times \nerve U 2
\ar [dddd] _ { \id F \times \del j 0 }
\ar [rrrr] ^ { \id F \times \pair g { \del j 1 } } &&&&
F \times G \times U \ar [rrr] ^ { r \times \id U } &&&
F \times U \ar [dd] ^ { \hat F \times \id U } \\\\
&&&& &&& U \ar [dd] ^ j \\\\
F \times U \ar [rrrr] _ { \hat F \times \id U } &&&& U \ar [rrr] _ j &&& B }
\end {equation}
in fact, both sides are simply \( \hat F \times \nerve j 2 \).
Since \( E \) is \a quotient,
this defines \a map \( E \maplongto p B \)
such that (\ref {summary pullback}) commutes.
\par Therefore, \( E \)
is \a \G bundle over \( B \) with fibre \( F \)
associated with the transition \( g \).
\end {prf}
\par Now, \I would like to say
that \( E \) is \emph {the} associated \G bundle,
but this terminology is appropriate only if \( E \) is unique---
that is, unique up to isomorphism of \G bundles.
This is true; in fact, it is \a special case of this result:
\begin {prop} \label {associated bundle morphism}
Given \G bundles \( E \) and \( E ' \)
(both over \( B \) and with fibre \( F \))
associated with \G transition morphisms \( g \) and \( g ' \) (respectively),
and given \a \G transition morphism \( b \) from \( g \) to \( g ' \),
there is \a unique \G bundle morphism from \( E \) to \( E ' \)
associated with \( b \).
\end {prop}
\begin {prf}
First note that this diagram commutes:
\begin {equation}
\label {bundle morphism construction intermediate quotient cone}
\matrixxy
{ F \times U \cap \nerve { U ' } 2
\ar [dddddd] _ { \id F \times \del j { 0 1 } }
\ar [rrrrrrrrr] ^ { \id F \times \del j { 0 2 } } &&&& &&& &&
F \times U \cap U ' \ar [dd] ^ { \id F \times \pair b { \del j 1 } } \\\\
&&&& &&& && F \times G \times U ' \ar [dd] ^ { r \times \id U } \\\\
&&&& &&& && F \times U \ar [dd] ^ { \tilde j ' } \\\\
F \times U \cap U ' \ar [rrrr] _ { \id F \times \pair b { \del j 1 } } &&&&
F \times G \times U ' \ar [rrr] _ { r \times \id { U ' } } &&&
F \times U ' \ar [rr] _ { \tilde j ' } && E ' }
\end {equation}
because \( \jptildesub { \rsub w { \bsub x { y ' } } } { y ' }
= \jptildesub { \rsub w
{ \msub { \bsub x { x ' } } { \bsub { x ' } { y ' } } } }
{ y ' }
= \jptildesub { \rsub { \rsub w { \bsub x { x ' } } }
{ \bsub { x ' } { y ' } } }
{ y ' }
= \jptildesub { \rsub w { \bsub x { x ' } } } { x ' } \),
as also seen in this string diagram:
\begin {equation}
\label {bundle morphism construction intermediate quotient cone string}
\matrixxy
{ \dddStringl w && && \dstring b && && \ddString \\ && && && & y ' \\
& x & \ddstring b \rrstring \delta && && \dstring { g ' } \\\\
&& && x ' && \rString && \theta ' \dString \\\\ && && && && }
\end {equation}
Since \( \id F \times \del j 0 \), being \a cover,
is \a quotient of its kernel pair,
\I can construct from this \a map \( F \times U \maplongto { \tilde f } E \)
making this diagram commute:
\begin {equation} \label {bundle morphism construction intermediate quotient}
\matrixxy
{ F \times U \cap U '
\ar [ddddrr] _ { \id F \times \del j 0 }
\ar [rrrr] ^ { \id F \times \pair b { \del j 1 } } &&&&
F \times G \times U ' \ar [dd] ^ { r \times \id { U ' } } \\
&& \relax \ddricell { \tilde \zeta } \\
&& && F \times U \ar [dd] ^ { \tilde j ' } \\ && & \\
&& F \times U \ar [rr] _ { \tilde f } && E ' }
\end {equation}
In terms of elements
\( w \) of \( F \) and \( \pair x { x ' } \) of \( U \cap U ' \),
this result states
that \( \jptildesub { \rsub w { \bsub x { x ' } } } { x ' } \)
is equal to \( \pair w x \) composed with \( \tilde f \),
which can be drawn as this string diagram:
\begin {equation}
\label {bundle morphism construction intermediate quotient string}
\matrixxy { \ddStringl w && \dstring b && \dString \\ && & x ' \\
&& \rString && \tilde \zeta \\ && x \\ & }
\end {equation}
where \( \tilde f \)
is applied to the portion of the diagram open on the right.
\par Next, notice that this diagram commutes:
\begin {equation} \label {bundle morphism construction quotient cone}
\matrixxy
{ F \times \nerve U 2
\ar [dd] _ { \id F \times \del j 0 }
\ar [rrrr] ^ { \id F \times \pair g { \del j 1 } } &&&&
F \times G \times U \ar [rrr] ^ { r \times \id U } &&&
F \times U \ar [dd] ^ { \tilde f } \\\\
F \times U \ar [rrrrrrr] _ { \tilde f } &&&& &&& E ' }
\end {equation}
because \( \pair { \rsub w { \gsub x y } } y \) composed with \( \tilde f \)
is \( \jptildesub { \rsub { \rsub w { \gsub x y } } { \bsub y { x ' } } }
{ x ' }
= \jptildesub { \rsub w { \msub { \gsub x y } { \bsub y { x ' } } } } { x ' }
= \jptildesub { \rsub w { \bsub x { x ' } } } { x ' } \),
which is the composite with \( \tilde f \) of \( \pair w x \);
as also seen in this string diagram:
\begin {equation} \label {bundle morphism construction quotient cone string}
\matrixxy
{ \ddddStringl w \\\\
&& x && \dstring b \rrString && && \tilde \zeta \ddString \\
&& && && & x ' \\ && \ddstring g \rrstring \sigma && && \dstring b \\\\
&& && y && \rString && \tilde \zeta \\\\ && }
\end {equation}
(Notice that \I can drop \( x ' \)
because \( F \times \nerve U 2 \cap U '
\maplongto { \del j { 0 1 } } F \times \nerve U 2 \)
is epic.)
Since \( \tilde j \), being the pullback of \a cover,
is \a cover and hence the quotient of its kernel pair,
\I can construct from this \a map \( E \maplongto f E ' \)
making the diagram (\ref {Gbundle morphism}) commute.
\par In short, I've constructed \a \G bundle morphism \( E \maplongto f E ' \)
associated with \( b \).
Given any other \G bundle morphism \( E \maplongto { f ' } E ' \)
associated with \( b \),
let \( \tilde f ' \)
be the composite \( F \times U \maplongto { \tilde j } E \maplongto f E ' \).
Then \( \tilde f ' = \tilde f \)
by the unicity of maps out of quotients,
so \( f ' = f \) by that same principle.
\end {prf}
\par In particular, the identity \G transition morphism on \( g \)
becomes \a \G bundle morphism between any \G bundles associated with~\( g \).
Also, notice that the composite of \G bundle morphisms
associated with \a composable pair of \G transition morphisms
is itself associated with the composite \G transition morphism.
By the unicity clause of the previous proposition,
it is the same as the associated \G bundle morphism
constructed by the existence clause.
This shows, in particular, that any \G bundle morphism
associated with an identity \G transition morphism
is \a \G bundle isomorphism (an equivalence of \G bundles).
\par In fact, there is \a functor
from the category \( \Bun \C G B F \)
of \G bundles over \( B \) with fibre \( F \)
to the category \( B ^ G \) of \G transitions over \( B \),
defined by simply forgetting the total space \( E \).
The propositions above state (respectively)
that this functor is surjective and fully faithful.
In other words, \I have proved this theorem:
\begin {thm} \label {bundle cat equals transition cat}
Given any \G space \( F \),
the category \( \Bun \C G B F \) is equivalent to the category \( B ^ G \).
\end {thm}
So depending on the point of view desired,
you can think of \G transitions over \( B \) (a local view),
principal \G bundles over \( B \) (a global view),
or \G bundles over \( B \) with some convenient fibre \( F \)
(such as when \( G \) is defined as \a group of transformations of \( F \),
as is particularly common for linear groups);
the concepts are all equivalent.
\par It remains to explain my use of the terminology \q {associated}.
Normally, one begins with \a principal bundle
(\a bundle with fibre \( G \))
\( E _ G \)
and asks for the bundle \( E _ F \) with fibre \( F \)
associated with the given \( E _ G \).
But this is constructed by looking at the transition map of \( E _ G \)
and building \( E _ F \) out of that as \a certain quotient space.
So I've adapted the terminology
to say that the bundle is associated with the transition map \emph {directly}
rather than merely associated with the principal bundle
\emph {through} the transition map---
not because \I thought that this concept deserves the term better,
but because \I needed \emph {some} terminology to refer to it!
But really, all of the bundles associated with \a given transition map
should be considered to be associated with one another;
that is the only really fair way to look at it.
\par Theorem~\ref {bundle cat equals transition cat} is not \a new theorem.
Although I've found no reference that expresses it just like this,
it may be found implicitly in references on fibre bundles,
such as \cite [3.2, 8.2, 9.1] {Steenrod}.
However, its categorification,
Theorem~\ref {2bundle 2cat equals 2transition 2cat},
is the central internal result of this paper.
Accordingly, \I now turn to part~\ref 2
with \a categorification of Theorem~\ref {bundle cat equals transition cat}
as my goal.
\part {Categorified bundles} \label 2
Now \I categorify the above
to construct the theory of \two bundles with \a structure \two group.
\section {\Two categorical preliminaries} \label {2prelim}
First, \I will turn the category theory from section~\ref {1prelim}
into \two category theory.
\par Just as part~\ref 1 was about \a specific category \( \C \),
so this part~\ref 2 is about \a specific \two category \( \CCC \).
But just as the details of \( \C \) were largely irrelevant in part~\ref 1,
so the details of \( \CCC \) will be largely irrelevant here.
\I will explain what \( \CCC \) is in section~\ref {2C},
but all that matters outside that section
is that \( \CCC \) supports the structures described here.
\par \I should mention, however,
that the \two category to be described in section~\ref {2C}
is (unlike in \cite {HDA5}, \cite {HDA6},
and some earlier versions of this paper)
\a \emph {weak} \two category (that is \a bicategory),
while \I will mostly treat it as if it were \a \emph {strict} \two category
(the more familiar sort of \two category).
This is all right, however,
since \a generalisation of the Mac Lane Coherence Theorem (see \cite {bicat})
shows that every bicategory is equivalent to \a strict \two category.
\subsection {Notation and terminology} \label {2term}
\A \strong {\two space} is an object in \( \CCC \),
and \a \strong {\two map} is \a morphism in \( \CCC \),
and \a \strong {natural transformation} is \a \two morphism in \( \CCC \).
Uppercase calligraphic letters like \( \XX \) denote \two spaces,
lowercase Fraktur letters like \( \xx \) denote \two maps,
and lowercase Greek letters like \( \psi \) denote natural transformations.
However, the identity \two map on \( \XX \) will be denoted \( \id \XX \),
and the identity natural transformation on \( \xx \)
will be denoted \( \iid \xx \).
(Combining these, even \( \idid \XX \) is possible.)
\par In general, the names of \two maps and natural transformations
will label arrows directly,
perhaps in \a small inline diagram
like \( \XX \maplongto \xx \YY \maplongto \yy \ZZ \),
which denotes \a composition of \two maps,
or perhaps in \a huge \two cell diagram
which describes the composition of several natural transformations.
\I will endeavour to make such diagrams easy to read
by orienting \two maps to the right when convenient, or if not then downwards;
natural transformations will always go down and/or to the left.
\par As in part~\ref 1,
my discussion is in almost purely arrow-theoretic terms.
However, \I will again use the generalised arrow-theoretic concept of element.
Specifically, an \strong {element} \( \xx \) of \a \two space \( \YY \)
is \a \two space \( \XX \)
together with \a \two map \( \XX \maplongto \xx \YY \).
In certain contexts, the \two map \( \xx \)
has the same information in it
as \a point in the space of objects of \( \YY \).
On the other hand, if \( \XX \) is chosen differently,
then the map \( \xx \) may described more complicated features,
such as \a curve of objects, \a point in the space of arrows,
or even (say) \a surface of arrows between two curves of objects.
The most general element is in fact the identity \two map
\( \XX \maplongto { \id \XX } \XX \).
You (as reader) may imagine the elements as set-theoretic points if that helps,
but their power in proofs lies in their complete generality.
\I will also use morphism-elements, or \strong {arrows},
that is natural transformations between these generalised elements.
\subsection {\Two products} \label {2prod}
Given \two spaces \( \XX \) and \( \YY \),
there is \a \two space \( \XX \times \YY \),
the \strong {Cartesian \two product} of \( \XX \) and \( \YY \).
One can also form more general Cartesian \two products,
like \( \XX \times \YY \times \ZZ \) and so on.
There is also \a trivial \two space \( 1 \),
which is the Cartesian \two product of no \two spaces.
\par The generalised elements of \( \XX \times \YY \)
may be taken to be ordered pairs;
that is, given \two maps \( \XX \maplongto \xx \YY \)
and \( \XX \maplongto { \xx ' } \YY \),
there is \a \strong {pairing \two map}
\( \XX \maplongto { \pair \xx { \xx ' } } \YY \times \YY ' \);
conversely, given \a \two map \( \XX \maplongto \yy \YY \times \YY ' \),
there is \a pair of \two maps
\( \XX \maplongto \xx \YY \) and \( \XX \maplongto { \xx ' } \YY ' \)
such that \( \yy \) is naturally isomorphic to \( \pair \xx { \xx ' } \).
(That is, there is an invertible natural transformation
between \( \yy \) and \( \pair \xx { \xx ' } \).)
Furthermore, this pair is unique up to natural isomorphism;
and the automorphisms of such \a pair
are in bijective correspondence with the automorphisms of \( \yy \).
This extends to \two products of multiple \two spaces.
Additionally, given any \two space \( \XX \),
there is \a \strong {trivial \two map} \( \XX \maplongto { \hat \XX } 1 \);
it's unique in the sense that
every map \( \XX \maplongto \yy 1 \) is isomorphic to \( \hat \XX \),
and this isomorphism is itself unique.
\par There are also product \two maps;
given \two maps
\( \XX \maplongto \xx \YY \) and \( \XX ' \maplongto { \xx ' } \YY ' \),
the \strong {product \two map}
is \( \XX \times \XX ' \maplongto { \xx \times \xx ' } \YY \times \YY ' \).
If furthermore there are \two maps
\( \XX \maplongto \yy \YY \) and \( \XX ' \maplongto { \yy ' } \YY ' \),
with natural transformations
from \( \xx \) to \( \yy \) and from \( \xx ' \) to \( \yy ' \),
then there is \a product natural transformation
from \( \xx \times \xx ' \) to \( \yy \times \yy ' \).
This respects all of the \two category operations.
The Mac Lane Coherence Theorem even applies,
allowing me to make
all of the same abuses of notation that appeared in part~\ref 1.
\subsection {\Two pullbacks} \label {2pull}
Unlike products, pullbacks do not always exist in \( \C \),
hence neither do \two pullbacks in \( \CCC \).
Thus, \I will have to treat this in more detail,
so that \I can discuss exactly what it means for \a \two pullback to exist.
\par \A \strong {\two pullback diagram}
consists of \two spaces \( \XX \), \( \YY \), and \( \ZZ \),
and \two maps \( \XX \maplongto \xx \ZZ \) and \( \YY \maplongto \yy \ZZ \):
\begin {equation} \label {2pullback diagram}
\matrixxy { & & \XX \ar [dd] ^ \xx \\ \\ \YY \ar [rr] _ \yy & & \ZZ }
\end {equation}
Given \a \two pullback diagram, \a \strong {\two pullback cone}
is \a \two space \( \CC \)
together with \two maps
\( \CC \maplongto \zz \XX \) and \( \CC \maplongto \ww \YY \)
and \a natural isomorphism \( \omega \):
\begin {equation} \label {2pullback cone}
\matrixxy
{ \CC \ar [dd] _ \ww \ar [rr] ^ \zz \ddrriicell \omega & &
\XX \ar [dd] ^ \xx \\ \\
\YY \ar [rr] _ \yy & & \ZZ }
\end {equation}
(So while in section~\ref {1pull},
\( \omega \) merely labelled \a commutative diagram,
here it refers to \a specific natural transformation!)
Given \two pullback cones \( \CC \) and \( \CC ' \),
\a \strong {\two pullback cone morphism} from \( \CC \) to \( \CC ' \)
is \a \two map \( \CC \maplongto \uu \CC ' \)
and natural isomorphisms \( \chi \) and \( \psi \):
\begin {equation} \label {2pullback cone morphism}
\matrixxy
{ \CC \ar [dr] _ \uu \ar `r _dr [drrr] ^ \zz [drrr] &
\relax \driicell \chi \\
& \CC ' \ar [rr] _ { \zz ' } && \XX }
\qquad
\matrixxy
{ \CC \ar `d ^dr [dddr] _ \ww [dddr] \ar [dr] ^ \uu \\
\relax \driicell \psi & \CC ' \ar [dd] ^ { \ww ' } \\ & \\ & \YY }
\end {equation}
Furthermore, to have \a \two pullback cone morphism,
the following composite of natural transformations
must be equal to \( \omega \):
\begin {equation} \label {2pullback cone morphism coherence}
\matrixxy
{ \CC
\ar `d ^dr [dddr] _ \ww [dddr] \ar [dr] ^ \uu \ar `r _dr [drrr] ^ \zz [drrr] &
\relax \driicell \chi \\
\relax \driicell \psi &
\CC ' \ar [dd] ^ { \ww ' } \ar [rr] _ { \zz ' } \ddrriicell { \omega ' } &&
\XX \ar [dd] ^ \xx \\ & \\
& \YY \ar [rr] _ \yy && \ZZ }
\end {equation}
This is \a typical example of \a \emph {coherence law};
it's the specification of coherence laws
that makes \a categorified theory logically richer than the original theory.
\par \A \strong {\two pullback} of the given \two pullback diagram
is \a \two pullback cone \( \PP \)
that is \emph {universal}
in the sense that, given any other \two pullback cone \( \CC \),
there is \a \two pullback cone morphism \( \CC \maplongto \uu \PP \)
such that, given any other \two pullback cone morphism
\( \CC \maplongto { \uu ' } \PP \),
there is \a unique natural isomorphism \( \nu \):
\begin {equation} \label {2pullback universality}
\matrixxy { \CC \drfulliicell { \uu ' } \uu \nu \\ & \PP }
\end {equation}
such that the following natural transformations:
\begin {equation} \label {2pullback universality coherence}
\matrixxy
{ \CC
\ar [dr] ^ { \uu ' } \ar `r _dr [drrr] ^ \zz [drrr] \drloweriicell \uu \nu &
\relax \driicell \chi \\
& \PP \ar [rr] _ { \zz _ \PP } && \XX }
\qquad
\matrixxy
{ \CC
\ar `d ^dr [dddr] _ \ww [dddr] \ar [dr] _ \uu
\drupperiicell { \uu ' } \nu \\
\relax \driicell \psi & \PP \ar [dd] ^ { \ww _ \PP } \\ & \\ & \YY }
\end {equation}
are equal (respectively) to \( \chi ' \) and \( \psi ' \).
\par As with pullbacks in \( \C \),
so \two pullbacks in \( \CCC \) needn't always exist.
\par In the rest of this paper,
\I will want to define certain spaces as pullbacks, when they exist.
If such definitions are to be sensible,
then \I must show that it doesn't matter which pullback one uses.
\begin {prop} \label {2pullback unicity}
Given the pullback diagram (\ref {2pullback diagram})
and pullbacks \( \PP \) and \( \PP ' \),
the \two spaces \( \PP \) and \( \PP ' \) are equivalent in \( \CCC \).
\end {prop}
\begin {prf}
Since \( \PP \) is universal,
there is \a pullback cone morphism \( \PP ' \maplongto \uu \PP \):
\begin {equation} \label {2pullback morphism}
\matrixxy
{ \PP '
\ar `d ^dr [dddr] _ { \ww ' } [dddr] \ar [dr] ^ \uu
\ar `r _dr [drrr] ^ { \zz ' } [drrr] &
\relax \driicell \chi \\
\relax \driicell \psi &
\PP \ar [dd] ^ \ww \ar [rr] _ \zz \ddrriicell \omega & &
\XX \ar [dd] ^ \xx \\ & \\
& \YY \ar [rr] _ \yy & & \ZZ }
\end {equation}
Since \( \PP ' \) is universal,
there is also \a pullback cone morphism
\( \PP \maplongto { \bar \uu } \PP ' \):
\begin {equation} \label {competing 2pullback morphism}
\matrixxy
{ \PP
\ar `d ^dr [dddr] _ \ww [dddr] \ar [dr] ^ { \bar \uu }
\ar `r _dr [drrr] ^ \zz [drrr] &
\relax \driicell { \bar \chi } \\
\relax \driicell { \bar \psi } &
\PP ' \ar [dd] ^ { \ww ' } \ar [rr] _ { \zz ' } \ddrriicell { \omega ' } & &
\XX \ar [dd] ^ \xx \\ & \\
& \YY \ar [rr] _ \yy & & \ZZ }
\end {equation}
Composing these one way,
\I get \a pullback cone morphism from \( \PP \) to itself:
\begin {equation} \label {composed 2pullback automorphism}
\matrixxy
{ \PP
\ar `d ^r [ddddrr] _ \ww [ddddrr] \ar [dr] ^ { \bar \uu }
\ar `r _d [ddrrrr] ^ \zz [ddrrrr] & &
\relax \driicell { \bar \chi } \\
& \PP '
\ar `d ^dr [dddr] _ { \ww ' } [dddr] \ar [dr] ^ \uu
\ar `r _dr [drrr] ^ { \zz ' } [drrr] &
\relax \driicell \chi & \\
\relax \driicell { \bar \psi } & \relax \driicell \psi &
\PP \ar [dd] ^ \ww \ar [rr] _ \zz \ddrriicell \omega & & \XX \ar [dd] ^ \xx \\
& & \\ & & \YY \ar [rr] _ \yy & & \ZZ }
\end {equation}
But since the identity on \( \PP \) is also \a pullback cone morphism,
the universal property
shows that \( \PP \maplongto { \bar \uu } \PP ' \maplongto \uu \PP \)
is naturally isomorphic to the identity on \( \PP \).
Similarly, \( \PP ' \maplongto \uu \PP \maplongto { \bar \uu } \PP ' \)
is naturally isomorphic to the identity on \( \PP ' \).
Therefore, \( \PP \) and \( \PP ' \) are equivalent.
\end {prf}
\subsection {Equivalence \two relations} \label {2eqrel}
Given \a \two space \( \UU \),
\a \strong {binary \two relation} on \( \UU \)
is \a \two space \( \nerve \RR 2 \) equipped with \two maps
\( \nerve \RR 2 \maplongto { \del \jj 0 } \UU \)
and \( \nerve \RR 2 \maplongto { \del \jj 1 } \UU \)
that are \emph {jointly \two monic};
this means that given generalised elements
\( \XX \maplongto \xx \nerve \RR 2 \) and \( \XX \maplongto \yy \nerve \RR 2 \)
of \( \nerve \RR 2 \)
and natural isomorphisms \( \del \chi 0 \) and \( \del \chi 1 \) as follows:
\begin {equation} \label {2relation}
\matrixxy
{ \XX \ar [dd] _ \yy \ar [drr] ^ \xx \\
\relax \drriicell { \del \chi 0 } &&
\relax \nerve \RR 2 \ar [dd] ^ { \del \jj 0 } \\
\relax \nerve \RR 2 \ar [drr] _ { \del \jj 0 } && \\ && \UU }
\qquad
\matrixxy
{ \XX \ar [ddr] _ \yy \ar [rr] ^ \xx & \relax \ddriicell { \del \chi 1 } &
\relax \nerve \RR 2 \ar [ddr] ^ { \del \jj 1 } \\\\
& \relax \nerve \RR 2 \ar [rr] _ { \del \jj 1 } & & \UU }
\end {equation}
then there is \a unique natural isomorphism \( \xx \maplongTo \chi \yy \)
such that the following composites
are equal, respectively, to \( \del \chi 0 \) and \( \del \chi 1 \):
\begin {equation} \label {2relation coherence}
\matrixxy
{ \XX \ar [dd] _ \yy \ar [drr] ^ \xx \ddriicell \chi \\
& \relax \ddreqcell &
\relax \nerve \RR 2
\ar [dll] ^ { \id { \nerve \RR 2 } } \ar [dd] ^ { \del \jj 0 } \\
\relax \nerve \RR 2 \ar [drr] _ { \del \jj 0 } & \\ && \UU }
\qquad
\matrixxy
{ \XX \ar [ddr] _ \yy \ar [rr] ^ \xx \drriicell \chi & &
\relax \nerve \RR 2
\ar [ddl] ^ { \id { \nerve \RR 2 } } \ar [ddr] ^ { \del \jj 1 } \\
& \relax \drreqcell & \\
& \relax \nerve \RR 2 \ar [rr] _ { \del \jj 1 } & & \UU }
\end {equation}
The upshot of this is that an element of \( \nerve \RR 2 \)
is determined, up to unique isomorphism, by two elements of \( \UU \).
\par \A binary \two relation is \strong {reflexive}
if it is equipped with \a \emph {reflexivity \two map}
\( \UU \maplongto { \del \jj { 0 0 } } \nerve \RR 2 \)
and natural isomorphisms \( \del \omega 0 \) and \( \del \omega 1 \):
\begin {equation} \label {2reflexivity}
\matrixxy
{ \UU \ar `d ^dr [dddr] _ { \id \UU } [dddr] \ar [dr] ^ { \del \jj { 0 0 } } \\
\relax \driicell { \del \omega 0 } &
\relax \nerve \RR 2 \ar [dd] ^ { \del \jj 0 } \\ & \\
& \UU }
\qquad
\matrixxy
{ \UU \ar [dr] _ { \del \jj { 0 0 } } \ar `r _dr [drrr] ^ { \id \UU } [drrr] &
\relax \driicell { \del \omega 1 } \\
& \relax \nerve \RR 2 \ar [rr] _ { \del \jj 1 } && \UU }
\end {equation}
In terms of generalised elements,
given any element \( \XX \maplongto \xx \UU \) of \( \UU \),
composing with \( \del \jj { 0 0 } \)
gives an element of \( \nerve \RR 2 \)
that corresponds to \( \pair \xx \xx \) in \( \UU \times \UU \).
By \two monicity, the reflexivity \two map is unique up to unique isomorphism.
\par Assume that the \two kernel pair of \( \del \jj 0 \) exists,
and let it be the \two space \( \nerve \RR 3 \):
\begin {equation} \label {equivalence 2relation triple}
\matrixxy
{ \relax \nerve \RR 3
\ar [dd] _ { \del \jj { 0 1 } } \ar [drr] ^ { \del \jj { 0 2 } } \\
\relax \drriicell { \del \omega { 0 0 } } &&
\relax \nerve \RR 2 \ar [dd] ^ { \del \jj 0 } \\
\relax \nerve \RR 2 \ar [drr] _ { \del \jj 0 } && \\ && \UU }
\end {equation}
\A binary \two relation is \strong {right Euclidean}
if it is equipped with \a (right) \emph {Euclideanness \two map}
\( \nerve \RR 3 \maplongto { \del \jj { 1 2 } } \nerve \RR 2 \)
and natural isomorphisms
\( \del \omega { 0 1 } \) and \( \del \omega { 1 1 } \):
\begin {equation} \label {2Euclideanness}
\matrixxy
{ \relax \nerve \RR 3
\ar [dd] _ { \del \jj { 0 1 } } \ar [rr] ^ { \del \jj { 1 2 } }
\ddrriicell { \del \omega { 0 1 } } &&
\relax \nerve \RR 2 \ar [dd] ^ { \del \jj 1 } \\\\
\relax \nerve \RR 2 \ar [rr] _ { \del \jj 0 } && \UU }
\qquad
\matrixxy
{ \relax \nerve \RR 3
\ar [ddr] _ { \del \jj { 0 2 } } \ar [rr] ^ { \del \jj { 1 2 } } &
\relax \ddriicell { \del \omega { 1 1 } } &
\relax \nerve \RR 2 \ar [ddr] ^ { \del \jj 1 } \\\\
& \relax \nerve \RR 2 \ar [rr] _ { \del \jj 1 } & & \UU }
\end {equation}
Like the reflexivity \two map, any Euclideanness \two map
is unique up to unique isomorphism.
An element of \( \nerve \RR 3 \)
is given by elements \( \xx \), \( \yy \), and \( \zz \) of \( \UU \)
such that \( \pair \xx \yy \) and \( \pair \xx \zz \)
give elements of \( \nerve \RR 2 \);
then by the Euclideanness \two map,
\( \pair \yy \zz \) also gives an element of \( \nerve \RR 2 \).
\par Putting these together,
an \strong {equivalence \two relation}
is \a binary \two relation that is both reflexive and Euclidean.
\par Given an equivalence \two relation as above,
\a \strong {\two quotient} of the equivalence \two relation
is \a \two space \( \nerve \RR 0 \)
and \a map \( \UU \maplongto \jj \nerve \RR 0 \)
that \emph {\two coequalises} \( \del \jj 0 \) and \( \del \jj 1 \).
This means, first, \a natural isomoprhism \( \omega \):
\begin {equation} \label {2quotient}
\matrixxy
{ \relax \nerve \RR 2
\ar [dd] _ { \del \jj 0 } \ar [rr] ^ { \del \jj 1 } \ddrriicell \omega &&
\UU \ar [dd] ^ \jj \\\\
\UU \ar [rr] _ \jj && \relax \nerve \RR 0 }
\end {equation}
such that these two composite isomorphisms are equal:
\begin {equation} \label {2quotient coherence reflexivity}
\matrixxy
{ \UU
\ar `d ^dr [ddddrr] [ddddrr] _ { \id \UU } \ar [ddrr] |- { \del \jj { 0 0 } }
\ar `r _dr [ddrrrr] [ddrrrr] ^ { \id \UU } &
\relax \ddrriicell { \del \omega 1 } \\
\relax \ddrriicell { \del \omega 0 } \\
&& \relax \nerve \RR 2
\ar [dd] ^ { \del \jj 0 } \ar [rr] _ { \del \jj 1 } \ddrriicell \omega &&
\UU \ar [dd] ^ \jj \\ && \\
&& \UU \ar [rr] _ \jj && \nerve \RR 0 }
\qquad = \qquad
\matrixxy
{ \UU
\ar `d ^dr [ddddrr] [ddddrr] _ { \id \UU }
\ar `r _dr [ddrrrr] [ddrrrr] ^ { \id \UU } \\
& \relax \ddrreqcell \\ && && \UU \ar [dd] ^ \jj \\
&& & \\ && \UU \ar [rr] _ \jj && \relax \nerve \RR 0 }
\end {equation}
and also these two composite isomorphisms are equal:
\begin {equation} \label {2quotient coherence Euclid}
\matrixxy
{ \relax \nerve \RR 3
\ar [dd] _ { \del \jj { 0 1 } } \ar [ddrr] |- { \del \jj { 0 2 } }
\ar [rr] ^ { \del \jj { 1 2 } } &
\relax \ddrriicell { \del \omega { 0 0 } } &
\relax \nerve \RR 2 \ar [ddrr] ^ { \del \jj 1 } \\
\relax \ddrriicell { \del \omega { 1 1 } } \\
\relax \nerve \RR 2 \ar [ddrr] _ { \del \jj 0 } &&
\relax \nerve \RR 2
\ar [dd] ^ { \del \jj 0 } \ar [rr] _ { \del \jj 1 } \ddrriicell \omega &&
\UU \ar [dd] ^ \jj \\
&& \\ && \UU \ar [rr] _ \jj && \relax \nerve \RR 0 }
\qquad = \qquad
\matrixxy
{ \relax \nerve \RR 3
\ar [dd] _ { \del \jj { 0 1 } } \ar [rr] ^ { \del \jj { 1 2 } } &&
\relax \nerve \RR 2 \ar [ddrr] ^ { \del \jj 1 } \\
& \relax \driicell { \del \omega { 0 1 } } \\
\relax \nerve \RR 2 \ar [ddrr] _ { \del \jj 0 } && &&
\UU \ar [ddll] _ { \id \UU } \ar [dd] ^ \jj \\ && & \relax \dreqcell \\
&& \UU \ar [rr] _ \jj && \relax \nerve \RR 0 }
\end {equation}
Secondly, given any \two map \( \UU \maplongto \xx \nerve \RR 0 \)
and natural isomorphism \( \omega _ \xx \):
\begin {equation} \label {2quotient cone}
\matrixxy
{ \relax \nerve \UU 2
\ar [dd] _ { \del \jj 0 } \ar [rr] ^ { \del \jj 1 }
\ddrriicell { \omega _ \xx } &&
\UU \ar [dd] ^ \xx \\\\
\UU \ar [rr] _ \xx && \XX }
\end {equation}
if these two composite isomorphisms are equal:
\begin {equation} \label {2quotient cone coherence reflexivity}
\matrixxy
{ \UU
\ar `d ^dr [ddddrr] [ddddrr] _ { \id \UU } \ar [ddrr] |- { \del \jj { 0 0 } }
\ar `r _dr [ddrrrr] [ddrrrr] ^ { \id \UU } &
\relax \ddrriicell { \del \omega 1 } \\
\relax \ddrriicell { \del \omega 0 } \\
&& \relax \nerve \RR 2
\ar [dd] ^ { \del \jj 0 } \ar [rr] _ { \del \jj 1 }
\ddrriicell { \omega _ \xx } &&
\UU \ar [dd] ^ \xx \\ && \\
&& \UU \ar [rr] _ \xx && \nerve \RR 0 }
\qquad = \qquad
\matrixxy
{ \UU
\ar `d ^dr [ddddrr] [ddddrr] _ { \id \UU }
\ar `r _dr [ddrrrr] [ddrrrr] ^ { \id \UU } \\
& \relax \ddrreqcell \\ && && \UU \ar [dd] ^ \xx \\
&& & \\ && \UU \ar [rr] _ \xx && \relax \nerve \RR 0 }
\end {equation}
and also these two composite isomorphisms are equal:
\begin {equation} \label {2quotient cone coherence Euclid}
\matrixxy
{ \relax \nerve \RR 3
\ar [dd] _ { \del \jj { 0 1 } } \ar [ddrr] |- { \del \jj { 0 2 } }
\ar [rr] ^ { \del \jj { 1 2 } } &
\relax \ddrriicell { \del \omega { 0 0 } } &
\relax \nerve \RR 2 \ar [ddrr] ^ { \del \jj 1 } \\
\relax \ddrriicell { \del \omega { 1 1 } } \\
\relax \nerve \RR 2 \ar [ddrr] _ { \del \jj 0 } &&
\relax \nerve \RR 2
\ar [dd] ^ { \del \jj 0 } \ar [rr] _ { \del \jj 1 }
\ddrriicell { \omega _ \xx } &&
\UU \ar [dd] ^ \xx \\
&& \\ && \UU \ar [rr] _ \xx && \relax \nerve \RR 0 }
\qquad = \qquad
\matrixxy
{ \relax \nerve \RR 3
\ar [dd] _ { \del \jj { 0 1 } } \ar [rr] ^ { \del \jj { 1 2 } } &&
\relax \nerve \RR 2 \ar [ddrr] ^ { \del \jj 1 } \\
& \relax \driicell { \del \omega { 0 1 } } \\
\relax \nerve \RR 2 \ar [ddrr] _ { \del \jj 0 } && &&
\UU \ar [ddll] _ { \id \UU } \ar [dd] ^ \xx \\ && & \relax \dreqcell \\
&& \UU \ar [rr] _ \xx && \relax \nerve \RR 0 }
\end {equation}
then there is \a \two map \( \nerve \RR 0 \maplongto { \tilde \xx } \XX \)
and \a natural isomorphism \( \tilde \omega _ \xx \):
\begin {equation} \label {2quotient universality}
\matrixxy
{ \UU \ar [dr] _ \jj \ar `r _dr [drrr] ^ \xx [drrr] &
\relax \driicell { \tilde \omega _ \xx } \\
& \relax \nerve \RR 0 \ar [rr] _ { \tilde \xx } && \XX }
\end {equation}
such that these composite isomorphisms are equal:
\begin {equation} \label {2quotient universality coherence}
\matrixxy
{ \relax \nerve \UU 2
\ar [dd] _ { \del \jj 0 } \ar [rr] ^ { \del \jj 1 } &
\relax \ddriicell { \omega _ \xx } &
\UU \ar `dr _d [ddddrr] ^ \xx [ddddrr] \\\\
\UU \ar [ddrr] _ \jj \ar `r _dr [ddrrrr] ^ \xx [ddrrrr] &
\relax \ddrriicell { \tilde \omega _ \xx } & \\\\
&& \relax \nerve \RR 0 \ar [rr] _ { \tilde \xx } && \XX }
\qquad = \qquad
\matrixxy
{ \relax \nerve \UU 2
\ar [dd] _ { \del \jj 0 } \ar [rr] ^ { \del \jj 1 } &&
\UU \ar [dddd] _ j \ar `dr _d [ddddrr] ^ \xx [ddddrr] \\
\relax \drriicell \omega \\ \UU \ar [ddrr] _ \jj &&
\relax \drriicell { \tilde \omega _ \xx } \\
&& && \\ && \relax \nerve \RR 0 \ar [rr] _ { \tilde \xx } && \XX }
\end {equation}
Finally, given any alternative \two map
\( \nerve \RR 0 \maplongto { \tilde \xx ' } \XX \)
and natural isomorphism \( \tilde \omega _ \xx ' \):
\begin {equation} \label {2quotient universality alternative}
\matrixxy
{ \UU \ar [dr] _ \jj \ar `r _dr [drrr] ^ \xx [drrr] &
\relax \driicell { \tilde \omega _ \xx ' } \\
& \relax \nerve \RR 0 \ar [rr] _ { \tilde \xx ' } && \XX }
\end {equation}
if these composite isomorphisms are equal:
\begin {equation} \label {2quotient universality alternative coherence}
\matrixxy
{ \relax \nerve \UU 2
\ar [dd] _ { \del \jj 0 } \ar [rr] ^ { \del \jj 1 } &
\relax \ddriicell { \omega _ \xx } &
\UU \ar `dr _d [ddddrr] ^ \xx [ddddrr] \\\\
\UU \ar [ddrr] _ \jj \ar `r _dr [ddrrrr] ^ \xx [ddrrrr] &
\relax \ddrriicell { \tilde \omega _ \xx ' } & \\\\
&& \relax \nerve \RR 0 \ar [rr] _ { \tilde \xx ' } && \XX }
\qquad = \qquad
\matrixxy
{ \relax \nerve \UU 2
\ar [dd] _ { \del \jj 0 } \ar [rr] ^ { \del \jj 1 } &&
\UU \ar [dddd] _ j \ar `dr _d [ddddrr] ^ \xx [ddddrr] \\
\relax \drriicell \omega \\ \UU \ar [ddrr] _ \jj &&
\relax \drriicell { \tilde \omega _ \xx ' } \\
&& && \\ && \relax \nerve \RR 0 \ar [rr] _ { \tilde \xx ' } && \XX }
\end {equation}
then there exists \a unique isomorphism \( \xx \maplongTo \nu \xx ' \)
such that this isomorphism:
\begin {equation} \label {2quotient 2universality}
\matrixxy
{ \UU \ar [dr] _ \jj \ar `r _dr [drrr] ^ \xx [drrr] &
\relax \driicell { \tilde \omega _ \xx } \\
& \relax \nerve \RR 0
\ar [rr] ^ { \tilde \xx } \rrloweriicell { \xx ' } \nu && \XX }
\end {equation}
is equal to~\( \tilde \omega _ \xx ' \).
\par Thus if an equivalence \two relation has \a \two quotient,
then \I can define \a \two map out of this \two quotient,
up to unique natural isomorphism,
by defining \a \two map out of the relation's base space
with certain attendant isomorphisms satisfying appropriate coherence relations.
\subsection {\Two covers} \label {2cover}
Analogously to section~\ref {1cover}, here are the axioms for \two covers:
\begin {closeitemize}
\item All equivalences are \two covers;
\item \A composite of \two covers is \a \two cover;
\item The \two pullback of \a \two cover along any \two map
exists and is \a \two cover;
\item The \two quotient
of every equivalence \two relation involving \a \two cover
exists and is \a \two cover;
and \item Every \two cover is the \two quotient of its \two kernel pair,
with the same natural isomorphism \( \omega \) involved in each.
\end {closeitemize}
\par \I will be particularly interested
in the \two pullback of \( \UU \) along itself,
that is the \two pullback of this diagram:
\begin {equation} \label {2cover 2kernel pair diagram}
\matrixxy { & & \UU \ar [dd] ^ \jj \\ \\ \UU \ar [rr] _ \jj & & \BB }
\end {equation}
This \two pullback, \( \nerve \UU 2 \),
is the \strong {\two kernel pair} of \( \UU \),
giving the \two pullback diagram:
\begin {equation} \label {2cover 2kernel pair}
\matrixxy
{ \relax \nerve \UU 2
\ar [dd] _ { \del \jj 1 } \ar [rr] ^ { \del \jj 0 } \ddrriicell \omega & &
\UU \ar [dd] ^ \jj \\ \\
\UU \ar [rr] _ \jj & & \BB }
\end {equation}
\par As before,
the \two space \( \nerve \UU 3 \) may also be defined as \a \two pullback:
\begin {equation} \label {2cover 2kernel triple}
\matrixxy
{ \relax \nerve \UU 3
\ar [dd] _ { \del \jj { 1 2 } } \ar [rr] ^ { \del \jj { 0 1 } } \ddrreqcell & &
\relax \nerve \UU 2 \ar [dd] ^ { \del \jj 1 } \\ \\
\relax \nerve \UU 2 \ar [rr] _ { \del \jj 0 } & & \UU }
\end {equation}
And the \two space \( \nerve \UU 4 \):
\begin {equation} \label {2cover 2kernel quadruple}
\matrixxy
{ \relax \nerve \UU 4
\ar [dd] _ { \del \jj { 1 2 3 } } \ar [rr] ^ { \del \jj { 0 1 2 } }
\ddrreqcell & &
\relax \nerve \UU 3 \ar [dd] ^ { \del \jj { 1 2 } } \\ \\
\relax \nerve \UU 3 \ar [rr] _ { \del \jj { 0 1 } } & & \relax \nerve \UU 2 }
\end {equation}
And so forth, very much as in section~\ref {1cover}.
\par In fact, all of these natural isomorphisms
based on \( \omega \)
can be ignored using (again) the Mac Lane Coherence Theorem,
applied to the \two products in the slice \two category \( \slice \CCC \BB \).
\section {The \two category of \two spaces} \label {2C}
To define my \two category \( \CCC \),
\I need to categorify the notion of space to get \a notion of \two space.
If \a space is \a smooth manifold,
then \a \two space is \a \emph {Lie category}, or \emph {smooth category},
\a category whose objects and morphisms form smooth manifolds.
Then \two spaces
(together with \two maps and \a new feature, natural transformations)
will form my \two category~\( \CCC \).
\subsection {\Two spaces} \label {2space}
It was Charles Ehresmann~\cite {Ehresmann}
who first defined the notion of differentiable category;
today, one understands this as \a special case
of the general notion of \emph {internal category}~\cite [chapter~8] {Borceux}.
\par Following this,
\a \strong {\two space} consists of
\a space \( \cat \XX 1 \) (the space of objects, or points)
and \a space \( \cat \XX 2 \) (the space of morphisms, or arrows)
together with maps
\( \cat \XX 2 \maplongto { \dom d 0 } \cat \XX 1 \) (the source map)
and \( \cat \XX 2 \maplongto { \dom d 1 } \cat \XX 1 \) (the target map).
This creates \a pullback diagram:
\begin {equation} \label {comppairdiagram}
\matrixxy
{ & & \relax \cat \XX 2 \ar [dd] ^ { \dom d 1 } \\ \\
\relax \cat \XX 2 \ar [rr] _ { \dom d 0 } & & \relax \cat \XX 1 }
\end {equation}
To have \a \two space,
this pullback diagram must have \a pullback
\( \cat \XX 3 \) (the space of composable pairs of arrows),
giving \a commutative diagram:
\begin {equation} \label {comppair}
\matrixxy
{ \relax \cat \XX 3
\ar [dd] _ { \dom d { 1 2 } } \ar [rr] ^ { \dom d { 0 1 } } & &
\relax \cat \XX 2 \ar [dd] ^ { \dom d 1 } \\ \\
\relax \cat \XX 2 \ar [rr] _ { \dom d 0 } & & \relax \cat \XX 1 }
\end {equation}
Finally, consider the pullback diagram
\begin {equation} \label {comptriplediagram}
\matrixxy
{ & & \relax \cat \XX 3 \ar [dd] ^ { \dom d { 1 2 } } \\ \\
\relax \cat \XX 3 \ar [rr] _ { \dom d { 0 1 } } & & \relax \cat \XX 2 }
\end {equation}
\A pullback \( \cat \XX 4 \) (the space of composable triples) must exist:
\begin {equation} \label {comptriple}
\matrixxy
{ \relax \cat \XX 4
\ar [dd] _ { \dom d { 1 2 3 } } \ar [rr] ^ { \dom d { 0 1 2 } } & &
\relax \cat \XX 3 \ar [dd] ^ { \dom d { 1 2 } } \\ \\
\relax \cat \XX 3 \ar [rr] _ { \dom d { 0 1 } } & & \relax \cat \XX 2 }
\end {equation}
\par \A \two space is additionally equipped with
maps \( \cat \XX 1 \maplongto { \dom d { 0 0 } } \cat \XX 2 \)
(taking an object to its identity morphism, or \a point to its null arrow)
and \( \cat \XX 3 \maplongto { \dom d { 0 2 } } \cat \XX 2 \)
(composing \a composable pair of arrows),
such that each of these diagrams commutes:
\begin {equation} \label {2space dom}
\matrixxy
{ \\ \relax \cat \XX 3
\ar [ddrr] _ { \dom d { 0 2 } } \ar [rr] ^ { \dom d { 0 1 } } & &
\relax \cat \XX 2 \ar [ddrr] ^ { \dom d 0 } \\ \\
& & \relax \cat \XX 2 \ar [rr] ^ { \dom d 0 } & & \relax \cat \XX 1 }
\qquad
\matrixxy
{ \relax \cat \XX 3
\ar [dd] _ { \dom d { 1 2 } } \ar [ddrr] ^ { \dom d { 0 2 } } \\ \\
\relax \cat \XX 2 \ar [ddrr] _ { \dom d 1 } & &
\relax \cat \XX 2 \ar [dd] _ { \dom d 1 } \\ \\
& & \relax \cat \XX 1 }
\qquad
\matrixxy
{ \\ \relax \cat \XX 1
\ar [ddrr] _ { \id { \cat \XX 1 } } \ar [rr] ^ { \dom d { 0 0 } } & &
\relax \cat \XX 2 \ar [dd] ^ { \dom d 1 } \\ \\
& & \relax \cat \XX 1 }
\qquad
\matrixxy
{ \\ \relax \cat \XX 1
\ar [dd] _ { \dom d { 0 0 } } \ar [ddrr] ^ { \id { \cat \XX 1 } } \\ \\
\relax \cat \XX 2 \ar [rr] _ { \dom d 0 } & & \relax \cat \XX 1 }
\end {equation}
(These just say that identities and compositions
all have the proper sources and targets.)
Next, consider the commutative diagram:
\begin {equation} \label {2space d013}
\matrixxy
{ \relax \cat \XX 4
\ar [dd] _ { \dom d { 1 2 3 } } \ar [rr] ^ { \dom d { 0 1 2 } } & &
\relax \cat \XX 3
\ar [dd] ^ { \dom d { 1 2 } } \ar [ddrr] ^ { \dom d { 0 1 } } \\ \\
\relax \cat \XX 3
\ar [ddrr] _ { \dom d { 0 2 } } \ar [rr] _ { \dom d { 0 1 } } & &
\relax \cat \XX 2 \ar [ddrr] ^ { \dom d 0 } & &
\relax \cat \XX 2 \ar [dd] ^ { \dom d 1 } \\ \\
& & \relax \cat \XX 2 \ar [rr] _ { \dom d 0 } & & \relax \cat \XX 1 }
\end {equation}
This is \a pullback cone into the diagram defining \( \cat \XX 3 \),
so the universal property of that pullback
gives \a map \( \cat \XX 4 \maplongto { \dom d { 0 1 3 } } \cat \XX 3 \).
Similarly, the commutative diagram
\begin {equation} \label {2space d023}
\matrixxy
{ \relax \cat \XX 4
\ar [dd] _ { \dom d { 1 2 3 } } \ar [rr] ^ { \dom d { 0 1 2 } } & &
\relax \cat \XX 3
\ar [dd] ^ { \dom d { 1 2 } } \ar [ddrr] ^ { \dom d { 0 2 } } \\ \\
\relax \cat \XX 3
\ar [ddrr] _ { \dom d { 1 2 } } \ar [rr] _ { \dom d { 0 1 } } & &
\relax \cat \XX 2 \ar [ddrr] ^ { \dom d 1 } & &
\relax \cat \XX 2 \ar [dd] ^ { \dom d 1 } \\ \\
& & \relax \cat \XX 2 \ar [rr] _ { \dom d 0 } & & \relax \cat \XX 1 }
\end {equation}
defines \a map \( \cat \XX 4 \maplongto { \dom d { 0 2 3 } } \cat \XX 3 \).
In \a \two space, this diagram must also commute:
\begin {equation} \label {2space ass}
\matrixxy
{ \relax \cat \XX 4
\ar [dd] _ { \dom d { 0 2 3 } } \ar [rr] ^ { \dom d { 0 1 3 } } & &
\relax \cat \XX 3 \ar [dd] ^ { \dom d { 0 2 } } \\ \\
\relax \cat \XX 3 \ar [rr] _ { \dom d { 0 2 } } & & \relax \cat \XX 2 }
\end {equation}
(This is the \emph {associative law} for an category,
since the composable triples in \( \cat \XX 4 \)
are being composed in different ways to give the same result.)
Finally, consider the commutative diagram:
\begin {equation} \label {2space d001}
\matrixxy
{ \relax \cat \XX 2
\ar [ddrr] _ { \id { \cat \XX 2 } } \ar [rr] ^ { \dom d 0 } & &
\relax \cat \XX 1
\ar [ddrr] _ { \id { \cat \XX 1 } } \ar [rr] ^ { \dom d { 0 0 } } & &
\relax \cat \XX 2 \ar [dd] ^ { \dom d 1 } \\ \\
& & \relax \cat \XX 2 \ar [rr] _ { \dom d 0 } & & \relax \cat \XX 1 }
\end {equation}
This too is \a pullback cone
into the pullback diagram defining \( \cat \XX 3 \),
so it defines \a map
\( \cat \XX 2 \maplongto { \dom d { 0 0 1 } } \cat \XX 3 \).
Similarly, the commutative diagram
\begin {equation} \label {2space d011}
\matrixxy
{ \relax \cat \XX 2
\ar [dd] _ { \dom d 1 } \ar [ddrr] ^ { \id { \cat \XX 2 } } \\ \\
\relax \cat \XX 1
\ar [dd] _ { \dom d { 0 0 } } \ar [ddrr] ^ { \id { \cat \XX 1 } } & &
\relax \cat \XX 2 \ar [dd] ^ { \dom d 1 } \\ \\
\relax \cat \XX 2 \ar [rr] _ { \dom d 0 } & & \relax \cat \XX 1 }
\end {equation}
defines \a map \( \cat \XX 2 \maplongto { \dom d { 0 1 1 } } \cat \XX 3 \).
The last requirement of \a \two space
is that both parts of this diagram commute:
\begin {equation} \label {2space unit}
\matrixxy
{ \relax \cat \XX 2
\ar [ddrr] _ { \id { \cat \XX 2 } } \ar [rr] ^ { \dom d { 0 0 1 } } & &
\relax \cat \XX 3 \ar [dd] ^ { \dom d { 0 2 } } \\ \\
& & \relax \cat \XX 2 }
\qquad
\matrixxy
{ \relax \cat \XX 2
\ar [dd] _ { \dom d { 0 1 1 } } \ar [ddrr] ^ { \id { \cat \XX 2 } } \\ \\
\relax \cat \XX 3 \ar [rr] _ { \dom d { 0 2 } } & & \relax \cat \XX 2 }
\end {equation}
(This expresses the \emph {unit laws} of the internal category.)
\par This definition is \a complicated mess of superscripts and subscripts,
but \a simple combinatorial pattern results.
The space \( \cat \XX n \)
is the space of composable \prefix { \lr ( { n - 1 } ) }tuples of arrows---
which is to say that \( n \) points are involved
when you include all of the endpoints.
These points can be numbered (in order) from \( 0 \) to \( n - 1 \).
Then the various \( d \) maps from \( \cat \XX n \) to \( \cat \XX { n ' } \)
are given by listing \( n ' \) numbers from \( 0 \) to \( n - 1 \)
(in order, allowing repetitions and omissions),
indicating which points are retained.
Repeated numbers
indicate an identity arrow
in the composable \prefix { \lr ( { n ' - 1 } ) }tuple,
while skipped numbers in the middle
indicate that some of the original \( n - 1 \) arrows have been composed.
The axioms above are sufficient to yield the natural combinatorial behaviour.
\par For example, consider the final commutative triangle
in the final axiom (\ref {2space unit}).
Following that around the long way,
you start with an arrow \( \dom \xx 0 \longto \dom \xx 1 \)
in \( \cat \XX 2 \),
then repeat position \( 1 \) (introducing an identity arrow)
to get \a composable pair
\( \dom \xx 0 \longto \dom \xx 1 \longto \dom \xx 1 \)
in \( \cat \XX 3 \),
and then skip the middle position (composing across it)
to get an arrow \( \dom \xx 0 \longto \dom \xx 1 \) in \( \cat \XX 2 \) again.
If the right unit law is to hold in this category,
then this must be the very arrow that you started with,
and that is precisely what the commutative triangle guarantees.
\subsection {\Two maps} \label {2map}
Just as maps go between spaces, so \two maps go between \two spaces.
The usual theory of internal categories,
as described in \cite [chapter~8] {Borceux},
uses internal functors,
and early versions of this paper did the same.
However, these are not sufficient,
and \I will need instead \a notion of internal \emph {ana}functors;
anafunctors (in the case of ordinary categories)
were first described in \cite {Makkai}.
(In the category \( \C \) of smooth manifolds and smooth functions,
these \two maps will be \emph {smooth} anafunctors.)
\par Given \two spaces \( \XX \) and \( \YY \),
\a \strong {\two map} \( \xx \) from \( \XX \) to \( \YY \)
consists in part of \a cover
\( \cat { \src \xx } 1 \maplongto { \cat j 1 } \cat \XX 1 \).
Since \( \cat { \src \xx } 1 \) is \a cover,
so is \( \cat { \src \xx } 1 \times \cat { \src \xx } 1 \),
so this pullback exists:
\begin {equation} \label {2mapcover}
\matrixxy
{ \relax \cat { \src \xx } 2
\ar [dd] _ { \pair { \dom d 0 } { \dom d 1 } } \ar [rrr] ^ { \cat j 2 } &&&
\relax \cat \XX 1 \ar [dd] ^ { \pair { \dom d 0 } { \dom d 1 } } \\\\
\relax \cat { \src \xx } 1 \times \cat { \src \xx } 1
\ar [rrr] _ { \cat j 1 \times \cat j 1 } &&&
\relax \cat \XX 0 \times \cat \XX 0 }
\end {equation}
Spaces like \( \cat { \src \xx } 3 \)
and maps like \( \cat { \src \xx } 3
\maplongto { \dom d { 0 2 } } \cat { \src \xx } 2 \)
can all be defined similarly.
\par The \two map \( \xx \) additionally consists of
\a map \( \cat { \src \xx } 1 \maplongto { \cat \xx 1 } \cat \YY 1 \)
and \a map \( \cat { \src \xx } 2 \maplongto { \cat \xx 2 } \cat \YY 2 \).
The following diagrams must commute:
\begin {equation} \label {2mapdom}
\matrixxy
{ \relax \cat { \src \xx } 2
\ar [ddrr] _ { \cat \xx 2 } \ar [rr] ^ { \dom d 0 } &&
\relax \cat { \src \xx } 1 \ar [ddrr] ^ { \cat \xx 1 } \\\\
&& \relax \cat \YY 2 \ar [rr] _ { \dom d 0 } & & \relax \cat \YY 1 }
\qquad
\matrixxy
{ \relax \cat { \src \xx } 2
\ar [dd] _ { \dom d 1 } \ar [ddrr] ^ { \cat \xx 2 } \\\\
\relax \cat { \src \xx } 1 \ar [ddrr] _ { \cat \xx 1 } &&
\relax \cat \YY 2 \ar [dd] ^ { \dom d 1 } \\ \\
&& \relax \cat \YY 1 }
\end {equation}
(These indicate the proper source and target
for an arrow in the range of the \two map.)
\par Now consider this commutative diagram:
\begin {equation} \label {2map2}
\matrixxy
{ \relax \cat { \src \xx } 3
\ar [dd] _ { \dom d { 1 2 } } \ar [rr] ^ { \dom d { 0 1 } } & &
\relax \cat { \src \xx } 2
\ar [dd] ^ { \dom d 1 } \ar [ddrr] ^ { \cat \xx 2 } \\ \\
\relax \cat { \src \xx } 2
\ar [ddrr] _ { \cat \xx 2 } \ar [rr] _ { \dom d 0 } & &
\relax \cat { \src \xx } 1 \ar [ddrr] ^ { \cat \xx 1 } & &
\relax \cat \YY 2 \ar [dd] ^ { \dom d 1 } \\ \\
& & \relax \cat \YY 2 \ar [rr] _ { \dom d 0 } & & \relax \cat \YY 1 }
\end {equation}
This is \a pullback cone into the pullback diagram defining \( \cat \YY 3 \),
so it defines \a map
\( \cat { \src \xx } 3 \maplongto { \cat \xx 3 } \cat \YY 3 \).
The final requirement for \a \two map
is that the following \emph {functoriality} diagrams commute:
\begin {equation} \label {2mapnat}
\matrixxy
{ \relax \cat { \src \xx } 1
\ar [dd] _ { \dom d { 0 0 } } \ar [rr] ^ { \cat \xx 1 } &&
\relax \cat \YY 1 \ar [dd] ^ { \dom d { 0 0 } } \\\\
\relax \cat { \src \xx } 2 \ar [rr] _ { \cat \xx 2 } && \relax \cat \YY 2 }
\qquad
\matrixxy
{ \relax \cat { \src \xx } 3
\ar [dd] _ { \dom d { 0 2 } } \ar [rr] ^ { \cat \xx 3 } &&
\relax \cat \YY 3 \ar [dd] ^ { \dom d { 0 2 } } \\\\
\relax \cat { \src \xx } 2 \ar [rr] _ { \cat \xx 2 } && \relax \cat \YY 2 }
\end {equation}
\par These axioms are enough to ensure
that you can construct arbitrary \( \cat \xx n \)
and that \( \xx \) always commutes with \( d \).
\subsection {Natural transformations} \label {nat}
There is \a new twist with \two spaces.
Not only are there \two maps between \two spaces,
but there is also \a kind of mapping between \two maps.
Since these have no name in ordinary (uncategorified) geometry,
\I borrow the name \q {natural transformation} from pure category theory.
But most properly, these are \emph {smooth} natural tranformations.
\par So, given \two spaces \( \XX \) and \( \YY \),
and given \two maps \( \xx \) and \( \yy \), both from \( \XX \) to \( \YY \),
\a \strong {natural transformation} \( \chi \) from \( \xx \) to \( \yy \)
consists of \a single map
\( \cat { \src \xx } 1 \cap \cat { \src \yy } 1
\maplongto { \src \chi } \cat \YY 2 \).
Of course, there are conditions.
First, these diagrams must commute,
as usual specifying the proper source and target:
\begin {equation} \label {natdom}
\matrixxy
{ \relax \cat { \src \xx } 1 \cap \cat { \src \yy } 1
\ar [dd] _ { \src \chi } \ar [rrr] ^ { \cat { \del j 0 } 1 } &&&
\relax \cat { \src \xx } 1 \ar [dd] ^ { \cat \xx 1 } \\ \\
\relax \cat \YY 2 \ar [rrr] _ { \dom d 0 } &&& \relax \cat \YY 1 }
\qquad
\matrixxy
{ \relax \cat { \src \xx } 1 \cap \cat { \src \yy } 1
\ar [dd] _ { \cat { \del j 1 } 1 } \ar [rrr] ^ { \src \chi } &&&
\relax \cat \YY 2 \ar [dd] ^ { \dom d 1 } \\ \\
\relax \cat { \src \yy } 1 \ar [rrr] _ { \cat \yy 1 } &&& \relax \cat \YY 1 }
\end {equation}
Now the commutative diagram
\begin {equation} \label {nat0[1]}
\matrixxy
{ \relax \cat { \src \xx } 2 \cap \cat { \src \yy } 2
\ar [dd] _ { \cat { \del j 1 } 2 } \ar [rrr] ^ { \dom d 0 } &&&
\relax \cat { \src \xx } 1 \cap \cat { \src \yy } 1
\ar [dd] _ { \cat { \del j 1 } 1 } \ar [ddrr] ^ { \src \chi } \\ \\
\relax \cat { \src \yy } 2
\ar [ddrr] _ { \cat \yy 2 } \ar [rrr] ^ { \dom d 0 } &&&
\relax \cat { \src \yy } 1 \ar [ddrr] ^ { \cat \yy 1 } &&
\relax \cat \YY 2 \ar [dd] ^ { \dom d 1 } \\ \\
&& \relax \cat \YY 2 \ar [rrr] _ { \dom d 0 } &&& \relax \cat \YY 1 }
\end {equation}
defines \a map
\( \cat { \src \xx } 2 \cap \cat { \src \yy } 2
\maplongto { \deldom { \src \chi } 0 1 } \cat \YY 3 \),
and the commutative diagram
\begin {equation} \label {nat[0]1}
\matrixxy
{ \relax \cat { \src \xx } 2 \cap \cat { \src \yy } 2
\ar [dd] _ { \dom d 1 } \ar [rrr] ^ { \cat { \del j 0 } 2 } &&&
\relax \cat { \src \xx } 2
\ar [dd] _ { \dom d 1 } \ar [ddrr] ^ { \cat \xx 2 } \\ \\
\relax \cat { \src \xx } 1 \cap \cat { \src \yy } 1
\ar [ddrr] _ { \src \chi } \ar [rrr] ^ { \cat { \del j 0 } 1 } &&&
\relax \cat { \src \xx } 1 \ar [ddrr] ^ { \cat \xx 1 } &&
\cat \YY 2 \ar [dd] ^ { \dom d 1 } \\ \\
&& \relax \cat \YY 2 \ar [rrr] _ { \dom d 0 } &&& \relax \cat \YY 2 }
\end {equation}
defines \a map
\( \cat { \src \xx } 2 \cap \cat { \src \yy } 2
\maplongto { \domdel { \src \chi } 0 1 } \cat \YY 3 \).
Then the final requirement for \a natural transformation
is that the following diagram commute:
\begin {equation} \label {natnat}
\matrixxy
{ \relax \cat { \src \xx } 2 \cap \cat { \src \yy } 2
\ar [dd] _ { \deldom { \src \chi } 0 1 }
\ar [rrr] ^ { \domdel { \src \chi } 0 1 } &&&
\relax \cat \YY 3 \ar [dd] ^ { \dom d { 0 2 } } \\ \\
\relax \cat \YY 3 \ar [rrr] _ { \dom d { 0 2 } } &&& \relax \cat \YY 2 }
\end {equation}
This condition is the \emph {naturality square};
although because of the internalised presentation,
it's not the same shape as the familiar square with this name.
\begin {prop} \label {2space 2cat}
Using suitable operations,
\two spaces, \two maps, and natural transformations
now form \a bicategory (\a weak \two category) \( \CCC \).
\end {prop}
\begin {prf}
The argument is \a combination
of the argument in \cite {Ehresmann}
that internal categories, (strict) internal functors,
and internal natural transformations
form \a (strict) \two category
and the argument in \cite {Makkai}
that (small) categories, (small) anafunctors, and natural transformations
form \a bicategory.
\par Let me start by proving that, given \two spaces \( \XX \) and \( \YY \),
the \two maps from \( \XX \) to \( \YY \)
and the natural transformations between them
form \a category.
First, given \a \two map \( \XX \maplongto \xx \YY \),
\I want to define
the identity natural transformation \( \iid \xx \) of \( \xx \).
Consider the space \( \nervecat { \src \xx } 1 2 \)
(that is the kernel pair
of \( \cat { \src \xx } 1 \maplongto { \cat j 1 } \cat \XX 1 \),
which should be the source of \( \src { \iid \xx } \)).
This diagram commutes:
\begin {equation} \label {identity nattrans source pullback cone}
\matrixxy
{ \relax \nervecat { \src \xx } 1 2
\ar [dd] _ { \pair { \del { \cat j 1 } 0 } { \del { \cat j 2 } 1 } }
\ar [rr] ^ { \nervecat j 1 2 } &&
\relax \cat \XX 1 \ar [rr] ^ { \dom d { 0 0 } } &&
\relax \cat \XX 2 \ar [dd] ^ { \pair { \dom d 0 } { \dom d 1 } } \\\\
\relax \cat { \src \xx } 1 \times \cat { \src \xx } 1
\ar [rrrr] _ { \cat j 1 \times \cat j 1 } &&&& \cat \XX 1 }
\end {equation}
Since \( \cat { \src \xx } 2 \) is given as \a pullback,
this diagram defines \a map
\( \nervecat { \src \xx } 1 2 \longto \cat { \src \xx } 2 \);
let
\( \nervecat { \src \xx } 1 2 \maplongto { \src { \iid \xx } } \cat \YY 2 \)
be the composite of this map with \( \cat \xx 2 \).
Then the diagrams (\ref {natdom}, \ref {natnat}) easily commute.
\par Now given \two maps \( \xx \), \( \yy \), and \( \zz \),
all from \( \XX \) to \( \YY \),
let \( \chi \) be \a natural transformation from \( \xx \) to \( \yy \),
and let \( \psi \) be \a natural transformation from \( \yy \) to \( \zz \).
To define the composite of \( \chi \) and \( \psi \),
the key is this commutative diagram:
\begin {equation} \label {pulledback composite nattrans}
\matrixxy
{ \relax \cat { \src \xx } 1
\cap \nervecat { \src \yy } 1 2 \cap \cat { \src \zz } 1
\ar [dddddd] _ { \del { \cat j 1 } { 0 2 3 } }
\ar [rrrrrrrrrrrrr] ^ { \del { \cat j 1 } { 0 1 3 } } &&&&&& &&&&& &&
\relax \cat { \src \xx } 1 \cap \cat { \src \yy } 1 \cap \cat { \src \zz } 1
\ar [dd] ^
{ \pair { \del { \cat j 1 } { 0 1 } } { \del { \cat j 1 } { 1 2 } } } \\\\
&&&&&& &&&&& && \relax \cat { \src \xx } 1 \cap \cat { \src \yy } 1
\cap \cat { \src \yy } 1 \cap \cat { \src \zz } 1
\ar [dd] ^ { \src \chi \cap \src \psi } \\\\
&&&&&& &&&&& && \relax \nerve \YY 3
\ar [dd] ^ { \cat { \del j { 0 2 } } 1 } \\\\
\relax \cat { \src \xx } 1 \cap \cat { \src \yy } 1 \cap \cat { \src \zz } 1
\ar [rrrrrr] _
{ \pair { \del { \cat j 1 } { 0 1 } } { \del { \cat j 1 } { 1 2 } } } &&&&&&
\relax \cat { \src \xx } 1 \cap \cat { \src \yy } 1
\cap \cat { \src \yy } 1 \cap \cat { \src \zz } 1
\ar [rrrrr] _ { \src \chi \cap \src \psi } &&&&&
\relax \nerve \YY 3 \ar [rr] _ { \cat { \del j { 0 2 } } 1 } && \nerve \YY 2 }
\end {equation}
This diagram is an example of~(\ref {quotient cone})
for the cover
\( \cat { \src \xx } 1 \cap \cat { \src \yy } 1 \cap \cat { \src \zz } 1
\maplongto { \cat { \del j { 0 2 } } 1 }
\cat { \src \xx } 1 \cap \cat { \src \zz } 1 \),
so (\ref {quotient universality}) defines
the map \( \cat { \src \xx } 1 \cap \cat { \src \zz } 1 \longto \cat \YY 2 \)
that defines the composite natural transformation.
It is straightforward (but tedious)
to check that this composition is associative
and that the identity natural transformations are indeed identities.
\par For the identity \two map \( \id \XX \) on the \two space \( \XX \),
take \( \cat { \src { \id \XX } } 1 : = \cat \XX 1 \),
\( \cat { \src { \id \XX } } 2 : = \cat \XX 2 \),
\( \cat { \id \XX } 1 : = \id { \cat \XX 1 } \),
and \( \cat { \id \XX } 2 : = \id { \cat \XX 2 } \).
All of the conditions to be met are trivial.
\par If \( \XX \maplongto \xx \YY \) and \( \YY \maplongto \yy \ZZ \)
are \two maps,
then \I need to define \a \two map \( \XX \maplongto \zz \ZZ \).
Let \( \cat { \src \zz } 1 \)
be given by the pullback
of \( \cat { \src \xx } 1 \maplongto { \cat \xx 1 } \cat \YY 1 \)
along the cover
\( \cat { \src \yy } 1 \maplongto { \cat j 1 _ \yy } \cat \YY 1 \).
This pullback comes with \a map
\( \cat { \src \zz } 1 \longto \cat { \src \xx } 1 \),
which composes
with \( \cat { \src \xx } 1 \maplongto { \cat j 1 _ \xx } \cat \XX 1 \)
to form the cover
\( \cat { \src \zz } 1 \maplongto { \cat j 1 _ \zz } \cat \XX 1 \);
and it comes with \a map \( \cat { \src \zz } 1 \longto \cat { \src \yy } 1 \),
which composes with
\( \cat { \src \yy } 1 \maplongto { \cat \yy 1 } \cat \ZZ 1 \)
to form \( \cat { \src \zz } 1 \maplongto { \cat \zz 1 } \cat \ZZ 1 \).
Form the pullback \( \cat { \src \zz } 2 \);
the cover \( \cat { \src \zz } 2 \maplongto { \cat j 2 _ \zz } \cat \XX 1 \)
factors through \( \cat { \src \xx } 2 \)
(using the latter's universal property)
to define \a map \( \cat { \src \zz } 2 \longto \cat { \src \xx } 2 \),
which composes with
\( \cat { \src \xx } 2 \maplongto { \cat \xx 2 } \cat \YY 2 \)
to form \a map \( \cat { \src \zz } 2 \longto \cat \YY 2 \).
Now the universal property of \( \cat { \src \yy } 2 \)
defines \a map \( \cat { \src \zz } 2 \longto \cat { \src \yy } 2 \),
which composes
with \( \cat { \src \yy } 2 \maplongto { \cat j 2 _ \yy } \cat \YY 2 \)
to form \( \cat { \src \zz } 2 \maplongto { \cat \zz 2 } \cat \ZZ 2 \).
\par The various coherence conditions in \a (weak) \two category
are now tedious but straightforward to check.
\end {prf}
\par An internal \strong {natural isomorphism}
is simply an internal natural transformation
equipped with an inverse in the sense of this \two category \( \CCC \).
\subsection {Spaces as \two spaces} \label {disc2space}
The category \( \C \) of spaces
may be embedded into the \two category \( \CCC \) of \two spaces.
That is, each space may be interpreted as \a \two space,
and maps between spaces
will be certain \two maps between the corresponding \two spaces.
Specifically, if \( X \) is \a space,
then let \( \cat \XX 1 \) and \( \cat \XX 2 \) each be \( X \),
and let all the \( d \) maps in the definition of \two space
be the identity on \( X \);
this defines \a \two space \( \XX \).
To define \a \two map out of this \two space \( \XX \),
you need \a cover \( \cat { \src \xx } 1 \maplongto { \cat j 1 } X \);
but \I will call this simply \( \src \xx \maplongto j X \) now,
because in fact \( \cat { \src \xx } n \)
is simply \( \nerve { \src \xx } n \) in this case.
So if \( \YY \) is \a \two space,
then \a \two map \( \XX \maplongto \xx \YY \)
consists of \a cover \( \src \xx \maplongto j X \),
\a map \( \src \xx \maplongto { \cat \xx 1 } \cat \YY 1 \)
and \a map \( \nerve { \src \xx } 2 \maplongto { \cat \xx 2 } \cat \YY 2 \),
such that these diagrams commute:
\begin {equation} \label {2discmap}
\matrixxy
{ \relax \nerve { \src \xx } 2
\ar [ddrr] _ { \cat \xx 2 } \ar [rr] ^ { \del j 0 } &&
\src \xx \ar [ddrr] ^ { \cat \xx 1 } \\\\
&& \relax \cat \YY 2 \ar [rr] _ { \dom d 0 } & & \relax \cat \YY 1 }
\qquad
\matrixxy
{ \relax \nerve { \src \xx } 2
\ar [dd] _ { \del j 1 } \ar [ddrr] ^ { \cat \xx 2 } \\\\
\src \xx \ar [ddrr] _ { \cat \xx 1 } &&
\relax \cat \YY 2 \ar [dd] ^ { \dom d 1 } \\ \\
&& \relax \cat \YY 1 }
\qquad
\matrixxy
{ \src \xx
\ar [dd] _ { \del j { 0 0 } } \ar [rr] ^ { \cat \xx 1 } &&
\relax \cat \YY 1 \ar [dd] ^ { \dom d { 0 0 } } \\\\
\relax \nerve { \src \xx } 2 \ar [rr] _ { \cat \xx 2 } && \relax \cat \YY 2 }
\qquad
\matrixxy
{ \relax \nerve { \src \xx } 3
\ar [dd] _ { \del j { 0 2 } } \ar [rr] ^ { \cat \xx 3 } &&
\relax \cat \YY 3 \ar [dd] ^ { \dom d { 0 2 } } \\\\
\relax \nerve { \src \xx } 2 \ar [rr] _ { \cat \xx 2 } && \relax \cat \YY 2 }
\end {equation}
\par Note that if \( \YY \) is derived from \a space \( Y \),
then the map \( \src \xx \maplongto { \cat \xx 1 } Y \)
defines \a map \( X \maplongto x Y \)
(because \( X \) is the quotient of \( \nerve { \src \xx } 2 \));
conversely, any such map defines \a \two map where \( \src \xx : = X \).
Furthermore, if you start with \( \XX \maplongto \xx \YY \),
turn it into \( X \maplongto x Y \),
and turn that back into \a new \( \XX \maplongto { \xx ' } \YY \),
then \( \xx \) and \( \xx ' \) are naturally isomorphic,
with the natural isomorphism \( \chi \)
given by \( \src \chi : = \cat \xx 2 \).
\par Finally, \a group \( G \) defines \a \two space \( \GG \) as follows.
Let \( \cat \GG 1 \) be the trivial space \( 1 \),
and let \( \cat \GG 2 \) be \( G \);
then composition in \( \GG \) is given by multiplication in \( G \).
(If \( G \) is abelian,
then \( \GG \) will also have the structure of \a \two group
---see section~\ref {2group}---,
but not in general.)
\par Comparing the previous paragraphs,
\I get this interesting result:
\begin {prop}
Given \a space \( B \) and \a group \( G \),
interpret \( B \) as \a \two space \( \BB \)
where \( \cat \BB 1 , \cat \BB 2 : = B \),
and interpret \( G \) as \a \two space \( \GG \)
where \( \cat \GG 1 : = 1 \) and \( \cat \GG 2 : = G \).
Then the category of \two maps from \( \BB \) to \( \GG \)
is equivalent to the category of \G transitions on \( B \).
\end {prop}
\begin {prf}
The cover \( \src \xx \) of the \two map is the cover \( U \) of the bundle,
the map \( \cat \xx 1 \) is trivial,
the map \( \cat \xx 2 \) is the transition map \( g \),
the laws~(\ref {2mapdom}) are trivial,
and the functoriality diagrams~(\ref {2mapnat})
are (respectively) the laws \( \eta \) and \( \gamma \)
described in the diagrams~(\ref {transition}).
Given \a pair of \two maps and the corresponding pair of transitions,
the map \( \src \chi \) of \a natural transformation
is the map \( b \) of \a transition morphism;
checking that composition is well behaved as well
is tedious but straightforward.
\end {prf}
Thus the theory of (uncategorified) bundles
is intimately tied up in the general notion of (categorified) \two space;
this is just an interesting aside now,
but it will be important in part~\ref {example},
where \I link the to the notions of gerbe and bundle gerbe.
\subsection {\Two covers in the \two category of \two spaces} \label {222}
If \( \BB \) is \a \two space,
then let \a \strong {\two cover} of \( \BB \)
be \a \two space \( \UU \)
equipped with \a \two map \( \UU \maplongto \jj \BB \)
such that the maps
\( \cat { \src \jj } 1 \maplongto { \cat \jj 1 } \cat \BB 1 \)
and \( \cat { \src \jj } 2 \maplongto { \cat \jj 2 } \cat \BB 2 \)
are both covers in the category~\( \C \).
\par Just as pullbacks along covers always exist,
so \two pullbacks along \two covers always exist.
This is because
you can pullback along the covers \( \cat \jj 1 \) and \( \cat \jj 2 \)
to define the spaces in the \two pullback,
and the universal property of these pullbacks
will give the coherence laws of the \two pullback as \a \two space.
\par In the case of \two quotients,
the map \( \cat \jj 1 \) in the \two cover itself
becomes the cover for \a \two map out of \( B \),
so in this way \two covers have \two quotients as well.
\par In the case that the \two spaces involved are spaces
(\two spaces that are discrete in the category-theoretic sense),
then their \two products, \two pullbacks, and \two quotients in \( \CCC \)
are the same as in~\( \C \) (when they exist).
\section {\Two bundles} \label {2bun}
The most general notion of \two bundle is quite simple;
if \( \BB \) is \a \two space,
then \a \strong {\two bundle} over \( \BB \)
is simply \a \two space \( \EE \)
together with \a \two map \( \EE \maplongto \pp \BB \).
Before long, of course,
\I will restrict attention to \two bundles with more structure than this---
eventually, to locally trivial \two bundles with \a structure \two group.
\par In physical applications where \( \BB \) is spacetime,
the \two space \( \BB \) will be discrete in the category-theoretic sense;
that is, it comes from \a space.
This would simplify some of the following presentation.
\subsection {The \two category of \two bundles} \label {2C/B}
For \a proper notion of equivalence of \two bundles,
\I should define the \two category \( \slice \CCC \BB \)
of \two bundles over \( \BB \).
\par Given \two bundles \( \EE \) and \( \EE ' \) over \( \BB \),
\a \strong {\two bundle morphism} from \( \EE \) to \( \EE ' \)
is \a \two map \( \EE \maplongto \ff \EE ' \)
equipped with \a natural isomorphism \( \pi \):
\begin {equation} \label {2bundle morphism}
\matrixxy
{ \EE \ar [dddr] _ \pp \ar [rr] ^ \ff \ddrriicell \pi & &
\EE ' \ar [dddl] ^ { \pp ' } \\ \\
& & \\ & \BB }
\end {equation}
\par Given \two bundle morphisms
\( \ff \) from \( \EE \) to \( \EE ' \)
and \( \ff ' \) from \( \EE ' \) to \( \EE ' ' \),
the composite \two bundle morphism
is the composite \two map
\( \EE \maplongto \ff \EE ' \maplongto { \ff ' } \EE ' ' \)
equipped with this natural transformation:
\begin {equation} \label {composite 2bundle morphism}
\matrixxy
{ \EE \ar [dddrr] _ \pp \ar [rr] ^ \ff & \relax \ddriicell \pi &
\EE '\ar [ddd] ^ { \pp ' } \ar [rr] ^ { \ff ' } \ddriicell { \pi ' } & &
\EE ' ' \ar [dddll] ^ { \pp ' ' } \\ \\ & & & \\
& & \BB }
\end {equation}
Also, the identity \two bundle morphism on \( \EE \)
is simply the identity \two map on \( \EE \)
equipped with its identity natural transformation.
\par So far, this section is quite analogous to section~\ref {1C/B},
with only \a few natural isomorphisms added in.
But now the categorified case has something distinctly new.
Given \two bundle morphisms \( \ff \) and \( \ff ' \)
both from \( \EE \) to \( \EE ' \),
\a \strong {\two bundle \two morphism} from \( \ff \) to \( \ff ' \)
is \a natural transformation \( \kappa \) from \( \ff \) to \( \ff ' \):
\begin {equation} \label {2bundle 2morphism}
\matrixxy { \EE \rrfulliicell \ff { \ff ' } \kappa & & \EE ' }
\end {equation}
such that the following natural transformation:
\begin {equation} \label {2bundle 2morphism coherence}
\matrixxy
{ \EE
\ar `d ^dr [dddr] _ \pp [dddr] \ar [dr] _ { \ff ' }
\drupperiicell \ff \kappa \\
\relax \driicell { \pi ' } & \EE ' \ar [dd] ^ { \pp ' } \\ & \\ & \BB }
\end {equation}
is equal to \( \pi \).
\begin {prop} \label {2bundle 2cat}
Given \a \two space \( \BB \),
\two bundles over \( \BB \), \two bundle morphisms,
and \two bundle \two morphisms
form \a \two category.
\end {prop}
\begin {prf}
\( \slice \CCC \BB \)
is just \a slice \two category of the \two category \( \CCC \) of \two spaces.
\end {prf}
\par \I now know what it means
for \two bundles \( \EE \) and \( \EE ' \)
to be \strong {equivalent \two bundles}:
weakly isomorphic objects in the \two category \( \slice \CCC \BB \).
That is, there must be
\a \two bundle morphism \( \ff \) from \( \EE \) to \( \EE ' \)
equipped with \a weak inverse,
which is \a \two bundle morphism \( \bar \ff \) from \( \EE ' \) to \( \EE \),
an invertible \two bundle \two morphism \( \kappa \)
from the identity \two bundle morphism on \( \EE \)
to the composite of \( \ff \) and \( \bar \ff \),
and an invertible \two bundle \two morphism \( \bar \kappa \)
from the identity \two bundle morphism on \( \EE ' \)
to the composite of \( \bar \ff \) and \( \ff \).
What this amounts to, then,
are \two maps \( \EE \maplongto \ff \EE ' \)
and \( \EE ' \maplongto { \bar \ff } \EE \)
equipped with natural isomorphisms
\( \pi \), \( \bar \pi \), \( \kappa \), and \( \bar \kappa \):
\begin {equation} \label {2bundle equivalence}
\matrixxy
{ \EE \ar [dddr] _ \pp \ar [rr] ^ \ff \ddrriicell \pi & &
\EE ' \ar [dddl] ^ { \pp ' } \\ \\
& & \\ & \BB }
\qquad
\matrixxy
{ \EE '
\ar [dddr] _ { \pp ' } \ar [rr] ^ { \bar \ff } \ddrriicell { \bar \pi } & &
\EE \ar [dddl] ^ \pp \\ \\
& & \\ & \BB }
\qquad
\matrixxy
{ \EE \ar [drr] _ \ff \ar [rrrr] ^ \EE & \relax \drriicell \kappa & & & \EE \\
& & \EE ' \ar [urr] _ { \bar \ff } & }
\qquad
\matrixxy
{ \EE ' \ar [rrrr] ^ { \EE ' } \ar [drr] _ { \bar \ff } &
\relax \drriicell { \bar \kappa } & & & \EE ' \\
& & \EE \ar [urr] _ \ff & }
\end {equation}
such that the following natural transformations
\begin {equation} \label {2bundle equivalence coherence}
\matrixxy
{ \EE \ar `d ^dr [dddrr] [dddrr] _ \pp \ar [drr] _ \ff \ar [rrrr] ^ \EE &
\relax \drriicell \kappa & & & \EE \ar `d _dl [dddll] [dddll] ^ \pp \\
& \relax \driicell \pi &
\EE ' \ar [dd] ^ { \pp ' } \ar [urr] _ { \bar \ff } \driicell { \bar \pi } & \\
& & & \\ & & \BB & & }
\qquad
\matrixxy
{ \EE '
\ar `d ^dr [dddrr] [dddrr] _ { \pp ' }
\ar [drr] _ { \bar \ff } \ar [rrrr] ^ { \EE ' } &
\relax \drriicell { \bar \kappa } & & &
\EE ' \ar `d _dl [dddll] [dddll] ^ { \pp ' } \\
& \relax \driicell { \bar \pi } &
\EE \ar [dd] ^ \pp \ar [urr] _ \ff \driicell \pi & \\
& & & \\ & & \BB & & }
\end {equation}
are identity natural transformations
on (respectively) \( \EE \maplongto \pp \BB \)
and \( \EE ' \maplongto { \pp ' } \BB \).
In particular, \( \EE \) and \( \EE ' \)
are equivalent as \two spaces.
\subsection {Trivial \two bundles} \label {2trivbun}
If \( \BB \) is \a space and \( \FF \) is \a \two space,
then the Cartesian \two product \( \FF \times \BB \)
is automatically \a \two bundle over \( \BB \).
Just let the \two map be
\( \FF \times \BB \maplongto { \hat \FF \times \BB } \BB \),
the projection onto the right factor in the product.
\par Of this \two bundle, it can rightly be said that \( \FF \) is its fibre.
I'll want to define \a more general notion, however,
of \two bundle over \( \BB \) with fibre \( \FF \).
To begin, I'll generalise this example to the notion of \a trival \two bundle.
Specifically, \a \strong {trivial \two bundle} over \( \BB \)
with fibre \( \FF \)
is simply \a \two bundle over \( \BB \) equipped with
\a \two bundle equivalence to it from \( \FF \times \BB \).
\par In more detail,
this is \a \two space \( \EE \)
equipped with \two maps \( \EE \maplongto \pp \BB \),
\( \FF \times \BB \maplongto \ff \EE \),
and \( \EE \maplongto { \bar \ff } \FF \times \BB \)
equipped with natural isomorphisms
\( \pi \), \( \bar \pi \), \( \kappa \), and \( \bar \kappa \):
\begin {equation} \label {trivial 2bundle}
\matrixxy
{ \FF \times \BB
\ar [dddr] _ { \hat \FF \times \BB } \ar [rr] ^ \ff \ddrriicell \pi & &
\EE \ar [dddl] ^ \pp \\ \\
& & \\ & \BB }
\qquad
\matrixxy
{ \EE \ar [dddr] _ \pp \ar [rr] ^ { \bar \ff } \ddrriicell { \bar \pi } & &
\FF \times \BB \ar [dddl] ^ { \hat \FF \times \BB } \\ \\
& & \\ & \BB }
\qquad
\matrixxy
{ \FF \times \BB \ar [drr] _ \ff \ar [rrrr] ^ { \FF \times \BB } &
\relax \drriicell \kappa & & & \FF \times \BB \\
& & \EE \ar [urr] _ { \bar \ff } & }
\qquad
\matrixxy
{ \EE \ar [rrrr] ^ \EE \ar [drr] _ { \bar \ff } &
\relax \drriicell { \bar \kappa } & & & \EE \\
& & \FF \times \BB \ar [urr] _ \ff & }
\end {equation}
such that the following natural transformations
\begin {equation} \label {trivial 2bundle coherence}
\matrixxy
{ \FF \times \BB
\ar `d ^dr [dddrr] [dddrr] _ { \hat \FF \times \BB }
\ar [drr] _ \ff \ar [rrrr] ^ { \FF \times \BB } &
\relax \drriicell \kappa & & &
\FF \times \BB \ar `d _dl [dddll] [dddll] ^ { \hat \FF \times \BB } \\
& \relax \driicell \pi &
\EE \ar [dd] ^ \pp \ar [urr] _ { \bar \ff } \driicell { \bar \pi } & \\
& & & \\ & & \BB & & }
\qquad
\matrixxy
{ \EE
\ar `d ^dr [dddrr] [dddrr] _ \pp
\ar [drr] _ { \bar \ff } \ar [rrrr] ^ \EE &
\relax \drriicell { \bar \kappa } & & &
\EE \ar `d _dl [dddll] [dddll] ^ \pp \\
& \relax \driicell { \bar \pi } &
\FF \times \BB
\ar [dd] ^ { \hat \FF \times \BB } \ar [urr] _ \ff \driicell \pi & \\
& & & \\ & & \BB & & }
\end {equation}
are identity natural transformations
on (respectively) \( \FF \times \BB \maplongto { \hat \FF \times \BB } \BB \)
and \( \EE \maplongto \pp \BB \).
\subsection {Restrictions of \two bundles} \label {2restrict}
Ultimately, I'll want to deal with \emph {locally} trivial \two bundles,
so \I need \a notion of restriction to \a \two cover.
Analogously to section~\ref {1restrict},
I'll initially model the \two cover
by any \two map \( \UU \maplongto \jj \BB \).
\I call such \a map an (unqualified) \strong {\two subspace} of \( \BB \).
(Thus logically, there is no difference between
\a \two subspace of \( \BB \) and \a \two bundle over \( \BB \);
but one does different things with them,
and they will be refined in different ways.)
\par Given \a \two bundle \( \EE \)
and \a \two subspace \( \UU \) (equipped with the \two map \( \jj \) as above),
\I get \a \two pullback diagram:
\begin {equation} \label {2restriction diagram}
\matrixxy { & & \EE \ar [dd] ^ \pp \\ \\ \UU \ar [rr] _ \jj & & \BB }
\end {equation}
Then the \strong {restriction}
\( \restrict \EE \UU \) of \( \EE \) to \( \UU \)
is \emph {any} \two pullback of this diagram, if any exists.
\I name the associated \two maps and natural isomorphism
as in this commutative diagram:
\begin {equation} \label {2restriction}
\matrixxy
{ \relax \restrict \EE \UU
\ar [dd] _ { \restrict \pp \UU } \ar [rr] ^ { \tilde \jj } \ddrreqcell & &
\EE \ar [dd] ^ \pp \\ \\
\UU \ar [rr] _ \jj & & \BB }
\end {equation}
Notice that \( \restrict \EE \UU \) becomes
both \a \two subspace of \( \EE \) and \a \two bundle over \( \UU \).
\par By the unicity of \two pullbacks,
the restriction \( \restrict \EE \UU \), if it exists,
is well defined up to \two diffeomorphism (equivalence in \( \CCC \));
but is it well defined up to equivalence of \two bundles over \( \UU \)?
The answer is yes,
because the \two pullback cone morphisms
in diagram (\ref {composed 2pullback automorphism})
become \two bundle morphisms when applied to this situation.
\subsection {Locally trivial \two bundles} \label {2loctrivbun}
\I can now combine the preceding ideas
to define \a locally trivial \two bundle with \a given fibre,
which is the general notion of fibre \two bundle
without \a fixed structure \two group.
\par Given \a \two space \( \BB \) and \a \two space \( \FF \),
suppose that \( \BB \) has been supplied
with \a \two cover \( \UU \maplongto \jj \BB \);
this \two cover is \a \two subspace.
Then \a \strong {locally trivial \two bundle}
over \( \BB \) with fibre \( \FF \) subordinate to \( \UU \)
is \a \two bundle \( \EE \) over \( \BB \)
such that the trivial \two bundle \( \FF \times \UU \) over \( \UU \)
is \a restriction \( \restrict \EE \UU \)
(once an appropriate \( \tilde \jj \) has been specified).
\par So to sum up, this consists of the following items:
\begin {closeitemize}
\item \a \strong {base \two space} \( \BB \);
\item \a \strong {cover \two space} \( \UU \);
\item \a \strong {fibre \two space} \( \FF \);
\item \a \strong {total \two space} \( \EE \);
\item \a \strong {cover \two map} \( \UU \maplongto \jj \BB \);
\item \a \strong {projection \two map} \( \EE \maplongto \pp \BB \);
and \item \a \strong {pulled-back \two map}
\( \FF \times \UU \maplongto { \tilde \jj } \EE \);
\end {closeitemize}
equipped with \a \strong {\two pullback} natural isomorphism
\( \tilde \omega \):
\begin {equation} \label {summary 2pullback}
\matrixxy
{ \FF \times \UU \ar [dd] _ { \hat \FF \times \UU } \ar [rr] ^ { \tilde \jj }
\ddrriicell { \tilde \omega } & &
\EE \ar [dd] ^ \pp \\ \\
\UU \ar [rr] _ \jj & & \BB }
\end {equation}
\par Even though \( \omega \) could be ignored in section~\ref {2cover},
the natural transformation \( \tilde \omega \)
incorporates the equivalence
of \( \FF \times \UU \) and \( \restrict \EE \UU \),
so it is an essential ingredient.
(However, since most of the action takes place back in \( \UU \),
it doesn't come up much.)
\par In studying fibre \two bundles,
\I regard \( \BB \) as \a fixed structure on which the \two bundle is defined.
The fibre \( \FF \) (like the \two group \( \GG \) in the next section)
indicates the type of problem that one is considering;
while the \two cover \( \UU \) (together with \( \jj \))
is subsidiary structure that is not preserved by \two bundle morphisms.
Accordingly, \( \EE \) (together with \( \pp \)) \q {is} the \two bundle;
\I may use superscripts or primes on \( \pp \)
if I'm studying more than one \two bundle.
\section {\GGG \two spaces} \label {2Gspace}
Just as \a group can act on \a space,
so \a \two group can act on \a \two space.
\subsection {\Two groups} \label {2group}
\I define here what \cite {HDA5} calls \a \q {coherent} \two group.
For me, they are the default type of \two group.
My notation differs slightly from that of \cite {HDA5}
in order to maintain the conventions indicated earlier.
But my definition is the same as \cite [Definition~19] {HDA5}.
\par First, \a \strong {\two monoid} is \a \two space \( \GG \)
together with \two maps
\( 1 \maplongto \ee \GG \) (the identity, or unit, object)
and \( \GG \times \GG \maplongto \mm \GG \) (multiplication),
as well as natural isomorphisms \( \alpha \), \( \lambda \), and \( \rho \):
\begin {equation} \label {2monoid}
\matrixxy
{ \GG \times \GG \times \GG
\ar [dd] _ { \id \GG \times \mm } \ar [rrr] ^ { \mm \times \id \GG } &
\relax \ddriicell \alpha && \GG \times \GG \ar [dd] ^ \mm \\\\
\GG \times \GG \ar [rrr] _ \mm &&& \GG }
\qquad
\matrixxy
{ \GG \ar [ddrrr] _ { \id \GG } \ar [rrr] ^ { \ee \times \id \GG } &&
\relax \ddriicell \lambda & \GG \times \GG \ar [dd] ^ \mm \\\\
&&& \GG }
\qquad
\matrixxy
{ \GG \ar [ddrrr] _ { \id \GG } \ar [rrr] ^ { \id \GG \times \ee } &&
\relax \ddriicell \rho & \GG \times \GG \ar [dd] ^ \mm \\\\
&&& \GG }
\end {equation}
These isomorphisms
are (respectively) the \emph {associator}
and the \emph {left and right unitors}.
Additionally, in each of the following pairs of diagrams,
the diagram on the left must define the same composite natural transformation
as the diagram on the right:
\begin {equation} \label {pentagon}
\matrixxy
{ \GG \times \GG \times \GG \times \GG
\ar [dd] _ { \id \GG \times \id \GG \times \mm }
\ar [ddrrrr] ^ { \id \GG \times \mm \times \id \GG }
\ar [rrrr] ^ { \mm \times \id \GG \times \id \GG } &&&
\relax \ddrriicell { \alpha \times \idid \GG } &
\GG \times \GG \times \GG \ar [ddrrr] ^ { \mm \times \id \GG } \\
& \relax \ddrriicell { \idid \GG \times \alpha } \\
\GG \times \GG \times \GG \ar [ddrrrr] _ { \id \GG \times \mm } &&&&
\GG \times \GG \times \GG
\ar [dd] _ { \id \GG \times \mm } \ar [rrr] ^ { \mm \times \id \GG } &
\relax \ddriicell \alpha && \GG \times \GG \ar [dd] ^ \mm \\ &&& \\
&&&& \GG \times \GG \ar [rrr] _ \mm &&& \GG }
\quad = \quad
\matrixxy
{ \GG \times \GG \times \GG \times \GG
\ar [dd] _ { \id \GG \times \id \GG \times \mm }
\ar [rrrr] ^ { \mm \times \id \GG \times \id \GG } &
\relax \ddrreqcell &&&
\GG \times \GG \times \GG
\ar [dd] ^ { \id \GG \times \mm } \ar [ddrr] ^ { \mm \times \id \GG } \\
&&&& \relax \ddrriicell \alpha \\
\GG \times \GG \times \GG
\ar [ddrrrr] _ { \id \GG \times \mm } \ar [rrrr] _ { \mm \times \id \GG } &&&
\relax \ddriicell \alpha & \GG \times \GG \ar [ddrr] ^ \mm &&
\GG \times \GG \ar [dd] ^ \mm \\ &&&& && \\
&&&& \GG \times \GG \ar [rr] _ \mm && \GG }
\end {equation}
and:
\begin {equation} \label {triangle}
\matrixxy
{ \GG \times \GG
\ar [ddrrrr] _ { \id \GG \times \id \GG }
\ar [rrrr] ^ { \id \GG \times \ee \times \id \GG } &&&
\relax \ddriicell { \idid \GG \times \lambda } &
\GG \times \GG \times \GG
\ar [dd] ^ { \id \GG \times \mm } \ar [ddrrr] ^ { \mm \times \id \GG } \\
&&&& & \relax \ddriicell \alpha \\
&&&& \GG \times \GG \ar [ddrrr] _ \mm &&& \GG \times \GG \ar [dd] ^ \mm \\
&&&& && \\ &&&& &&& \GG }
\quad = \quad
\matrixxy
{ \GG \times \GG
\ar [ddrrrr] _ { \id \GG \times \id \GG }
\ar [rrrr] ^ { \id \GG \times \ee \times \id \GG } &&&
\relax \ddrriicell { \rho \times \idid \GG } &
\GG \times \GG \times \GG \ar [ddrrr] ^ { \mm \times \id \GG } \\\\
&&&& \GG \times \GG \ar [ddrrr] _ \mm \ar [rrr] ^ { \id \GG \times \id \GG } &&
\relax \ddreqcell &
\GG \times \GG \ar [dd] ^ \mm \\\\
&&&& &&& \GG }
\end {equation}
These coherence requirements
are known (respectively)
as the \emph {pentagon identity} and the \emph {triangle identity};
these names ultimately derive
from the number of natural isomorphisms that appear in each.
\par In terms of generalised elements,
if \( \XX \maplongto \xx \GG \) and \( \XX \maplongto \yy \GG \)
are elements of \( \GG \),
then \I will (again) write \( \mmsub \xx \yy \)
for the composite \( \XX
\maplongto { \pair \xx \yy } \GG \times \GG \maplongto \mm \GG \);
and I'll write \( \eesub \)
for \( \XX \maplongto { \hat \XX } 1 \maplongto \ee \GG \)
(suppressing \( \XX \)).
Then if \( \xx \), \( \yy \), and \( \zz \) are elements of \( \GG \),
then the associator
defines an arrow \( \mmsub { \mmsub \xx \yy } \zz
\maplongTo { \alphasub \xx \yy \zz } \mmsub \xx { \mmsub \yy \zz } \),
and the left and right unitors
define (respectively) arrows \( \mmsub \eesub \xx
\maplongTo { \lambdasub \xx } \xx \)
and \( \mmsub \xx \eesub \maplongTo { \rhosub \xx } \xx \).
The pentagon identity then says
that \( \mmsub { \mmsub { \mmsub \ww \xx } \yy } \zz
\maplongTo { \mmsub { \alphasub \ww \xx \yy } { \iid \zz } }
\mmsub { \mmsub \ww { \mmsub \xx \yy } } \zz
\maplongTo { \alphasub \ww { \mmsub \xx \yy } \zz }
\mmsub \ww { \mmsub { \mmsub \xx \yy } \zz }
\maplongTo { \mmsub { \iid \ww } { \alphasub \xx \yy \zz } }
\mmsub \ww { \mmsub \xx { \mmsub \yy \zz } } \)
is equal to \( \mmsub { \mmsub { \mmsub \ww \xx } \yy } \zz
\maplongTo { \alphasub { \mmsub \ww \xx } \yy \zz }
\mmsub { \mmsub \ww \xx } { \mmsub \yy \zz }
\maplongTo { \alphasub \ww \xx { \mmsub \yy \zz } }
\mmsub \ww { \mmsub \xx { \mmsub \yy \zz } } \);
and the triangle identity
says that \( \mmsub { \mmsub \xx \eesub } \yy
\maplongTo { \alphasub \xx \eesub \yy }
\mmsub \xx { \mmsub \eesub \yy }
\maplongTo { \mmsub { \iid \xx } { \lambdasub \yy } }
\mmsub \xx \yy \)
is equal to \( \mmsub { \rhosub \xx } { \iid \yy } \).
\par \Two monoids have been studied extensively
as \q {weak monoidal categories};
see \cite [chapter~11] {CWM} for \a summary of the noninternalised case.
But by thinking of them as categorifications of monoids,
\cite {HDA5} motivates the following definition to categorify groups:
\A \strong {\two group} is \a \two monoid \( \GG \)
together with \a \two map \( \GG \maplongto \ii \GG \) (the inverse operator),
as well as natural isomorphisms \( \epsilon \) and \( \iota \):
\begin {equation} \label {2inverse}
\matrixxy
{ \GG \ar [dd] _ { \hat \GG } \ar [rrr] ^ { \pair \ii { \id \GG } } &
\relax \ddriicell \epsilon &&
\GG \times \GG \ar [dd] ^ \mm \\\\
1 \ar [rrr] _ \ee &&& \GG }
\qquad
\matrixxy
{ \GG \ar [dd] _ { \pair { \id \GG } \ii } \ar [rrr] ^ { \hat \GG } &
\relax \ddriicell \iota &&
1 \ar [dd] ^ \ee \\\\
\GG \times \GG \ar [rrr] _ \mm &&& \GG }
\end {equation}
These are (respectively) the \emph {left and right invertors},
but they're also called (respectively) the \emph {counit} and \emph {unit},
and this is not merely an analogy---
they in fact form an (internal) adjunction;
compare \cite [Definition~7] {HDA5}.
\par These must also satisfy coherence laws:
\begin {equation} \label {first zig-zag}
\matrixxy
{ \GG
\ar `d ^dr [ddddrr] [ddddrr] _ { \id \GG } \ar [ddrr] ^ { \check \GG }
\ar `r _dr [ddrrrrr] [ddrrrrr] ^ { \triple { \id \GG } \ii { \id \GG } } &&
\relax \ddreqcell \\
\relax \ddrreqcell \\
&& \GG \times \GG
\ar [dd] ^ { \id \GG \times \hat \GG }
\ar [rrr] ^ { \id \GG \times \pair \ii { \id \GG } } &
\relax \ddriicell { \idid \GG \times \epsilon } &&
\GG \times \GG \times \GG
\ar [dd] _ { \id \GG \times \mm } \ar [ddrr] ^ { \mm \times \id \GG } \\
&& &&& \relax \ddrriicell \alpha \\
&& \GG
\ar `dr ^r [ddrrrrr] [ddrrrrr] _ { \id \GG }
\ar [rrr] ^ { \id \GG \times \ee } &&
\relax \ddriicell \rho &
\GG \times \GG \ar [ddrr] ^ \mm && \GG \times \GG \ar [dd] ^ \mm \\
&& &&& && \\
&& &&& && \GG }
\quad = \quad
\matrixxy
{ \GG
\ar `d ^dr [ddddrr] [ddddrr] _ { \id \GG } \ar [ddrr] ^ { \check \GG }
\ar `r _dr [ddrrrrr] [ddrrrrr] ^ { \triple { \id \GG } \ii { \id \GG } } &&
\relax \ddreqcell \\
\relax \ddrreqcell \\
&& \GG \times \GG
\ar [dd] ^ { \hat \GG \times \id \GG }
\ar [rrr] ^ { \pair { \id \GG } \ii \times \id \GG } &
\relax \ddrriicell { \iota \times \idid \GG } &&
\GG \times \GG \times \GG \ar [ddrr] ^ { \mm \times \id \GG } \\
&& && \\
&& \GG
\ar `dr ^r [ddrrrrr] [ddrrrrr] _ { \id \GG }
\ar [rrrrr] ^ { \ee \times \id \GG } &&
\relax \ddrriicell \lambda & &&
\GG \times \GG \ar [dd] ^ \mm \\\\
&& &&& && \GG }
\end {equation}
and:
\begin {equation} \label {second zig-zag}
\matrixxy
{ \GG
\ar `d ^dr [ddddrr] [ddddrr] _ { \id \GG }
\ar [ddrr] ^ { \pair \ii { \id \GG } }
\ar `r _dr [ddrrrrr] [ddrrrrr] ^ { \triple \ii { \id \GG } \ii } &&
\relax \ddreqcell \\
\relax \ddrreqcell \\
&& \GG \times \GG
\ar [dd] ^ { \id \GG \times \hat \GG }
\ar [rrr] ^ { \id \GG \times \pair { \id \GG } \ii } &
\relax \ddriicell { \idid \GG \times \iota } &&
\GG \times \GG \times \GG
\ar [dd] _ { \id \GG \times \mm } \ar [ddrr] ^ { \mm \times \id \GG } \\
&& &&& \relax \ddrriicell \alpha \\
&& \GG
\ar `dr ^r [ddrrrrr] [ddrrrrr] _ { \id \GG }
\ar [rrr] ^ { \id \GG \times \ee } &&
\relax \ddriicell \rho &
\GG \times \GG \ar [ddrr] ^ \mm && \GG \times \GG \ar [dd] ^ \mm \\
&& &&& && \\
&& &&& && \GG }
\quad = \quad
\matrixxy
{ \GG
\ar `d ^dr [ddddrr] [ddddrr] _ { \id \GG }
\ar [ddrr] ^ { \pair { \id \GG } \ii }
\ar `r _dr [ddrrrrr] [ddrrrrr] ^ { \triple \ii { \id \GG } \ii } &&
\relax \ddreqcell \\
\relax \ddrreqcell \\
&& \GG \times \GG
\ar [dd] ^ { \hat \GG \times \id \GG }
\ar [rrr] ^ { \pair \ii { \id \GG } \times \id \GG } &
\relax \ddrriicell { \epsilon \times \idid \GG } &&
\GG \times \GG \times \GG \ar [ddrr] ^ { \mm \times \id \GG } \\
&& && \\
&& \GG
\ar `dr ^r [ddrrrrr] [ddrrrrr] _ { \id \GG }
\ar [rrrrr] ^ { \ee \times \id \GG } &&
\relax \ddrriicell \lambda & &&
\GG \times \GG \ar [dd] ^ \mm \\\\
&& &&& && \GG }
\end {equation}
These laws are known as the \emph {zig-zag identities},
because of their representation in string diagrams;
see section~\ref {2string}.
\par In terms of generalised elements,
if \( \XX \maplongto \xx \GG \) is an element of \( \GG \),
then \I will (again) write \( \iisub \xx \)
for the composite \( \XX \maplongto \xx \GG \maplongto \ii \GG \).
Then the left and right invertors
define (respectively) arrows \( \mmsub { \iisub \xx } \xx
\maplongTo { \epsilonsub \xx } \eesub \)
and \( \mmsub \xx { \iisub \xx } \maplongTo { \iotasub \xx } \eesub \).
The zig-zag identities
state, respecitively, that \( \mmsub { \mmsub \xx { \iisub \xx } } \xx
\maplongTo { \alphasub \xx { \iisub \xx } \xx }
\mmsub \xx { \mmsub { \iisub \xx } \xx }
\maplongTo { \mmsub { \iid \xx } { \epsilonsub \xx } }
\mmsub \xx \eesub
\maplongTo { \rhosub \xx } \xx \)
equals \( \mmsub { \mmsub \xx { \iisub \xx } } \xx
\maplongTo { \mmsub { \iotasub \xx } { \iid \xx } } \mmsub \eesub \xx
\maplongTo { \lambdasub \xx } \xx \);
and that \( \mmsub { \mmsub { \iisub \xx } \xx } { \iisub \xx }
\maplongTo { \alphasub { \iisub \xx } \xx { \iisub \xx } }
\mmsub { \iisub \xx } { \mmsub \xx { \iisub \xx } }
\maplongTo { \mmsub { \iid { \iisub \xx } } { \iotasub \xx } }
\mmsub { \iisub \xx } \eesub
\maplongTo { \rhosub { \iisub \xx } } \iisub \xx \)
equals \( \mmsub { \mmsub { \iisub \xx } \xx } { \iisub \xx }
\maplongTo { \mmsub { \epsilonsub \xx } { \iid { \iisub \xx } } }
\mmsub \eesub { \iisub \xx }
\maplongTo { \lambdasub { \iisub \xx } } \iisub \xx \).
\par Note that since maps in the category \( \C \) are smooth functions,
so that \two maps in the \two category \( \CCC \) are smooth functors,
\a \two group in that \two category is automatically \a \emph {Lie \two group}.
At this point, \I could define (internal) homomorphisms and \two homomorphisms
to get the \two category \( \CCC ^ \Grp \) of internal \two groups;
this is done in \cite {HDA5}.
But in the following, I'm only interested in \a single fixed \two group.
For the rest of part~\ref 2, let \( \GG \) be this fixed \two group.
\subsection {String diagrams for \two groups} \label {2string}
The formulas in the previous section were complicated,
and they'll only get worse further on.
To simplify these, I'll make use of string diagrams,
which are introduced in \cite {string} and used heavily in \cite {HDA5}.
The validity of these diagrams
rests once again on the Mac Lane Coherence Theorem,
which states, roughly speaking,
that when manipulating element expressions in any \two monoid,
one may safely assume that \( \alpha \), \( \lambda \), and \( \rho \)
are all identities,
which is the case in \q {strict} \two monoids.
(More precisely, there is \a unique coherent natural isomorphism
between any two functors intepreting an expression
from the free monoidal category on the variable names.)
Thus \I can prove coherence laws for general \two groups
by proving them for these \q {semistrict} \two groups.
(However, the invertors \( \epsilon \) and \( \iota \) are unavoidable!
Although one may prove ---see \cite [Proposition~45] {HDA5}---
that every \two group is equivalent to \a strict \two group
where \( \epsilon \) and \( \iota \) are also identities,
this depends essentially on the axiom of choice
and does not apply internally to the \two category~\( \CCC \).)
\par If \( \XX \maplongto \xx \GG \)
is \a generalised element of the \two group \( \GG \),
then this (or rather, its identity natural isomorphism)
is drawn as in the diagram on the left below.
Multiplication is shown by juxtaposition;
\( \mmsub \xx \yy \) is in the middle below.
The identity element \( \eesub \) is generally invisible,
but it may be shown as on the right below to stress its presence.
\begin {equation} \label {2monoid maps string}
\matrixxy { \dstring \xx \\\\ & }
\qquad \qquad \qquad
\matrixxy { \dstring \xx & \dstring \yy \\\\ & }
\qquad \qquad \qquad
\matrixxy { & \estring \\\\ & }
\end {equation}
Since string diagrams
use the Mac Lane Coherence Theorem to pretend that all \two monoids are strict,
there is no distinction between
\( \mmsub { \mmsub \xx \yy } \zz \) and \( \mmsub \xx { \mmsub \yy \zz } \),
both of which are shown on the left below.
Similarly, the picture on the right below
may be interpreted as \( \xx \),
\( \mmsub \eesub \xx \), or \( \mmsub \xx \eesub \);
\( \eesub \) is properly invivisible here.
You can also consider these maps to be pictures of the associator and unitors;
these are invisible to string diagrams.
\begin {equation} \label {2monoid string}
\matrixxy { \dstring \xx & \dstring \yy & \dstring \zz \\\\ & & }
\qquad \qquad \qquad
\matrixxy { & \dstring \xx \\\\ & }
\end {equation}
\par String diagrams are boring for \two monoids
(arguably, that's the point),
but the inverse operation is more interesting.
\I will denote \( \iisub \xx \) simply by the diagram at the left,
but the invertors are now interesting:
\begin {equation} \label {2inverse string}
\matrixxy { \ustring \xx \\\\ & }
\qquad \qquad \qquad
\matrixxy
{ & \estring & \\\\ \dstring \xx \rstring \iota && \ustring \xx \\\\ && }
\qquad \qquad \qquad
\matrixxy { \ustring \xx && \dstring \xx \\\\
\rstring \epsilon & \estring & \\\\ & }
\end {equation}
The zig-zag identities are far from invisible,
and you see where they get their name:
\begin {equation} \label {first zig-zag string}
\matrixxy
{ && && \ddstring \xx \\\\ \ddstring \xx \rstring \iota && \ustring \xx \\\\
&& \rstring \epsilon && \\\\ & }
\qquad = \qquad
\matrixxy { \dddstring \xx \\\\ \\\\ \\\\ & }
\end {equation}
and:
\begin {equation} \label {second zig-zag string}
\matrixxy
{ \uustring \xx \\\\ && \dstring \xx \rstring \iota && \uustring \xx \\\\
\rstring \epsilon && \\\\ && && }
\qquad = \qquad
\matrixxy { \uuustring \xx \\\\ \\\\ \\\\ & }
\end {equation}
\par Also keep in mind that \( \epsilon \) and \( \iota \) are isomorphisms,
so their inverses \( \bar \epsilon \) and \( \bar \iota \) also exist;
their string diagrams are simply
upside down the diagrams for \( \epsilon \) and \( \iota \).
These satisfy their own upside down coherence laws;
but more fundamentally, they cancel \( \epsilon \) and \( \iota \).
For example, this coherence law
simply reflects that \( \epsilon \) followed by \( \bar \epsilon \)
is an identity transformation:
\begin {equation} \label {right inverse cancellation string}
\matrixxy
{ \ustring \xx && \dstring \xx \\\\ \rstring \epsilon && \\\\
\ustring \xx \rstring { \bar \epsilon } && \dstring \xx \\\\ && }
\qquad = \qquad
\matrixxy { \uuustring \xx && \dddstring \xx \\\\ \\\\ \\\\ && }
\end {equation}
\par It's worth noticing, in the above examples,
how string diagrams accept natural transformations
whose domains or codomains are given using multiplication;
in general, an arrow (written horizontally)
might take any number of input strings
and produce any number of output strings.
For example, in the following string diagram,
\( \chi \) is an arrow
from \( \msub { \msub \xx { \xx ' } } { \xx ' ' } \)
(or \( \msub \xx { \msub { \xx ' } { \xx ' ' } } \))
to \( \msub \yy { \yy ' } \):
\begin {equation} \label {string 2complicated}
\matrixxy
{ \dstring \xx && \dstring { \xx ' } && \dstring { \xx ' ' } \\\\
\rrstring \chi & \dstring \yy & & \dstring { \yy ' } & \\\\ && && }
\end {equation}
\par As another point,
any equation between string diagrams
(giving an equation between natural transformations)
may be substituted into another string diagram;
this follows from the naturality of the arrows involved.
\subsection {Action on \a \two space} \label {2act}
\I want \( \GG \) to act on the fibre of \a \two bundle.
\( \ii \), \( \epsilon \), and \( \iota \) will play no role in this section,
which could be applied when \( \GG \) is just \a \two monoid
(see \cite [6.1] {Moerdijk1228}).
\par \A \strong {right \GGG \two space}
is \a \two space \( \FF \) equipped with
\a \two map \( \FF \times \GG \maplongto \rr \FF \)
and natural isomorphisms \( \mu \) and \( \upsilon \):
\begin {equation} \label {G2space}
\matrixxy
{ \FF \times \GG \times \GG
\ar [dd] _ { \id \FF \times \mm } \ar [rrr] ^ { \rr \times \id \GG } &
\relax \ddriicell \mu && \FF \times \GG \ar [dd] ^ \rr \\\\
\FF \times \GG \ar [rrr] _ \rr &&& \FF }
\qquad
\matrixxy
{ \FF \ar [ddrrr] _ { \id \FF } \ar [rrr] ^ { \id \FF \times \ee } &&
\relax \ddriicell \upsilon & \FF \times \GG \ar [dd] ^ \rr \\\\
&&& \FF }
\end {equation}
such that the following pairs of diagrams
each define equal natural transformations:
\begin {equation} \label {pentagonaction}
\matrixxy
{ \FF \times \GG \times \GG \times \GG
\ar [dd] _ { \id \FF \times \id \GG \times \mm }
\ar [ddrrrr] ^ { \id \FF \times \mm \times \id \GG }
\ar [rrrr] ^ { \rr \times \id \GG \times \id \GG } &&&
\relax \ddrriicell { \mu \times \idid \GG } &
\FF \times \GG \times \GG \ar [ddrrr] ^ { \rr \times \id \GG } \\
& \relax \ddrriicell { \idid \FF \times \alpha } \\
\FF \times \GG \times \GG \ar [ddrrrr] _ { \id \FF \times \mm } &&&&
\FF \times \GG \times \GG
\ar [dd] _ { \id \FF \times \mm } \ar [rrr] ^ { \rr \times \id \GG } &
\relax \ddriicell \mu && \FF \times \GG \ar [dd] ^ \rr \\ &&& \\
&&&& \FF \times \GG \ar [rrr] _ \rr &&& \FF }
\quad = \quad
\matrixxy
{ \FF \times \GG \times \GG \times \GG
\ar [dd] _ { \id \FF \times \id \GG \times \mm }
\ar [rrrr] ^ { \rr \times \id \GG \times \id \GG } &
\relax \ddrreqcell &&&
\FF \times \GG \times \GG
\ar [dd] ^ { \id \FF \times \mm } \ar [ddrr] ^ { \rr \times \id \GG } \\
&&&& \relax \ddrriicell \mu \\
\FF \times \GG \times \GG
\ar [ddrrrr] _ { \id \FF \times \mm } \ar [rrrr] _ { \rr \times \id \GG } &&&
\relax \ddriicell \mu & \FF \times \GG \ar [ddrr] ^ \rr &&
\FF \times \GG \ar [dd] ^ \rr \\ &&&& && \\
&&&& \FF \times \GG \ar [rr] _ \rr && \FF }
\end {equation}
and:
\begin {equation} \label {triangleaction}
\matrixxy
{ \FF \times \GG
\ar [ddrrrr] _ { \id \FF \times \id \GG }
\ar [rrrr] ^ { \id \FF \times \ee \times \id \GG } &&&
\relax \ddriicell { \idid \FF \times \lambda } &
\FF \times \GG \times \GG
\ar [dd] ^ { \id \FF \times \mm } \ar [ddrrr] ^ { \rr \times \id \GG } \\
&&&& & \relax \ddriicell \mu \\
&&&& \FF \times \GG \ar [ddrrr] _ \rr &&& \FF \times \GG \ar [dd] ^ \rr \\
&&&& && \\
&&&& &&& \FF }
\quad = \quad
\matrixxy
{ \FF \times \GG
\ar [ddrrrr] _ { \id \FF \times \id \GG }
\ar [rrrr] ^ { \id \FF \times \ee \times \id \GG } &&&
\relax \ddrriicell { \upsilon \times \idid \GG } &
\FF \times \GG \times \GG \ar [ddrrr] ^ { \rr \times \id \GG } \\\\
&&&& \FF \times \GG \ar [ddrrr] _ \rr \ar [rrr] ^ { \id \FF \times \id \GG } &&
\relax \ddreqcell &
\FF \times \GG \ar [dd] ^ \rr \\\\
&&&& &&& \FF }
\end {equation}
If more than one \GGG \two space is around at \a time,
then I'll use subscripts or primes on \( \rr \) to keep things straight.
\par In terms of generalised elements,
if \( \XX \maplongto \ww \FF \) is \a generalised element of \( \FF \)
and \( \XX \maplongto \xx \GG \) is \a generalised element of \( \GG \),
then \I will again write \( \rrsub \ww \xx \)
for the composite \( \XX
\maplongto { \pair \ww \xx } \FF \times \GG \maplongto \rr \FF \).
Then if \( \ww \) is an element of \( \FF \)
and \( \xx \) and \( \yy \) are elements of \( \GG \),
then \( \mu \)
defines an arrow \( \rrsub { \rrsub \ww \xx } \yy
\maplongTo { \musub \ww \xx \yy } \rrsub \ww { \mmsub \xx \yy } \),
and \( \upsilon \)
defines an arrow \( \rrsub \ww \eesub
\maplongTo { \upsilonsub \ww } \ww \).
Then the law (\ref {pentagonaction})
says that \( \rrsub { \rrsub { \rrsub \ww \xx } \yy } \zz
\maplongTo { \rrsub { \musub \ww \xx \yy } { \iid \zz } }
\rrsub { \rrsub \ww { \mmsub \xx \yy } } \zz
\maplongTo { \musub \ww { \mmsub \xx \yy } \zz }
\rrsub \ww { \mmsub { \mmsub \xx \yy } \zz }
\maplongTo { \rrsub { \iid \ww } { \alphasub \xx \yy \zz } }
\rrsub \ww { \mmsub \xx { \mmsub \yy \zz } } \)
is equal to \( \rrsub { \rrsub { \rrsub \ww \xx } \yy } \zz
\maplongTo { \musub { \rrsub \ww \xx } \yy \zz }
\rrsub { \rrsub \ww \xx } { \mmsub \yy \zz }
\maplongTo { \musub \ww \xx { \mmsub \yy \zz } }
\rrsub \ww { \mmsub \xx { \mmsub \yy \zz } } \),
for an element \( \ww \) of \( \FF \)
and elements \( \xx \), \( \yy \), and \( \zz \) of \( \GG \);
while the law (\ref {triangleaction})
says that \( \rrsub { \rrsub \ww \eesub } \xx
\maplongTo { \musub \ww \eesub \xx }
\rrsub \ww { \mmsub \eesub \xx }
\maplongTo { \rrsub { \iid \ww } { \lambdasub \xx } }
\rrsub \ww \xx \)
is equal to \( \rrsub { \upsilonsub \ww } { \iid \xx } \).
\par The natural isomorphisms \( \mu \) and \( \upsilon \)
are (respectively) the associator and right unitor of the right action,
while the coherence laws (\ref {pentagonaction}) and~(\ref {triangleaction})
are (respectively) the pentagon and triangle identities for the action,
analogous to~(\ref {pentagon}, \ref {triangle}).
This fits in with the laws for the \two group itself:
\begin {prop} \label {identity G2space}
The \two group \( \GG \) is itself \a right \GGG \two space,
acting on itself by right multiplication.
\end {prop}
\begin {prf}
Set \( \rr _ \GG : = \mm \), \( \mu _ \GG : = \alpha \),
and \( \upsilon _ \GG : = \rho \).
The requirements above
then reduce to (some of) the requirements for \a \two group.
\end {prf}
I'll use this to define \GGG \two torsors and principal \GGG \two bundles.
\par String diagrams also work with \two group actions.
For \a right action, where elements of \( \FF \) are written on the left,
\I write the leftmost string with \a double line to indicate \a boundary;
there is nothing more futher to the left.
For example, this string diagram
describes both \( \mmsub { \mmsub \ww \xx } \yy \)
and \( \mmsub \ww { \mmsub \xx \yy } \):
\begin {equation} \label {action example}
\matrixxy { \dStringl \ww && \dstring \xx && \dstring \yy \\\\ && && && }
\end {equation}
You can also think of this as \a picture of \( \musub \ww \xx \yy \);
\( \mu \) and \( \upsilon \) (like \( \alpha \), \( \lambda \), and \( \rho \))
are invisible to string diagrams.
(This again relies on the Mac Lane Coherence Theorem.)
\subsection {The \two category of \GGG \two spaces} \label {2C^G}
Just as \two spaces form \a \two category \( \CCC \),
so \GGG \two spaces form \a \two category \( \CCC ^ \GG \).
\par Given right \GGG \two spaces \( \FF \) and \( \FF ' \),
\a \strong {\GGG \two map} from \( \FF \) to \( \FF ' \)
is \a \two map \( \FF \maplongto \tt \FF ' \)
together with \a natural isomorphism \( \phi \):
\begin {equation} \label {G2map}
\matrixxy
{ \FF \times \GG \ar [dd] _ \rr \ar [rrr] ^ { \tt \times \id \GG } &
\relax \ddriicell \phi & & \FF ' \times \GG \ar [dd] ^ { \rr ' } \\ \\
\FF \ar [rrr] _ \tt & & & \FF ' }
\end {equation}
Note that \a composition of \GGG \two maps is \a \GGG \two map,
using the following composite natural transformation:
\begin {equation} \label {composite G2map}
\matrixxy
{ \FF \times \GG \ar [dd] _ \rr \ar [rrr] ^ { \tt \times \id \GG } &
\relax \ddriicell \phi & &
\FF ' \times \GG \ar [dd] _ { \rr ' } \ar [rrr] ^ { \tt ' \times \GG } &
\relax \ddriicell { \phi ' } & &
\FF ' ' \times \GG \ar [dd] ^ { \rr ' ' } \\ \\
\FF \ar [rrr] _ \tt & & & \FF ' \ar [rrr] _ { \tt ' } & & & \FF ' ' }
\end {equation}
The identity map on \( \FF \) is also \a \GGG \two map.
Indeed, \GGG \two spaces and \GGG \two maps form \a category \( \CCC ^ \GG \).
\par Given right \GGG \two spaces \( \FF \) and \( \FF ' \)
and \GGG \two maps \( \tt \) and \( \tt ' \) from \( \FF \) to \( \FF ' \),
\a \strong {\GGG natural transformation} from \( \tt \) to \( \tt ' \)
is \a natural transformation \( \tau \):
\begin {equation} \label {Gnattrans}
\matrixxy { \FF \rrfulliicell \tt { \tt ' } \tau & & \FF ' }
\end {equation}
such that the following pair of diagrams defines equal natural transformations:
\begin {equation} \label {Gnattrans coherence}
\matrixxy
{ \FF \times \GG
\ar [dd] _ \rr \ar [rr] ^ { \tt \times \id \GG } \ddrriicell \phi & &
\FF ' \times \GG \ar [dd] ^ { \rr ' } \\ \\
\FF \ar [rr] ^ \tt \rrloweriicell { \tt ' } \tau & & \FF ' }
\quad = \quad
\matrixxy
{ \FF \times \GG
\ar [dd] _ \rr \ar [rr] _ { \tt ' \times \GG }
\rrupperiicell { \tt \times \id \GG } { \tau \times \idid \GG }
\ddrriicell { \phi ' } & &
\FF ' \times \GG \ar [dd] ^ { \rr ' } \\ \\
\FF \ar [rr] _ { \tt ' } & & \FF ' }
\end {equation}
Now \( \CCC ^ \GG \) is in fact \a \two category.
\par Now \I know what it means
for \GGG \two spaces \( \FF \) and \( \FF ' \) to be \strong {equivalent}:
There must be maps
\( \FF \maplongto \tt \FF ' \) and \( \FF ' \maplongto { \bar \tt } \FF \)
equipped with natural isomorphisms:
\begin {eqnarray} \label {G2space equivalence}
&
\matrixxy
{ \FF \times \GG \ar [dd] _ \rr \ar [rrr] ^ { \tt \times \id \GG } &
\relax \ddriicell \phi & & \FF ' \times \GG \ar [dd] ^ { \rr ' } \\ \\
\FF \ar [rrr] _ \tt & & & \FF ' }
\qquad
\matrixxy
{ \FF ' \times \GG \ar [dd] _ { \rr ' } \ar [rrr] ^ { \GG \times \bar \tt } &
\relax \ddriicell { \bar \phi } & & \FF \times \GG \ar [dd] ^ \rr \\ \\
\FF ' \ar [rrr] _ { \bar \tt } & & & \FF }
& \nonumber \\ \\ &
\matrixxy
{ \FF \ar [drr] _ \tt \ar [rrrr] ^ \FF & \relax \drriicell \tau & & & \FF \\
& & \FF ' \ar [urr] _ { \bar \tt } & }
\qquad
\matrixxy
{ \FF ' \ar [drr] _ { \bar \tt } \ar [rrrr] ^ { \FF ' } &
\relax \drriicell { \bar \tau } & & & \FF ' \\
& & \FF \ar [urr] _ \tt & }
& \nonumber
\end {eqnarray}
such that each pair of diagrams below defines equal natural transformations:
\begin {equation} \label {G2space equivalence coherence left}
\matrixxy
{ \FF \times \GG \ar [dd] _ \rr \ar [rrrrrr] ^ { \FF \times \GG } & &
\relax \ddrreqcell & & & & \FF \times \GG \ar [dd] ^ \rr \\ \\
\FF \ar [drrr] _ \tt \ar [rrrrrr] ^ \FF & & \relax \drriicell \tau & & & &
\FF \\
& & & \FF ' \ar [urrr] _ { \bar \tt } & & }
\quad = \quad
\matrixxy
{ \FF \times \GG
\ar [dd] _ \rr \ar [drrr] _ { \tt \times \id \GG }
\ar [rrrrrr] ^ { \FF \times \GG } & &
\relax \drriicell { \tau \times \idid \GG } & & & &
\FF \times \GG \ar [dd] ^ \rr \\
& \relax \driicell \phi & &
\FF ' \times \GG \ar [dd] _ { \rr ' } \ar [urrr] _ { \GG \times \bar \tt } &
\relax \driicell { \bar \phi } \\
\FF \ar [drrr] _ \tt & & & & & & \FF \\ & & & \FF ' \ar [urrr] _ { \bar \tt } }
\end {equation}
and:
\begin {equation} \label {G2space equivalence coherence right}
\matrixxy
{ \FF ' \times \GG \ar [dd] _ { \rr ' } \ar [rrrrrr] ^ { \FF ' \times \GG } & &
\relax \ddrreqcell & & & & \FF ' \times \GG \ar [dd] ^ { \rr ' } \\ \\
\FF ' \ar [drrr] _ { \bar \tt } \ar [rrrrrr] ^ { \FF ' } & &
\relax \drriicell { \bar \tau } & & & & \FF ' \\
& & & \FF \ar [urrr] _ \tt & & }
\quad = \quad
\matrixxy
{ \FF ' \times \GG
\ar [dd] _ { \rr ' } \ar [drrr] _ { \GG \times \bar \tt }
\ar [rrrrrr] ^ { \FF ' \times \GG } & &
\relax \drriicell { \GG \times \bar \tau } & & & &
\FF ' \times \GG \ar [dd] ^ { \rr ' } \\
& \relax \driicell { \bar \phi } & &
\FF \times \GG \ar [dd] _ \rr \ar [urrr] _ { \tt \times \id \GG } &
\relax \driicell { \phi } \\
\FF ' \ar [drrr] _ { \bar \tt } & & & & & & \FF ' \\
& & & \FF \ar [urrr] _ \tt }
\end {equation}
\subsection {\GGG \two torsors} \label {2torsor}
When \I define principal \two bundles,
\I will want the fibre of the \two bundle to be \( \GG \).
However, it should be good enough
if the fibre is only \emph {equivalent} to \( \GG \).
So, let \a \strong {right \GGG \two torsor} be any right \GGG \two space
that is equivalent, as \a \GGG \two space, to \( \GG \) itself.
\par In more detail, this is \a \two space \( \FF \)
equipped with \two maps \( \FF \times \GG \maplongto \rr \FF \),
\( \FF \maplongto \tt \GG \), and \( \GG \maplongto { \bar \tt } \FF \),
and natural isomorphisms:
\begin {eqnarray} \label {G2torsor}
&
\matrixxy
{ \FF \times \GG \times \GG
\ar [dd] _ { \rr \times \id \GG } \ar [rrr] ^ { \id \FF \times \mm } &
\relax \ddriicell \mu & & \FF \times \GG \ar [dd] ^ \rr \\ \\
\FF \times \GG \ar [rrr] _ \rr & & & \FF }
\qquad
\matrixxy
{ \FF \ar [ddrr] _ \FF \ar [rr] ^ { \id \FF \times \ee } &
\relax \driicell \upsilon & \FF \times \GG \ar [dd] ^ \rr \\ & & \\
& & \FF }
\qquad
\matrixxy
{ \FF \times \GG \ar [dd] _ \rr \ar [rrr] ^ { \tt \times \id \GG } &
\relax \ddriicell \phi & & \GG \times \GG \ar [dd] ^ \mm \\ \\
\FF \ar [rrr] _ \tt & & & \GG }
& \nonumber \\ \\ &
\matrixxy
{ \GG \times \GG \ar [dd] _ \mm \ar [rrr] ^ { \GG \times \bar \tt } &
\relax \ddriicell { \bar \phi } & & \FF \times \GG \ar [dd] ^ \rr \\ \\
\GG \ar [rrr] _ { \bar \tt } & & & \FF }
\qquad
\matrixxy
{ \FF \ar [rrrr] ^ \FF \ar [drr] _ \tt & \relax \drriicell \tau & & & \FF \\
& & \GG \ar [urr] _ { \bar \tt } & }
\qquad
\matrixxy
{ \GG \ar [rrrr] ^ { \GG } \ar [drr] _ { \bar \tt } &
\relax \drriicell { \bar \tau } & & & \GG \\
& & \FF \ar [urr] _ \tt & }
& \nonumber
\end {eqnarray}
such that the coherence laws
(\ref {pentagonaction}, \ref {triangleaction}),
\begin {equation} \label {G2torsor coherence left}
\matrixxy
{ \FF \times \GG \ar [dd] _ \rr \ar [rrrrrr] ^ { \FF \times \GG } & &
\relax \ddrreqcell & & & & \FF \times \GG \ar [dd] ^ \rr \\ \\
\FF \ar [drrr] _ \tt \ar [rrrrrr] ^ \FF & &
\relax \drriicell \tau & & & & \FF \\
& & & \GG \ar [urrr] _ { \bar \tt } & & }
\quad = \quad
\matrixxy
{ \FF \times \GG
\ar [dd] _ \rr \ar [drrr] _ { \tt \times \id \GG }
\ar [rrrrrr] ^ { \FF \times \GG } & &
\relax \drriicell { \tau \times \idid \GG } & & & &
\FF \times \GG \ar [dd] ^ \rr \\
& \relax \driicell \phi & &
\GG \times \GG \ar [dd] _ \mm \ar [urrr] _ { \GG \times \bar \tt } &
\relax \driicell { \bar \phi } \\
\FF \ar [drrr] _ \tt & & & & & & \FF \\ & & & \GG \ar [urrr] _ { \bar \tt } }
\end {equation}
and
\begin {equation} \label {G2torsor coherence right}
\matrixxy
{ \GG \times \GG \ar [dd] _ \mm \ar [rrrrrr] ^ { \GG \times \GG } & &
\relax \ddrreqcell & & & & \GG \times \GG \ar [dd] ^ \mm \\ \\
\GG \ar [drrr] _ { \bar \tt } \ar [rrrrrr] ^ \GG & &
\relax \drriicell { \bar \tau } & & & & \GG \\
& & & \FF \ar [urrr] _ \tt & & }
\quad = \quad
\matrixxy
{ \GG \times \GG
\ar [dd] _ \mm \ar [drrr] _ { \GG \times \bar \tt }
\ar [rrrrrr] ^ { \GG \times \GG } & &
\relax \drriicell { \GG \times \bar \tau } & & & &
\GG \times \GG \ar [dd] ^ \mm \\
& \relax \driicell { \bar \phi } & &
\FF \times \GG \ar [dd] _ \rr \ar [urrr] _ { \tt \times \id \GG } &
\relax \driicell { \phi } \\
\GG \ar [drrr] _ { \bar \tt } & & & & & & \GG \\ & & & \FF \ar [urrr] _ \tt }
\end {equation}
are all satisfied.
\subsection {Crossed modules} \label {2crosmod}
As \I remarked in the Introduction,
previous work with gerbes~\cite {crosmodgerb}
and bundle gerbes~\cite {crosmodbungerb}
has been generalised to the case of crossed modules,
which are equivalent to strict \two groups,
but not to general \two groups.
In this extra section accordingly,
\I will explain how crossed modules work as \two groups,
including how to apply string diagrams to them.
\par \A \strong {strict \two group} is \a \two group
whose structural \two maps (\( \mm \), \( \ee \), and \( \ii \)) are all strict
and whose structural isomorphisms
(\( \alpha \), \( \lambda \), \( \rho \), \( \epsilon \), and \( \iota \))
are all identity isomorphisms.
(This is an unnatural condition
in that it forces certain diagrams of morphisms in \a \two category
---specifically, the backbones in (\ref {2monoid}, \ref {2inverse})---
to commute \q {on the nose}, rather than merely up to coherent isomorphism.)
\par \A \strong {left crossed module}
consists of the following data and conditions:
\begin {closeitemize}
\item \a \emph {base group} \( H \);
\item an \emph {automorphism group} \( D \);
\item \a group homomorphism \( d \) from \( H \) to \( D \);
and \item \a (left) group action \( l \) of \( D \) on \( H \);
such that \item the homomorphism \( d \)
is equivariant (an intertwiner)
between the action \( l \) of \( D \) on \( H \)
and the action by (left) conjugation of \( D \) on itself;
and \item the action of \( H \) on itself
(given by applying \( d \) and then \( l \))
is equal to (left) conjugation.
\end {closeitemize}
(The final item is the \emph {Peiffer identity}.)
Expanding on this, \a crossed module consists of
groups (as in section~\ref {1group}) \( D \) and \( H \)
and maps \( H \maplongto d D \) and \( D \times H \maplongto l H \)
such that the following diagrams commute:
\begin {equation} \label {crosmod homomorphism}
\matrixxy
{ H \times H \ar [dd] _ { m _ H } \ar [rrr] ^ { d \times d } &&&
D \times D \ar [dd] ^ { m _ D } \\\\
H \ar [rrr] _ d &&& D }
\qquad
\matrixxy
{ 1 \ar [dd] _ { e _ H } \ar [ddr] ^ { e _ D } \\\\ H \ar [r] _ d & D }
\end {equation}
(making \( d \) \a homomorphism),
and then:
\begin {eqnarray} \label {crosmod action}
&
\matrixxy
{ D \times D \times H
\ar [dd] _ { \id D \times l } \ar [rrrr] ^ { m _ D \times \id H } &&&&
D \times H \ar [dd] ^ l \\\\
D \times H \ar [rrrr] _ l &&&& H }
\qquad
\matrixxy
{ H \ar [ddrrr] _ { \id H } \ar [rrr] ^ { e _ D \times H } &&&
D \times H \ar [dd] ^ l \\\\
&&& H }
& \nonumber \\ \\ &
\matrixxy
{ D \times H \times H
\ar [dd] _ { \id D \times m _ H }
\ar [rrrrr] ^ { \pair l { \id D \times \hat H } \times \id H } &&&&&
H \times D \times H \ar [rrr] ^ { \id H \times l } &&&
H \times H \ar [dd] ^ { m _ H } \\\\
D \times H \ar [rrrrrrrr] _ l &&&&&&&& H }
\qquad
\matrixxy
{ D \ar [dd] _ { \id D \times e _ H } \ar [rr] ^ { \hat D } &&
1 \ar [dd] ^ { e _ H } \\\\
D \times H \ar [rr] _ l && H }
& \nonumber
\end {eqnarray}
(making \( l \) \a group action),
and next:
\begin {equation} \label {crosmod equivariance}
\matrixxy
{ D \times H
\ar [dd] _ { \pair l { \id D \times \hat H } }
\ar [rrrrr] ^ { \id D \times d } &&& &&
D \times D \ar [dd] ^ { m _ D } \\\\
H \times D \ar [rrr] _ { d \times \id D } &&&
D \times D \ar [rr] _ { m _ D } && D }
\end {equation}
(giving the equivariance of \( d \))
and finally:
\begin {equation} \label {crosmod Peiffer}
\matrixxy
{ H \times H
\ar [dd] _
{ \triple { d \times \hat H } { \hat H \times \id H } { \id H \times \hat H } }
\ar [ddrrrrr] ^ { m _ H } \\\\
D \times H \times H \ar [rrr] _ { l \times H } &&&
H \times H \ar [rr] _ { m _ H } && H }
\end {equation}
(giving the Peiffer identity).
(I've drawn these all to avoid using the inverse maps \( i _ D , i _ H \),
just in case anybody wants to try generalising to monoids.)
\par For me, crossed modules are just \a way to describe \two groups.
\begin {prop}
\A crossed module defines \a strict \two group.
\end {prop}
\begin {prf}
Let \( \cat \GG 1 \) be \( D \),
and let \( \cat \GG 2 \) be \( H \times D \)
(or if you prefer, \( H \semitimes D \),
because this semidirect product is quite relevant).
The structure maps \( \dom d 0 , \dom d 1 \)
are (respectively)
\( H \times D \maplongto { d \times \id D } D \times D
\maplongto { m _ D } D \)
and \( H \times D \maplongto { \hat H \times \id D } D \);
then \( \cat G 3 \) may be taken to be \( H \times H \times D \),
with the structure maps \( \dom d { 0 1 } \) and \( \dom d { 1 2 } \)
being (respectively)
\( H \times H \times D
\maplongto { \id H \times d \times \id D } H \times D \times D
\maplongto { \id H \times m _ D } H \times D \)
and \( H \times H \times D
\maplongto { \hat H \times \id H \times \id D } H \times D \).
Finally, \( \dom d { 0 0 } \) and \( \dom d { 0 2 } \)
are (respecitively) \( D \maplongto { e _ H \times \id D } H \times D \)
and \( H \times H \times D \maplongto { m _ H \times \id D } H \times D \);
it's easy to check that this makes \( \GG \) into \a \two space.
\par To define
the \two maps \( \mm \), \( \ee \), and \( \ii \) of the \two group \( \GG \),
first let their sources \( \src \mm \), \( \src \ee \), and \( \src \ii \)
be (repsectively)
\( D \times D \maplongto { \id { D \times D } } D \times D \),
\( 1 \maplongto { \id 1 } 1 \),
and \( D \maplongto { \id D } D \);
that is, they are all strict \two maps.
Since \( \cat \GG 1 \) is already \a group \( D \),
let the object part of these \two maps
be the corresponding maps for the group \( D \);
that is, \( \cat \mm 1 \) is \( D \times D \maplongto { m _ D } D \),
\( \cat \ee 1 \) is \( 1 \maplongto { e _ D } D \),
and \( \cat \ii 1 \) is \( D \maplongto { i _ D } D \).
Finally, let \( \cat \mm 2 \)
be \( H \times D \times H \times D
\maplongto { \id H \times \pair l { \id D \times \hat H } \times \id D }
H \times H \times D \times D \maplongto { m _ H \times m _ D } H \times D \),
let \( \cat \ee 2 \)
be \( 1 \maplongto { e _ H \times e _ D } H \times D \),
and let \( \cat \ii 2 \)
be \( H \times D \maplongto { i _ H \times i _ D } H \times D
\maplongto { \pair l { \hat H \times \id D } } H \times D \).
Then it is straightforward to check,
first that \( \mm \), \( \ee \), and \( \ii \) are all \two maps,
then that they make \( \GG \) into \a (strict) \two group.
\end {prf}
\par For the complete equivalence
between crossed modules and strict \two groups,
refer to \cite {crosmod};
for now, \I will consider how to translate
between statements about crossed modules and string diagrams about \two groups.
An element \( X \maplongto y D \) of \( D \)
defines to an element \( \XX \maplongto \yy \GG \) of \( \GG \),
where \( \XX \) and \( \yy \) are built out of \( X \) and \( y \)
as in section~\ref {disc2space}.
Given elements \( X \maplongto y D \) and \( X \maplongto { y ' } D \),
an element \( X \maplongto x H \)
may correspond to an arrow \( \yy \maplongTo \chi \yy ' \) in \( \GG \),
where the map \( \src \chi \)
is given by \( X \maplongto { \pair x { y ' } } H \times D \);
but of course this only really works
if \( X \maplongto { \pair x { y ' } } H \times D
\maplongto { d \times \id D } D \times D \maplongto { m _ D } D \)
is equal to \( y \).
Often in the literature on crossed modules
(and certainly in \cite {crosmodgerb},
which \I will make contact with in part~\ref {example}),
to indicate (in effect)
that \( x \) should be interpreted as an arrow from \( \yy \) to \( \yy ' \),
one reads that \( X \maplongto x H \maplongto d D \)
is equal to \( \msub y { \isub { y ' } } \).
\par Given an equation about elements of \( H \),
this can be turned into an equation about arrows in \( \GG \)
if the sources and targets can be matched properly,
but they don't have to match completely.
For example, suppose that \( y \), \( y ' \), and \( y ' ' \)
are elements of \( D \),
that \( x \), \( x ' \), and \( x ' ' \) are elements of \( H \),
and that these may be interpreted
(because \( d \) takes the appropriate values)
as \( \yy \maplongto \chi \yy ' \),
\( \yy ' \maplongto { \chi ' } \yy ' ' \),
and \( \yy \maplongto { \chi ' ' } \yy ' ' \).
Then the equation \( \msub x { x ' } = x ' ' \) (in \( H \))
may be interpreted using these string diagrams:
\begin {equation} \label {crosmod string}
\matrixxy
{ \dstring \yy \\\\ \dstring { \yy ' } \pstring \chi \\\\
\dstring { \yy ' ' } \pstring { \chi ' } \\\\ & }
\qquad = \qquad
\matrixxy { \dstring \yy \\\\
\dstring { \yy ' ' } \pstring { \chi ' ' } \\\\ & }
\end {equation}
\par But more interesting equations are possible.
Suppose instead that \( x \)
may still be interpreted as \( \yy \maplongto \chi \yy ' \)
but now \( x ' \) may be interpreted
as \( \msub \yy { \yy ' ' }
\maplongto { \chi ' } \msub { \yy ' } { \yy ' ' } \).
Then notice that \( x \)
can also be interpreted
as \( \msub \yy { \yy ' ' } \longTo \msub { \yy ' } { \yy ' ' } \),
because \( \msub y { \isub { y ' } }
= \msub { \msub y { y ' ' } } { \isub { \msub { y ' } { y ' ' } } } \)
in the group \( D \).
Thus, the equation \( x = x ' \) in \( H \)
may be interpreted using these string diagrams:
\begin {equation} \label {crosmod string right}
\matrixxy
{ \dstring \yy && \ddstring { \yy ' ' } \\\\
\dstring { \yy ' } \pstring \chi \\\\ && }
\qquad = \qquad
\matrixxy
{ \dstring \yy && \dstring { \yy ' ' } \\\\
\dstring { \yy ' } \rstring { \chi ' } && \dstring { \yy ' ' } \\\\ && }
\end {equation}
The lesson is that \( \yy ' ' \) can appear to the right of \( \chi \),
without messing things up.
\par On the other hand,
suppose that now \( x ' \) may be interpreted
as \( \msub { \yy ' ' } \yy
\maplongTo { \chi ' } \msub { \yy ' ' } { \yy ' } \).
Now \( x \) may \emph {not} be intepreted similarly,
because \( \msub y { \isub { y ' } }
= \msub { \msub { y ' ' } y } { \isub { \msub { y ' ' } { y ' } } } \)
is not valid in \a group;
instead, the left side
must be conjugated by \( y ' ' \) to produce the right side.
However, the action \( l \) of \( D \) on \( H \) now comes into play;
if \( y ' ' \) is applied to \( x \) to produce an element \( x ' ' \),
then, using (\ref {crosmod equivariance}),
\( x ' ' \) \emph {may} be interpreted
as \( \msub { \yy ' ' } \yy \longTo \msub { \yy ' ' } { \yy ' } \).
Thus, the following string diagrams,
rather than stating \( x = x ' \),
instead say that \( X
\maplongto { \pair { y ' ' } x } D \times H \maplongto l H \)
is equal to \( x ' \).
\begin {equation} \label {crosmod string left}
\matrixxy
{ \ddstring { \yy ' ' } && \dstring \yy \\\\
&& \dstring { \yy ' } \pstring \chi \\\\ && }
\qquad = \qquad
\matrixxy
{ \dstring \yy && \dstring { \yy ' ' } \\\\
\dstring { \yy ' } \rstring { \chi ' } && \dstring { \yy ' ' } \\\\ && }
\end {equation}
\par \I will apply these principles in part~\ref {example}.
\section {\GGG \two bundles} \label {2Gbun}
\I can now put the above ideas together to get the concept of \GGG \two bundle.
\subsection {Definition of \GGG \two bundle} \label {2def}
Suppose that \I have \a \two space \( \BB \)
with \a \two cover \( \UU \maplongto \jj \BB \),
as well as \a \two group \( \GG \).
In these circumstances,
\a \strong {\GGG \two transition} of the \two cover \( \UU \)
is \a \two map \( \nerve \UU 2 \maplongto \gg \GG \)
together with natural isomorphisms \( \gamma \) and \( \eta \):
\begin {equation} \label {2transition}
\matrixxy
{ \relax \nerve \UU 3
\ar [ddrrrr] _ { \del \jj { 0 2 } }
\ar [rrrr] ^ { \pair { \del \jj { 0 1 } } { \del \jj { 1 2 } } } &&&&
\relax \nerve \UU 2 \times \nerve \UU 2 \ar [rrr] ^ { \gg \times \gg }
\relax \ddriicell \gamma &&&
\GG \times \GG \ar [dd] ^ \mm \\\\
&&&& \relax \nerve \UU 2 \ar [rrr] _ \gg &&& \GG }
\qquad
\matrixxy
{ \UU \ar [dd] _ { \del \jj { 0 0 } }
\ar [rr] ^ { \hat \UU } \ddrriicell \eta &&
1 \ar [dd] ^ \ee \\\\
\relax \nerve \UU 2 \ar [rr] _ \gg && \GG }
\end {equation}
that satisfy three coherence laws.
\par To write these coherence laws in \a manner
that fits on the page while remaining legible,
\I will use string diagrams.
As in section~\ref {1Gbun},
if \( \pair \xx \yy \) defines an element of \( \nerve \UU 2 \),
then let \( \ggsub \xx \yy \) be the composite
\( \XX \maplongto { \pair \xx \yy } \nerve \UU 2 \maplongto \gg \GG \).
Then \( \gamma \)
defines an arrow \( \mmsub { \ggsub \xx \yy } { \ggsub \yy \zz }
\maplongTo { \gammasub \xx \yy \zz } \ggsub \xx \zz \)
(for an element \( \triple \xx \yy \zz \) of \( \nerve \UU 3 \)),
while \( \eta \)
defines an arrow \( \eesub \maplongTo { \etasub \xx } \ggsub \xx \xx \)
(for an element \( \XX \maplongto \xx \UU \) of \( \UU \)).
Then the string diagrams for (the identity on) \( \ggsub \xx \yy \),
for \( \gammasub \xx \yy \zz \), and for \( \etasub \xx \)
are (respectively):
\begin {equation} \label {2transition string}
\matrixxy { & \dstring \gg \\ \xx & & \yy \\ & }
\qquad \qquad
\matrixxy
{ \dstring \gg && && \dstring \gg \\ && \yy \\
\rrstring \gamma && \dstring \gg && \\ \xx && && \zz \\ && }
\qquad \qquad
\matrixxy { & \estring \\ \xx \\ & \pstring \eta \dstring \gg \\\\ & }
\end {equation}
(Here I've drawn in the invisible identity string diagram
to clarify the height of the diagram for \( \eta \).)
As in section~\ref {1Gbun}, only some diagrams can be drawn;
there is no diagram for, say,
\( \mmsub { \ggsub \ww \xx } { \ggsub \yy \zz } \)
---nor will there be any need for such \a diagram.
\par In these terms,
the coherence laws for \a \GGG \two transition are:
\begin {equation} \label {2transition ass}
\matrixxy
{ \dstring \gg && && \dstring \gg && \ddstring \gg \\
&& \xx \\ \rrstring \gamma && \dstring \gg && \\ && && & \yy \\
\ww && \rrstring \gamma && \dstring \gg && \\ && && && \zz \\ && && }
\qquad = \qquad
\matrixxy
{ \ddstring \gg && \dstring \gg && && \dstring \gg \\
&& && \yy \\ && \rrstring \gamma && \dstring \gg && \\ & \xx \\
\rrstring \gamma && \dstring \gg && && \zz \\ \ww \\ && }
\end {equation}
and:
\begin {equation} \label {2transition unit}
\matrixxy
{ && && \ddstring \gg \\ & \xx \\ \pstring \eta \dstring \gg \\\\
\rrstring \gamma && \dstring \gg && \\ && && \yy \\ && }
\qquad = \qquad
\matrixxy { & \dddstring \gg \\\\ \\ \xx & & \yy \\ \\\\ & }
\qquad = \qquad
\matrixxy
{ \ddstring \gg \\ && & \yy \\ && && \pstring \eta \dstring \gg \\\\
\rrstring \gamma && \dstring \gg && \\ \xx \\ && }
\end {equation}
\par The analogy between \a \G transition and \a group
is stronger now as an analogy between \a \GGG \two transition and \a group.
Again, \( \gg \) corresponds to the underlying set of the group,
\( \gamma \) corresponds to the operation of multiplication,
and \( \eta \) corresponds to the identity element.
But now the coherence law~(\ref {2transition ass})
corresponds to the associative law,
while the laws~(\ref {2transition unit})
correspond (respectively) to the left and right unit laws.
Law~(\ref {2transition ass}) is essentially the tetrahedron
in the definition of \two dimensional descent datum
in \cite [page~258] {Duskin}.
\par As in part~\ref 1, \( \eta \) may be derived from \( \gamma \).
But in fact, \I can replace any \GGG \two transition \( \gg \)
with an equivalent \GGG \two transition
for which \( \eta \) takes \a particularly simple form;
so \I defer this to section~\ref {2semistrict}
after the notion of equivalence between \two transitions has been defined.
\par If \( \GG \) acts on some (right) \GGG \two space \( \FF \),
and if \( \EE \) is
\a locally trivial \two bundle over \( \BB \) with fibre \( \FF \),
then \I can consider \a natural isomorphism \( \theta \):
\begin {equation} \label {G2bundle}
\matrixxy
{ \FF \times \nerve \UU 2
\ar [ddrrrr] _ { \id \FF \times \del \jj 0 }
\ar [rrrr] ^ { \id \FF \times \pair \gg { \del \jj 1 } } &&&&
\FF \times \GG \times \UU
\ar [rrr] ^ { \rr \times \id \UU } \ddriicell \theta &&&
\FF \times \UU \ar [dd] ^ { \tilde \jj } \\\\
&&&& \FF \times \UU \ar [rrr] _ { \tilde \jj } &&& \EE }
\end {equation}
\par In terms of elements,
if \( \XX \maplongto \ww \FF \) is an element of \( \FF \),
and \( \pair \xx \yy \) defines an element of \( \nerve \UU 2 \),
then \( \theta \)
defines an arrow \( \jjtildesub { \rrsub \ww { \ggsub \xx \yy } } \yy
\maplongTo { \thetasub \ww \xx \yy } \jjtildesub \ww \xx \).
\I will draw the string diagram for \( \thetasub \ww \xx \yy \)
the same way as in section~\ref {1Gbun}:
\begin {equation} \label {G2bundle string}
\matrixxy
{ \ddStringl \ww && \dstring \gg && \dString \\ && & \yy \\
&& \rString && \theta \dString \\ & \xx \\ && && }
\end {equation}
As in part~\ref 1, the bottom half of this diagram
represents \( \jjtildesub \ww \xx \), \a general notation for \( \tilde \jj \);
while the top half
represents \( \jjtildesub { \rrsub \ww { \ggsub \xx \yy } } \yy \),
which cannot be generalised to all valid expressions
(but does include all appropriate expressions that \I need).
\par By definition,
the natural isomorphism \( \theta \) defines \a \strong {\GGG \two bundle}
if the following coherence laws are satisfied:
\begin {equation} \label {G2bundle ass}
\matrixxy
{ \dddStringl \ww && \dstring \gg && && \dstring \gg && \ddString \\
&& && \yy \\ && \rrstring \gamma && \dstring \gg && \\ && && && & \zz \\
&& && \rrString && && \theta \dString \\ && & \xx \\ && && && && }
\qquad = \qquad
\matrixxy
{ \dddStringl \ww && \ddstring \gg && \dstring \gg && \dString \\
&& && & \zz \\ && && \rString && \theta \dString \\ && && \yy \\
&& \rrString && && \theta \dString \\ && & \xx \\ && && && }
\end {equation}
and:
\begin {equation} \label {G2bundle unit}
\matrixxy
{ \dddStringl \ww && && \ddString \\ & \xx \\
&& \pstring \eta \dstring \gg \\\\ && \rString && \theta \dString \\\\ && && }
\qquad = \qquad
\matrixxy { \dddStringl \ww && \dddString \\\\ \\ & \xx \\ \\\\ && }
\end {equation}
\par If you think of \( \theta \) as analogous to \a (left) group action,
then law~(\ref {G2bundle ass}) is analogous to the associative law,
while (\ref {G2bundle unit}) is the corresponding (left only) unit law.
Certainly the string diagrams are similar
to the associative and left unit laws
in (\ref {2transition ass}, \ref {2transition unit})!
\par So in summary,
\a \GGG \two bundle consists of
the data in bullet points in section~\ref {2loctrivbun},
and the data (\( \gg \), \( \gamma \), \( \eta \), \( \theta \))
involved in diagrams (\ref {2transition}, \ref {G2bundle}),
such that the coherence laws
(\ref {pentagonaction}, \ref {triangleaction},
\ref {2transition ass}, \ref {2transition unit},
\ref {G2bundle ass}, \ref {G2bundle unit})
are all satisfied.
\A \strong {principal \GGG \two bundle}
is simply \a \GGG \two bundle whose fibre \( \FF \) is \( \GG \) itself.
\subsection {The \two category of \GGG \two transitions} \label {2G^B}
To define \a morphism of \GGG \two bundles properly,
\I first need the notion of morphism of \GGG \two transitions.
\par Given \two covers \( \UU \) and \( \UU ' \) of \( \BB \),
\a \GGG \two transition \( \gg \) on \( \UU \),
and \a \GGG \two transition \( \gg ' \) on \( \UU ' \),
\a \strong {\GGG \two transition morphism} from \( \gg \) to \( \gg ' \)
is \a \two map \( \UU \cap \UU ' \maplongto \bb \GG \)
together with natural isomorphisms \( \sigma \):
\begin {equation} \label {2transition morphism left}
\matrixxy
{ \relax \nerve \UU 2 \cap \UU '
\ar [ddrrrr] _ { \del \jj { 0 2 } }
\ar [rrrr] ^ { \pair { \del \jj { 0 1 } } { \del \jj { 1 2 } } } &&&&
\relax \nerve \UU 2 \times \UU \cap \UU ' \ar [rrrr] ^ { \gg \times \bb }
\relax \ddriicell \sigma &&&&
\GG \times \GG \ar [dd] ^ \mm \\\\
&&&& \UU \cap \UU ' \ar [rrrr] _ \bb &&&& \GG }
\end {equation}
and \( \delta \):
\begin {equation} \label {2transition morphism right}
\matrixxy
{ \UU \cap \nerve { \UU ' } 2
\ar [ddrrrr] _ { \del \jj { 0 2 } }
\ar [rrrr] ^ { \pair { \del \jj { 0 1 } } { \del \jj { 1 2 } } } &&&&
\UU \cap \UU ' \times \nerve { \UU ' } 2 \ar [rrrr] ^ { \bb \times \gg ' }
\relax \ddriicell \delta &&&&
\GG \times \GG \ar [dd] ^ \mm \\\\
&&&& \UU \cap \UU ' \ar [rrrr] _ \bb &&&& \GG }
\end {equation}
satisfying five coherence laws.
\par If \( \pair \xx { \xx ' } \) defines an element of \( \UU \cap \UU ' \),
then let \( \bbsub \xx { \xx ' } \)
be \( \XX \maplongto { \pair \xx { \xx ' } } \UU \cap \UU '
\maplongto \bb \GG \).
Now given an element \( \triple \xx \yy { \xx ' } \)
of \( \nerve \UU 2 \cap \UU ' \),
\( \sigma \) defines an arrow
\( \mmsub { \ggsub \xx \yy } { \bbsub \yy { \xx ' } }
\maplongTo { \sigmasub \xx \yy { \xx ' } } \bbsub \xx { \xx ' } \);
and given an element \( \triple \xx { \xx ' } { \yy ' } \)
of \( \UU \cap \nerve { \UU ' } 2 \),
\( \delta \) defines an arrow
\( \mmsub { \bbsub \xx { \xx ' } } { \ggpsub { \xx ' } { \yy ' } }
\maplongTo { \deltasub \xx { \xx ' } { \yy ' } } \bbsub \xx { \yy ' } \).
As in section~\ref {1G^B},
\I draw \( \bbsub \xx { \xx ' } \), \( \sigmasub \xx \yy { \xx ' } \),
and \( \deltasub \xx { \xx ' } { \yy ' } \)
as these string diagrams:
\begin {equation} \label {2transition morphism string}
\matrixxy { & \dstring \bb \\ \xx & & \xx ' \\ & }
\qquad \qquad
\matrixxy
{ \dstring \gg && && \dstring \bb \\ && \yy \\
\rrstring \sigma && \dstring \bb && \\ \xx && && \yy ' \\ && }
\qquad \qquad
\matrixxy
{ \dstring \bb && && \dstring { \gg ' } \\ && \xx ' \\
\rrstring \delta && \dstring \bb && \\ \xx && && \yy ' \\ && }
\end {equation}
Then these are the coherence laws required of \a \GGG \two transition morphism:
\begin {equation} \label {2transition morphism ass left}
\matrixxy
{ \dstring \gg && && \dstring \gg && \ddstring \bb \\
&& \yy \\ \rrstring \gamma && \dstring \gg && \\ && && & \zz \\
\xx && \rrstring \sigma && \dstring \bb && \\ && && && \zz ' \\ && && }
\qquad = \qquad
\matrixxy
{ \ddstring \gg && \dstring \gg && && \dstring \bb \\
&& && \zz \\ && \rrstring \sigma && \dstring \bb && \\ & \yy \\
\rrstring \sigma && \dstring \bb && && \zz ' \\ \xx \\ && }
\end {equation}
and:
\begin {equation} \label {2transition morphism unit left}
\matrixxy
{ && && \ddstring \bb \\ & \xx \\ \pstring \eta \dstring \gg \\\\
\rrstring \sigma && \dstring \bb && \\ && && \xx ' \\ && }
\qquad = \qquad
\matrixxy { & \dddstring \bb \\\\ \\ \xx & & \xx ' \\ \\\\ & }
\end {equation}
and then:
\begin {equation} \label {2transition morphism ass right}
\matrixxy
{ \dstring \bb && && \dstring { \gg ' } && \ddstring { \gg ' } \\
&& \xx ' \\ \rrstring \delta && \dstring \bb && \\ && && & \yy ' \\
\xx && \rrstring \delta && \dstring \bb && \\ && && && \zz ' \\ && && }
\qquad = \qquad
\matrixxy
{ \ddstring \bb && \dstring { \gg ' } && && \dstring { \gg ' } \\
&& && \yy ' \\ && \rrstring { \gamma ' } && \dstring { \gg ' } && \\ & \xx ' \\
\rrstring \delta && \dstring \bb && && \zz ' \\ \xx \\ && }
\end {equation}
and next:
\begin {equation} \label {2transition morphism unit right}
\matrixxy { & \dddstring \bb \\\\ \\ \xx & & \xx ' \\ \\\\ & }
\qquad = \qquad
\matrixxy
{ \ddstring \bb \\ && & \xx ' \\
&& && \pstring { \eta ' } \dstring { \gg ' } \\\\
\rrstring \delta && \dstring \bb && \\ \xx \\ && }
\end {equation}
and finally:
\begin {equation} \label {2transition morphism ass bi}
\matrixxy
{ \dstring \gg && && \dstring \bb && \ddstring { \gg ' } \\
&& \yy \\ \rrstring \sigma && \dstring \bb && \\ && && & \yy ' \\
\xx && \rrstring \delta && \dstring \bb && \\ && && && \zz ' \\ && && }
\qquad = \qquad
\matrixxy
{ \ddstring \gg && \dstring \bb && && \dstring { \gg ' } \\
&& && \yy ' \\ && \rrstring \delta && \dstring \bb && \\ & \yy \\
\rrstring \sigma && \dstring \bb && && \zz ' \\ \xx \\ && }
\end {equation}
\par In the analogy
where the \GGG \two transitions \( \gg , \gg ' \) are analogous to groups,
the \two transition morphism \( \bb \)
is again analogous to
\a bihomogeneous space acted on by \( \gg \) and \( \gg ' \);
\( \sigma \) and \( \delta \)
are again analogous (respectively) to the left and right actions.
But now the coherence laws
(\ref {2transition morphism ass left}, \ref {2transition morphism ass right})
are analogous to the associative laws for the actions,
while (\ref {2transition morphism unit left},
\ref {2transition morphism unit right})
are analogous to the unit laws.
Finally, (\ref {2transition morphism ass bi})
is analogous to the associative law relating the two actions:
the law that it doesn't matter
whether you can multiply first on the left or the right
(so here it doesn't matter
whether you apply \( \sigma \) or \( \delta \) first).
\par Again, \a \GGG \two transition \( \gg \)
serves as its own identity \GGG \two transition morphism,
with \( \sigma \) and \( \delta \) both set to \( \gamma \).
\par To see how to compose \GGG \two transition morphisms,
let \( \quadruple \xx { \xx ' } { \yy ' } { \yy ' ' } \)
be an element of \( \UU \cap \nerve { \UU ' } 2 \cap \UU ' ' \).
Using that \( \sigma ' \) is invertible
(with inverse \( \bar \sigma ' \), say),
\I form this string diagram:
\begin {equation} \label {pulledback composite 2transition morphism string}
\matrixxy
{ & \ddstring \bb && && && \dstring { \bb ' } \\
& && & \xx ' & && && & \yy ' ' \\
& && && \dstring { \gg ' } \rrstring { \bar \sigma ' } && &&
\ddstring { \bb ' } \\\\
& \rrstring \delta && \dstring \bb && \\ \xx & && && & \yy ' \\ & && && && && }
\end {equation}
(Here \( \gg , \gg ' , \gg ' ' \) are all \GGG \two transitions,
\( \bb \) is \a \GGG \two transition morphism from \( \gg \) to \( \gg ' \),
and \( \bb ' \)
is \a \GGG \two transition morphism from \( \gg ' \) to \( \gg ' ' \).)
This diagram describes \a natural isomorphism of this form:
\begin {equation} \label {pulledback composite 2transition morphism}
\matrixxy
{ \UU \cap \nerve { \UU ' } 2 \cap \UU ' '
\ar [dddddd] _ { \del \jj { 0 2 3 } }
\ar [rrrrrrrrrrr] ^ { \del \jj { 0 1 3 } } &&&&& &&&& &&
\UU \cap \UU ' \cap \UU ' '
\ar [dd] ^ { \pair { \del \jj { 0 1 } } { \del \jj { 1 2 } } } \\\\
&&&&& \relax \ddriicell \empty
&&&& && \UU \cap \UU ' \times \UU ' \cap \UU ' '
\ar [dd] ^ { \bb \times \bb ' } \\\\
&&&&& &&&& && \GG \times \GG \ar [dd] ^ \mm \\\\
\UU \cap \UU ' \cap \UU ' '
\ar [rrrrr] _ { \pair { \del \jj { 0 1 } } { \del \jj { 1 2 } } } &&&&&
\UU \cap \UU ' \times \UU ' \cap \UU ' ' \ar [rrrr] _ { \bb \times \bb ' } &&&&
\GG \times \GG \ar [rr] _ \mm && \GG }
\end {equation}
This natural isomorphism is an example of~(\ref {2quotient cone})
for the cover \( \UU \cap \UU ' \cap \UU ' '
\maplongto { \del \jj { 0 2 } } \UU \cap \UU ' ' \),
so (\ref {2quotient universality}) defines
\a \two map \( \UU \cap \UU ' ' \maplongto { \bb ; \bb ' } \GG \)
together with \a natural isomorphism \( \beta \):
\begin {equation} \label {composite 2transition morphism universal property}
\matrixxy
{ \UU \cap \UU ' \cap \UU ' '
\ar [ddrrrrr] _ { \del \jj { 0 2 } }
\ar [rrrrr] ^ { \pair { \del \jj { 0 1 } } { \del \jj { 1 2 } } } &&&&&
\UU \cap \UU ' \times \UU ' \cap \UU ' '
\ar [rrrr] ^ { \bb \times \bb ' } \ddriicell \beta &&&&
\GG \times \GG \ar [dd] ^ \mm \\\\
&&&&& \UU \cap \UU ' ' \ar [rrrr] _ { \bb ; \bb ' } &&&& \GG }
\end {equation}
Given an element \( \pair \xx { \xx ' ' } \) of \( \UU \cap \UU ' ' \),
\I will denote its composite with \( \bb ; \bb ' \)
as \( \bbcbbpsub \xx { \xx ' ' } \).
Then given an element \( \triple \xx { \xx ' } { \xx ' ' } \)
of \( \UU \cap \UU ' \cap \UU ' ' \),
the transformation \( \beta \) above
defines an arrow
\( \mmsub { \bbsub \xx { \xx ' } } { \bbpsub { \xx ' } { \xx ' ' } }
\maplongTo { \betasub \xx { \xx ' } { \xx ' ' } }
\bbcbbpsub \xx { \xx ' ' } \).
As \a string diagram:
\begin {equation}
\label {composite 2transition morphism universal property string}
\matrixxy
{ & \dstring \bb && && \dstring { \bb ' } \\ & && \xx ' \\
& \rrstring \beta && \dstring { \bb ; \bb ' } && \\
\xx & && && & \xx ' ' \\ & && }
\end {equation}
\begin {prop} \label {G2transition morphism composite}
This \two map \( \bb ; \bb ' \)
is \a \GGG \two transition morphism from \( \gg \) to \( \gg ' ' \).
\end {prop}
\begin {prf}
Given an element \( \quadruple \xx \yy { \yy ' } { \yy ' ' } \)
of \( \nerve \UU 2 \cap \UU ' \cap \UU ' ' \),
\I construct this string diagram:
\begin {equation}
\label {composite 2transition morphism left action expansion string}
\matrixxy
{ \ddstring \gg && && && \dstring { \bb ; \bb ' } \\ && & \yy \\
&& && \dstring \bb \rrstring { \bar \beta } && && \ddstring { \bb ' } \\\\
\rrstring \sigma && \dstring \bb && && \yy ' \\ && && && && & \yy ' ' \\
\xx && \rrrstring \beta &&& \dstring { \bb ; \bb ' } &&& \\\\ && &&& }
\end {equation}
(Here \( \bar \beta \) is the inverse of \( \beta \).)
This describes \a natural isomorphism:
\begin {equation}
\label {composite 2transition morphism left action expansion}
\matrixxy
{ \relax \nerve \UU 2 \cap \UU ' \cap \UU ' '
\ar [dddddd] _ { \del \jj { 0 1 3 } }
\ar [rrrrr] ^ { \del \jj { 0 1 3 } } &&& &&
\relax \nerve \UU 2 \cap \UU ' '
\ar [dd] ^ { \pair { \del \jj { 0 1 } } { \del \jj { 1 2 } } } \\\\
&& \relax \ddriicell \empty
& && \relax \nerve \UU 2 \times \UU \cap \UU ' '
\ar [dd] ^ { \gg \times \bb ; \bb ' } \\\\
&&& && \GG \times \GG \ar [dd] ^ \mm \\\\
\relax \nerve \UU 2 \cap \UU ' ' \ar [rrr] _ { \del \jj { 0 2 } } &&&
\UU \cap \UU ' ' \ar [rr] _ { \bb ; \bb ' } && \GG }
\end {equation}
Because the \two cover \( \nerve \UU 2 \cap \UU ' \cap \UU ' '
\maplongto { \del \jj { 0 1 3 } } \UU \cap \UU ' \cap \UU ' ' \)
is an epimorphism, this natural isomorphism follows:
\begin {equation} \label {composite 2transition morphism left action}
\matrixxy
{ \relax \nerve \UU 2 \cap \UU ' '
\ar [ddrrrr] _ { \del \jj { 0 2 } }
\ar [rrrr] ^ { \pair { \del \jj { 0 1 } } { \del \jj { 1 2 } } } &&&&
\relax \nerve \UU 2 \times \UU \cap \UU ' '
\ar [rrrr] ^ { \gg ' \times \bb ; \bb ' }
\ddriicell { \sigma ; \bb ' } &&&&
\GG \times \GG \ar [dd] ^ \mm \\\\
&&&& \UU \cap \UU ' ' \ar [rrrr] _ { \bb ; \bb ' } &&&& \GG }
\end {equation}
This is precisely the natural isomorphism \( \sigma \) for \( \bb ; \bb ' \).
\A natural transformation \( \bb ; \delta ' \)
follows by an analogous argument.
\par To reason with \( \sigma ; \bb ' \) with string diagrams,
simply expand the middle diagram in (\ref {2transition morphism string})
---but for \( \bb ; \bb ' \) instead of \( \bb \), of course---
into (\ref {composite 2transition morphism left action expansion string}).
For example, to prove
the coherence law~(\ref {2transition morphism unit left}),
let \( \triple \xx { \xx ' } { \xx ' ' } \)
be an element of \( \UU \cap \UU ' \cap \UU ' ' \).
Then these string diagrams are all equal:
\begin {eqnarray} \label {composite 2transition morphism unit left}
&
\matrixxy
{ && && \ddstring { \bb ; \bb ' } \\ & \xx \\ \pstring \eta \dstring \gg \\\\
\rrstring { \sigma ; \bb ' } && \dstring { \bb ; \bb ' } && \\
&& && \xx ' ' \\ && }
\qquad = \qquad
\matrixxy
{ && && && \dstring { \bb ; \bb ' } \\ \xx \\
\pstring \eta \dstring \gg && &&
\dstring \bb \rrstring { \bar \beta } && && \ddstring { \bb ' } \\\\
\rrstring \sigma && \dstring \bb && && \xx ' \\ && &&& &&& & \xx ' ' \\
&& \rrrstring \beta &&& \dstring { \bb ; \bb ' } &&& \\\\ && &&& }
& \nonumber \\ \\ & = \qquad
\matrixxy
{ & && \dstring { \bb ; \bb ' } \\\\
\xx & \dstring \bb \rrstring { \bar \beta } && && \dstring { \bb ' } \\
& && \xx ' \\
& \rrstring \beta && \dstring { \bb ; \bb ' } && & \xx ' ' \\\\ & && }
\qquad = \qquad
\matrixxy { & \dddstring { \bb ; \bb ' } \\\\ \xx \\\\ & & \xx ' ' \\\\ & }
& \nonumber
\end {eqnarray}
Here \I use, in turn, the formula for \( \sigma ; \bb ' \),
the law~(\ref {2transition morphism unit left}) for \( \bb \)
(the same law that I'm trying to prove for \( \bb ; \bb ' \)),
and that \( \bar \beta \) is the inverse of \( \beta \).
\par Note that this does not \emph {quite} prove what \I want;
this proves the equality of two natural transformations
whose \two space domain is \( \UU \cap \UU ' \cap \UU ' ' \),
not \( \UU \cap \UU ' ' \) as desired.
However, since \( \UU \cap \UU ' \cap \UU ' '
\maplongto { \del \jj { 0 2 } } \UU \cap \UU ' ' \)
is \two epic,
the equation that \I want indeed follows.
\par The other coherence laws may all be derived similarly.
Therefore, \( \bb ; \bb ' \) is \a \GGG \two transition morphism, as desired.
\end {prf}
This \( \bb ; \bb ' \)
is the composite of \( \bb \) and \( \bb ' \)
in the \two category \( \GG ^ \BB \) of \GGG \two transitions in \( \BB \).
\par There is also \a notion of \two morphism
between \GGG \two transition morphisms.
Specifically, if \( \bb \) and \( \bb ' \)
are both \GGG \two transition morphisms from \( \gg \) to \( \gg ' \),
then \a \strong {\GGG \two transition \two morphism}
from \( \bb \) to \( \bb ' \)
is \a natural transformation \( \xi \):
\begin {equation} \label {2transition 2morphism}
\matrixxy
{ \nerve \UU 2 \rrfulliicell \bb { \bb ' } \xi & & \FF ' }
\end {equation}
where \I draw \( \bbsub \xx { \xx ' }
\maplongTo { \xisub \xx { \xx ' } } \bbpsub \xx { \xx ' } \)
in string diagrams as follows:
\begin {equation} \label {2transition 2morphism string}
\matrixxy
{ & \dstring \bb \\\\ \xx & \dstring { \bb ' } \pstring \xi & \xx ' \\\\ & }
\end {equation}
such that these two coherence laws are satisfied:
\begin {equation} \label {2transition 2morphism left}
\matrixxy
{ \dstring \gg && && \dstring \bb \\ && \yy \\
\rrstring \sigma && \dstring \bb && \\\\
\xx && \dstring { \bb ' } \pstring \xi && \yy ' \\\\ && }
\qquad = \qquad
\matrixxy
{ & \ddstring \gg && && \dstring \bb \\\\
& && \yy && \dstring { \bb ' } \pstring \xi \\\\
& \rrstring { \sigma ' } && \dstring { \bb ' } && \\
\xx & && && & \yy ' \\ & && }
\end {equation}
and:
\begin {equation} \label {2transition 2morphism right}
\matrixxy
{ \dstring \bb && && \dstring { \gg ' } \\ && \xx ' \\
\rrstring \delta && \dstring \bb && \\\\
\xx && \dstring { \bb ' } \pstring \xi && \yy ' \\\\ && }
\qquad = \qquad
\matrixxy
{ & \dstring \bb && && \ddstring { \gg ' } \\\\
& \dstring { \bb ' } \pstring \xi && \xx ' \\\\
& \rrstring { \delta ' } && \dstring { \bb ' } && \\
\xx & && && & \yy ' \\ & && }
\end {equation}
\A \GGG \two transition \two morphism like \( \xi \)
is analogous to \a homomorphism of bimodules;
these coherence laws are analogous to preserving left and right multiplication.
\begin {prop} \label {G2transition homcat}
Given \GGG \two transitions \( \gg \) and \( \gg ' \),
the \GGG \two transition morphisms between them
and the \GGG \two transition \two morphisms between those
form \a category.
\end {prop}
\begin {prf}
Composition of \GGG \two transition \two morphisms
is simply given by composition of their underlying natural transformations.
All \I must prove is that these composites (and the identity)
are indeed \GGG \two transition \two morphisms.
In the case of the identity, there is nothing to prove;
in the case of \a composite,
just apply each coherence law for \( \xi \) and \( \xi ' \) separately
to prove that coherence law for their composite.
\end {prf}
In fact, \I have:
\begin {prop} \label {G2transition 2cat}
Given \a \two group \( \GG \) and \a \two space \( \BB \),
the \GGG \two transitions, with their morphisms and \two morphisms,
form \a \two category \( \GG ^ \BB \).
\end {prop}
The various coherence laws of \a \two category
follow directly from the properties of \two quotients.
\par \I now know what it means that two \GGG \two transitions are equivalent:
that they are equivalent in this \two category.
\par When \( \GG \) is \a strict \two group,
it's possible to assume
that the isomorphism \( \eta \) in \a \GGG \two transition is an identity;
that is, every \GGG \two transition
is equivalent to such \a \q {semistrict} \GGG \two transition.
This will be necessary in part~\ref {example},
and \I will discuss it in section~\ref {2semistrict}.
\subsection {The \two category of \GGG \two bundles} \label {2Bun}
To classify \GGG \two bundles,
\I need \a proper notion of equivalence of \GGG \two bundles.
For this, \I should define the \two category
\( \Bun \CCC \GG \BB \FF \)
of \GGG \two bundles over \( \BB \) with fibre \( \FF \).
\par Assume \GGG \two bundles \( \EE \) and \( \EE ' \) over \( \BB \),
both with the given fibre \( \FF \),
and associated with the \GGG \two transitions
\( \gg \) and \( \gg ' \) (respectively).
Then \a \strong {\GGG \two bundle morphism} from \( \EE \) to \( \EE ' \)
is \a \two bundle morphism \( \ff \) from \( \EE \) to \( \EE ' \)
together with \a \GGG \two transition morphism
\( \bb \) from \( \gg \) to \( \gg ' \)
and \a natural isomorphism \( \zeta \):
\begin {equation} \label {G2bundle morphism}
\matrixxy
{ \FF \times \UU \cap \UU '
\ar [dd] _ { \id \FF \times \del \jj 0 }
\ar [rrrr] ^ { \id \FF \times \pair \bb { \del \jj 1 } } &&&&
\FF \times \GG \times \UU ' \ar [dd] ^ { \rr \times \id { \UU ' } } \\
& \relax \ddrriicell \zeta \\
\FF \times \UU \ar [ddrr] _ { \tilde \jj } && &&
\FF \times \UU ' \ar [dd] ^ { \tilde \jj } \\
&& & \\ && \EE \ar [rr] _ \ff && \EE ' }
\end {equation}
satisfying the two coherence laws given below.
\I think of \( \ff \) as \emph {being} the \GGG \two bundle morphism
and say that \( \ff \)
is \strong {associated} with the \two transition morphism \( \bb \).
\par In terms of generalised elements,
if \( \ww \) is an element of \( \FF \)
and \( \pair \xx { \xx ' } \) defines an element of \( \UU \cap \UU ' \),
then \( \zetasub \ww \xx { \xx ' } \)
is an arrow
from \( \jjptildesub { \rrsub \ww { \bbsub \xx { \xx ' } } } { \xx ' } \)
to the composite of \( \jjtildesub \ww \xx \) with \( \ff \).
\I can draw \( \zeta \) with this string diagram:
\begin {equation} \label {G2bundle morphism string}
\matrixxy
{ \ddStringl \ww && \dstring \bb && \dString \\ && & \xx ' \\
&& \rString && \zeta \dString \\ & \xx \\ && && }
\end {equation}
where it is understood that \( \ff \) is applied to the portion in \( \EE \)
so that the whole diagram may be interpreted in \( \EE ' \).
\par In these terms, the coherence laws that \( \zeta \) must satisfy are:
\begin {equation} \label {G2bundle morphism ass left}
\matrixxy
{ \dddStringl \ww && \dstring \gg && && \dstring \bb && \ddString \\
&& && \yy \\ && \rrstring \sigma && \dstring \bb && \\ && && && & \yy ' \\
&& && \rrString && && \zeta \dString \\ && & \xx \\ && && && && }
\qquad = \qquad
\matrixxy
{ \dddStringl \ww && \ddstring \gg && \dstring \bb && \dString \\
&& && & \yy ' \\ && && \rString && \zeta \dString \\ && && \yy \\
&& \rrString && && \theta \dString \\ && & \xx \\ && && && }
\end {equation}
and:
\begin {equation} \label {G2bundle morphism ass right}
\matrixxy
{ \dddStringl \ww && \dstring \bb && && \dstring { \gg ' } && \ddString \\
&& && \xx ' \\ && \rrstring \delta && \dstring \bb && \\ && && && & \yy ' \\
&& && \rrString && && \zeta \dString \\ && & \xx \\ && && && && }
\qquad = \qquad
\matrixxy
{ \dddStringl \ww && \ddstring \bb && \dstring { \gg ' } && \dString \\
&& && & \yy ' \\ && && \rString && \theta ' \dString \\ && && \xx ' \\
&& \rrString && && \zeta \dString \\ && & \xx \\ && && && }
\end {equation}
\par If \( \EE \) and \( \EE ' \) are the same \GGG \two bundle
(so associated with the same \GGG \two transition)
and \( \ff \) is the identity \two map \( \id \EE \),
then \( \zeta \) is \a special case of \( \theta \)
(with \( \bb \) taken to be
the identity \GGG \two transition morphism on \( \gg \),
which is \( \gg \) itself);
both coherence laws for \( \zeta \)
reduce to the associativity law for \( \theta \).
In this way, every \GGG \two bundle
has an identity \GGG \two bundle automorphism.
\par Given \GGG \two bundle morphisms
\( \ff \) from \( \EE \) to \( \EE ' \)
and \( \ff ' \) from \( \EE ' \) to \( \EE ' ' \),
the composite \GGG \two bundle morphism \( \ff ; \ff ' \)
is simply the composite \two bundle morphism
\( \EE \maplongto \ff \EE ' \maplongto { \ff ' } \EE ' ' \)
together with the composite \GGG \two transition \( \bb ; \bb ' \)
as described in section~\ref {2G^B}.
\begin {prop} \label {G2bundle morphism composite}
This \( \ff ; \ff ' \) really is \a \GGG \two bundle morphism.
\end {prop}
\begin {prf}
Let \( \ww \) be an element of \( \FF \),
and let \( \triple \xx { \xx ' } { \xx ' ' } \)
define an element of \( \UU \cap \UU ' \cap \UU ' ' \).
Then \I can form this string diagram:
\begin {equation} \label {G2bundle morphism composite string}
\matrixxy
{ \ddddStringl \ww && && \dstring { \bb ; \bb ' } && && \ddString \\
&& && && & \xx ' ' \\
&& \ddstring \bb \rrstring { \bar \beta } && && \dstring { \bb ' } && \\\\
&& && \xx ' && \rString && \zeta ' \dString \\\\
&& \rrrString && && && \zeta \dString \\ & \xx \\ && && && && }
\end {equation}
Because \( \nerve \UU 3 \maplongto { \del \jj { 0 2 } } \nerve \UU 2 \)
is \two epic,
\I can ignore \( \xx ' \)
(just as \I ignored \( \yy ' \)
in the proof of Proposition~\ref {G2transition morphism composite});
thus this natural isomorphism:
\begin {equation} \label {composite G2bundle morphism}
\matrixxy
{ \FF \times \UU \cap \UU ' '
\ar [dd] _ { \id \FF \times \del \jj 0 }
\ar [rrrrrr] ^ { \id \FF \times \pair { \bb ; \bb ' } { \del \jj 1 } } &&&&&&
\FF \times \GG \times \UU ' ' \ar [dd] ^ { \rr \times \id { \UU ' ' } } \\
&& \relax \ddrricell { \zeta ; \zeta ' } \\
\FF \times \UU \ar [ddrr] _ { \tilde \jj } && && &&
\FF \times \UU ' ' \ar [dd] ^ { \tilde \jj ' ' } \\
&& && \\ && \EE \ar [rr] _ \ff && \EE ' \ar [rr] _ { \ff ' } && \EE ' ' }
\end {equation}
This is simply the diagram (\ref {G2bundle morphism})
for the \GGG \two bundle morphism \( \ff ; \ff ' \).
The coherence laws for \( \zeta ; \zeta ' \)
may easily be proved using string diagrams,
by substituting (\ref {G2bundle morphism composite string})
and using the various coherence laws for \( \bb ; \bb ' \).
Therefore, \( \ff ; \ff ' \) really is \a \GGG \two bundle morphism.
\end {prf}
\par Given \GGG \two bundles \( \EE \) and \( \EE ' \)
and \GGG \two bundle morphisms \( \ff \) and \( \ff ' \),
both now from \( \EE \) to \( \EE ' \),
\a \strong {\GGG \two bundle \two morphism} from \( \ff \) to \( \ff ' \)
consists of \a bundle \two morphism \( \kappa \)
and \a \GGG \two transition \two morphism \( \xi \)
(between the underlying \GGG \two transition morphisms
\( \bb \) and \( \bb ' \))
satisfying this coherence law:
\begin {equation} \label {G2bundle 2morphism}
\matrixxy
{ \FF \times \UU \cap \UU '
\ar [dd] _ { \id \FF \times \del \jj 0 }
\ar [rrrr] ^ { \id \FF \times \pair \bb { \del \jj 1 } } &&&&
\FF \times \GG \times \UU ' \ar [dd] ^ { \rr \times \id { \UU ' } } \\
& \relax \ddrriicell \zeta \\
\FF \times \UU \ar [ddrr] _ { \tilde \jj } && &&
\FF \times \UU ' \ar [dd] ^ { \tilde \jj } \\
&& & \\ && \EE \ar [rr] ^ \ff \rrloweriicell { \ff ' } \kappa && \EE ' }
\qquad = \qquad
\matrixxy
{ \FF \times \UU \cap \UU '
\ar [dd] _ { \id \FF \times \del \jj 0 }
\ar [rrrr] _ { \id \FF \times \pair { \bb ' } { \del \jj 1 } } &
\relax \rrupperiicell { \id F \times \pair \bb { \del \jj 1 } }
{ \idid F \times \pair \xi { \iid { \del \jj 1 } } } &&&
\FF \times \GG \times \UU ' \ar [dd] ^ { \rr \times \id { \UU ' } } \\
& \relax \ddrriicell { \zeta ' } \\
\FF \times \UU \ar [ddrr] _ { \tilde \jj } && &&
\FF \times \UU ' \ar [dd] ^ { \tilde \jj } \\
&& & \\ && \EE \ar [rr] _ { \ff ' } && \EE ' }
\end {equation}
\begin {prop} \label {G2bundle 2cat}
Given \a space \( \BB \),
\GGG bundles over \( \BB \)
and their \GGG bundle morphisms and \GGG bundle \two morphisms
form \a \two category \( \Bun \CCC \GG \BB \FF \).
\end {prop}
\begin {prf}
\Two bundles, \two bundle morphisms, and \two bundle \two morphisms
form \a \two category \( \slice \CCC \BB \);
\GGG \two transitions, \GGG \two transition morphisms,
and \GGG \two transition \two morphisms
form \a \two category \( \GG ^ \BB \).
\( \Bun \CCC \GG \BB \FF \) is \a straightforward combination.
\end {prf}
\par \I now know what it means
for \GGG \two bundles \( \EE \) and \( \EE ' \)
to be \strong {equivalent \GGG \two bundles}:
equivalent objects in the \two category \( \Bun \CCC \GG \BB \FF \).
In particular, \( \EE \) and \( \EE ' \) are equivalent as \two bundles.
\subsection {Associated \GGG \two bundles} \label {2assocbun}
Just as \a \GGG \two bundle may be reconstructed from its \G transition,
so \a \GGG \two bundle may be reconstructed from its \GGG \two transition.
\begin {prop} \label {associated 2bundle}
Given \a \two cover \( \UU \maplongto \jj \BB \),
\a \GGG \two transition \( \nerve \UU 2 \maplongto \gg \GG \),
and \a \GGG \two space \( \FF \),
there is \a \GGG \two bundle \( \EE \) over \( \BB \) with fibre \( \FF \)
associated with the \two transition~\( \gg \).
\end {prop}
\begin {prf}
\I will construct \( \EE \)
as the \two quotient of an equivalence \two relation
from \( \FF \times \nerve \UU 2 \) to \( \FF \times \UU \).
One of the \two maps in the equivalence \two relation
is \( \id \FF \times \del \jj 0 \);
the other
is \( \FF \times \nerve \UU 2
\maplongto { \id \FF \times \pair \gg { \del \jj 1 } }
\FF \times \GG \times \UU
\maplongto { \rr \times \id \UU } \FF \times \UU \).
(You can see these \two maps in diagram~(\ref {G2bundle}),
where they appear with the \two quotient \two map
\( \FF \times \UU \maplongto { \tilde \jj } \EE \)
and natural isomorphism \( \theta \)
that this proof will construct.)
Since \( \id \FF \times \del \jj 0 \) is \a \two cover
and every equivalence \two relation
involving \a \two cover has \a \two quotient,
the desired \two quotient does exist, satisfying (\ref {G2bundle})---
at least, if this really is an equivalence \two relation!
\par To begin with, it is \a \two relation;
that is, the two \two maps are jointly \two monic.
Given two elements
\( \pair \ww { \pair \xx \yy } \)
and \( \pair { \ww ' } { \pair { \xx ' } { \yy ' } } \)
of \( \FF \times \nerve \UU 2 \),
the \two monicity diagrams (\ref {2relation})
give arrows \( \pair \ww \xx \longTo \pair { \ww ' } { \xx ' } \)
and \( \pair { \rrsub \ww { \ggsub \xx \yy } } \yy
\longTo \pair { \rrsub { \ww ' } { \ggsub { \xx ' } { \yy ' } } } { \yy ' } \);
these certainly give arrows
\( \ww \longTo \ww ' \), \( \xx \longTo \xx ' \), and \( \yy \longTo \yy ' \).
Furthermore, these are the only isomorphisms
that will satisfy (\ref {2relation coherence}),
since \( \del \jj 0 \) and \( \del \jj 1 \) are \two epic.
\par The reflexivity \two map of the equivalence \two relation
is \( \FF \times \UU
\maplongto { \id \FF \times \del \jj { 0 0 } } \FF \times \nerve \UU 2 \).
In (\ref {2reflexivity}), \( \del \omega 0 \) is trivial;
\( \del \omega 1 \) is given by this string diagram:
\begin {equation} \label {2bundle construction eq2rel reflexivity}
\matrixxy { \ddStringl \ww \\ && & \xx \\
&& \pstring \eta \dstring \gg \\\\ && }
\end {equation}
\par The \two kernel pair of \( \id \FF \times \del \jj 0 \)
is \( \FF \times \nerve \UU 3 \):
\begin {equation} \label {2bundle construction eq2rel 2kernel pair}
\matrixxy
{ \FF \times \nerve \UU 3
\ar [dd] _ { \id \FF \times \del \jj { 0 1 } }
\ar [rrr] ^ { \id \FF \times \del \jj { 0 2 } } &&&
\FF \times \nerve \UU 2 \ar [dd] ^ { \id \FF \times \del \jj 0 } \\\\
\FF \times \nerve \UU 2 \ar [rrr] _ { \id \FF \times \del \jj 0 } &&&
\FF \times \UU }
\end {equation}
Let the Euclideanness \two map
be \( \FF \times \nerve \UU 3
\maplongto { \id \FF \times \pair { \del \jj { 0 1 } } { \del \jj { 1 2 } } }
\FF \times \nerve \UU 2 \times \nerve \UU 2
\maplongto { \id \FF \times \gg \times \id { \nerve \UU 2 } }
\FF \times \GG \times \nerve \UU 2
\maplongto { \rr \times \id { \nerve \UU 2 } } \FF \times \nerve \UU 2 \).
In (\ref {2Euclideanness}),
the natural isomorphism \( \del \omega { 0 1 } \) is trivial,
while \( \del \omega { 1 1 } \)
is given by this string diagram:
\begin {equation} \label {2bundle construction eq2rel Euclidean string}
\matrixxy
{ \ddStringl \ww && && \dstring \gg \\ & \xx & && && & \zz \\
&& \dstring \gg \rrstring \gamma && && \dstring \gg \\ && && \yy \\ && && && }
\end {equation}
\par Therefore, \I really do have an equivalence \two relation,
so some \two quotient \( \FF \times \UU \maplongto { \tilde \jj } \EE \)
must exist.
The isomorphism (\ref {2quotient}) becomes \( \theta \);
while the coherence laws
(\ref {2quotient coherence reflexivity}, \ref {2quotient coherence Euclid})
become (respectively)
the coherence laws (\ref {G2bundle unit}, \ref {G2bundle ass}).
\par To define the bundle \two map \( \EE \maplongto \pp \BB \),
first note this natural isomorphism:
\begin {equation} \label {2bundle construction projection cone}
\matrixxy
{ \FF \times \nerve \UU 2
\ar [dddd] _ { \id \FF \times \del \jj 0 }
\ar [rrrr] ^ { \id \FF \times \pair \gg { \del \jj 1 } } &&&&
\FF \times \GG \times \UU \ar [rrr] ^ { \rr \times \id \UU } &&&
\FF \times \UU \ar [dd] ^ { \hat \FF \times \id \UU } \\
&&& \relax \ddriicell \empty \\ &&&& &&& \UU \ar [dd] ^ \jj \\ &&& & \\
\FF \times \UU \ar [rrrr] _ { \hat \FF \times \id \UU } &&&&
\UU \ar [rrr] _ \jj &&& \BB }
\end {equation}
which is simply one of the structural isomorphisms of the \two cover~\( \UU \).
Since \( \EE \) is \a \two quotient,
this defines \a \two map \( \EE \maplongto \pp \BB \)
and \a natural isomorphism as in (\ref {summary 2pullback}).
\par Therefore, \( \EE \)
is \a \two bundle over \( \BB \) with fibre \( \FF \)
associated with the \two transition \( \gg \).
\end {prf}
\par If I'm to conclude, as in part~\ref 1,
that this \( \EE \)
is unique up to an essentially unique equivalence of \GGG \two bundles,
then the next step is to consider \GGG \two bundle morphisms.
\begin {prop} \label {associated 2bundle morphism}
Given \GGG \two bundles \( \EE \) and \( \EE ' \)
(both over \( \BB \) and with fibre \( \FF \))
associated with \GGG \two transition morphisms
\( \gg \) and \( \gg ' \) (respectively),
and given \a \GGG \two transition morphism \( \bb \)
from \( \gg \) to \( \gg ' \),
there is \a \GGG \two bundle morphism \( \ff \) from \( \EE \) to \( \EE ' \)
associated with \( \bb \).
\end {prop}
\begin {prf}
First consider this natural isomorphism:
\begin {equation}
\label {2bundle morphism construction intermediate quotient cone}
\matrixxy
{ \FF \times \UU \cap \nerve { \UU ' } 2
\ar [dddddd] _ { \id \FF \times \del \jj { 0 1 } }
\ar [rrrrrrrrr] ^ { \id \FF \times \del \jj { 0 2 } } &&&& &&& &&
\FF \times \UU \cap \UU '
\ar [dd] ^ { \id \FF \times \pair \bb { \del \jj 1 } } \\\\
&&&& \relax \ddriicell \empty &&& &&
\FF \times \GG \times \UU ' \ar [dd] ^ { \rr \times \id \UU } \\\\
&&&& &&& && \FF \times \UU \ar [dd] ^ { \tilde \jj ' } \\\\
\FF \times \UU \cap \UU '
\ar [rrrr] _ { \id \FF \times \pair \bb { \del \jj 1 } } &&&&
\FF \times \GG \times \UU ' \ar [rrr] _ { \rr \times \id { \UU ' } } &&&
\FF \times \UU ' \ar [rr] _ { \tilde \jj ' } && \EE ' }
\end {equation}
given by this string diagram:
\begin {equation}
\label {2bundle morphism construction intermediate quotient cone string}
\matrixxy
{ \dddStringl \ww && && \dstring \bb && && \ddString \\ && && && & \yy ' \\
& \xx & \ddstring \bb \rrstring \delta && && \dstring { \gg ' } \\\\
&& && \xx ' && \rString && \theta ' \dString \\\\ && && && && }
\end {equation}
Since \( \id \FF \times \del \jj 0 \), being \a \two cover,
is \a \two quotient of its kernel pair,
\I can construct from this
\a \two map \( \FF \times \UU \maplongto { \tilde \ff } \EE \)
and \a natural isomorphism~\( \tilde \zeta \):
\begin {equation} \label {2bundle morphism construction intermediate quotient}
\matrixxy
{ \FF \times \UU \cap \UU '
\ar [ddddrr] _ { \id \FF \times \del \jj 0 }
\ar [rrrr] ^ { \id \FF \times \pair \bb { \del \jj 1 } } &&&&
\FF \times \GG \times \UU ' \ar [dd] ^ { \rr \times \id { \UU ' } } \\
&& \relax \ddricell { \tilde \zeta } \\
&& && \FF \times \UU \ar [dd] ^ { \tilde \jj ' } \\ && & \\
&& \FF \times \UU \ar [rr] _ { \tilde \ff } && \EE ' }
\end {equation}
In terms of elements
\( \ww \) of \( \FF \) and \( \pair \xx { \xx ' } \) of \( \UU \cap \UU ' \),
this gives an arrow \( \zetatildesub \ww \xx { \xx ' } \)
from \( \jjptildesub { \rrsub \ww { \bbsub \xx { \xx ' } } } { \xx ' } \)
to the composite of \( \pair \ww \xx \) with \( \tilde \ff \),
as in this string diagram:
\begin {equation}
\label {2bundle morphism construction intermediate quotient string}
\matrixxy { \ddStringl \ww && \dstring \bb && \dString \\ && & \xx ' \\
&& \rString && \tilde \zeta \\ && \xx \\ & }
\end {equation}
where \( \tilde \ff \)
is applied to the portion of the diagram left open on the right.
\par Next, notice this natural isomorphism:
\begin {equation} \label {2bundle morphism construction quotient cone}
\matrixxy
{ \FF \times \nerve \UU 2
\ar [dd] _ { \id \FF \times \del \jj 0 }
\ar [rrrr] ^ { \id \FF \times \pair \gg { \del \jj 1 } } &&&
\relax \ddriicell \empty &
\FF \times \GG \times \UU \ar [rrr] ^ { \rr \times \id \UU } &&&
\FF \times \UU \ar [dd] ^ { \tilde \ff } \\\\
\FF \times \UU \ar [rrrrrrr] _ { \tilde \ff } &&&& &&& \EE ' }
\end {equation}
given by this string diagram:
\begin {equation} \label {2bundle morphism construction quotient cone string}
\matrixxy
{ \ddddStringl \ww \\\\
&& \xx && \dstring \bb \rrString && && \tilde \zeta \ddString \\
&& && && & \xx ' \\ && \ddstring \gg \rrstring \sigma && && \dstring \bb \\\\
&& && \yy && \rString && \tilde \zeta \\\\ && }
\end {equation}
(Notice that \I can drop \( \xx ' \)
because \( \FF \times \nerve \UU 2 \cap \UU '
\maplongto { \del \jj { 0 1 } } \FF \times \nerve \UU 2 \)
is \two epic.)
Since \( \tilde \jj \), being the pullback of \a \two cover,
is \a \two cover and hence the \two quotient of its \two kernel pair,
\I can construct from this \a \two map \( \EE \maplongto \ff \EE ' \)
and \a natural isomorphism \( \zeta \) as in~(\ref {G2bundle morphism}).
As with \( \theta \) in the previous proposition,
the coherence laws for \( \zeta \)
come from the coherence laws for \a \two quotient.
\par In short,
I've constructed \a \GGG \two bundle morphism \( \EE \maplongto \ff \EE ' \)
associated with \( \bb \).
\end {prf}
\begin {prop}
Given \GGG \two bundle morphisms \( \ff \) and \( \ff ' \)
associated (respectively) with \( \bb \) and \( \bb ' \),
any \GGG \two transition \two morphism \( \xi \) from \( \bb \) to \( \bb ' \)
has associated with it \a unique \GGG \two bundle \two morphism
\( \kappa \) from \( \ff \) to~\( \ff ' \).
\end {prop}
\begin {prf}
Let \( \tilde \ff ' \)
be the composite
\( \FF \times \UU \maplongto { \tilde \jj } \EE \maplongto \ff \EE ' \).
Then \( \kappa \) is \( \tilde \omega _ \xx \)
in (\ref {2quotient universality}),
while (\ref {2quotient universality coherence})
gives (\ref {G2bundle 2morphism}).
\end {prf}
\par In fact, there is \a \two functor
from the \two category \( \Bun \CCC \GG \BB \FF \)
of \GGG \two bundles over \( \BB \) with fibre \( \FF \)
to the \two category \( \BB ^ \GG \) of \GGG \two transitions over \( \BB \),
defined by simply forgetting the total space \( \EE \).
The propositions above, taken altogether,
prove that this \two functor is an equivalence of \two categories.
That is:
\begin {thm} \label {2bundle 2cat equals 2transition 2cat}
Given \a \two space \( \BB \), \a \two group \( \GG \),
and \a \GGG \two space \( \FF \),
the \two category \( \Bun \CCC \GG \BB \FF \)
of \GGG \two bundles over \( \BB \) with fibre \( \FF \)
is equivalent to the \two category \( \BB ^ \GG \)
of \GGG \two transitions on \( \BB \).
\end {thm}
So depending on the point of view desired,
you can think locally of \GGG \two transitions over \( \BB \),
or globally of principal \GGG \two bundles over \( \BB \)
(or \GGG \two bundles over \( \BB \) with some other fibre \( \FF \)).
\subsection {Semistrictification} \label {2semistrict}
In this extra section,
\I show how to replace any \GGG \two transition
(and hence any \GGG \two bundle)
with \a \strong {semistrict} \GGG \two transition (or \GGG \two bundle),
defined as one where \( \eta \) is trivial,
more precisely where \( \etasub \xx = \iotasub { \ggsub \xx \xx } \).
Notice that if \( \GG \) is \a strict \two group,
then this \( \eta \) is an identity;
but reducing \( \eta \) to \a special case of \( \iota \)
is the best that you can hope for in general,
since you can't have an identity morphism between things that aren't equal!
Unlike most of this part,
it's \emph {not} sufficient that \( \GG \) be merely \a \two monoid;
\I will really need both \( \iota \) and \( \epsilon \)
(even in the strict case).
\par Given \a \GGG \two transition \( \gg \),
to define the \strong {semistrictification} \( \gg ' \) of \( \gg \),
let \( \ggpsub \xx \yy \)
be \( \mmsub { \ggsub \xx \yy } { \iisub { \ggsub \yy \yy } } \),
as in this string diagram:
\begin {equation} \label {2transition sstrict string}
\matrixxy { & \dstring { \gg ' } \\ \xx & & \yy \\ & }
\qquad : = \qquad
\matrixxy { & \dstring \gg && \ustring \gg \\ \xx & & \yy & & \yy \\ & && }
\end {equation}
Let \( \gammapsub \xx \yy \zz \) and \( \etapsub \xx \) be as follows:
\begin {equation} \label {2transition sstrict prod}
\matrixxy
{ \dstring { \gg ' } && && \dstring { \gg ' } \\ && \yy \\
\rrstring { \gamma ' } && \dstring { \gg ' } && \\ \xx && && \zz \\ && }
\qquad : = \qquad
\matrixxy
{ \dddstring \gg && \uustring \gg && && \dstring \gg && && \uuuustring \gg \\
&& & \yy \\
&& && \dstring \gg \rrstring { \bar \gamma } && && \ddstring \gg \\\\
&& \rstring \epsilon && && \yy \\ && && && && && & \zz \\
\rrrrstring \gamma && && \dstring \gg && && & \zz \\ \xx \\ && && && && && }
\end {equation}
and:
\begin {equation} \label {2transition sstrict unit}
\matrixxy
{ & \estring \\ \xx \\ & \pstring { \eta ' } \dstring { \gg ' } \\\\ & }
\qquad : = \qquad
\matrixxy
{ & \estring \\ \xx \\
\dstring \gg \rstring \iota && \ustring \gg \\ & \xx & & \\ && }
\end {equation}
That is, \( \etapsub \xx = \iotasub { \ggsub \xx \xx } \) as promised.
\begin {prop}
The semistrictification \( \gg ' \) is \a \GGG \two transition.
\end {prop}
\begin {prf}
The left unit law in (\ref {2transition unit}) is the simplest;
it follows using string diagrams from one application of \a zig-zag identity
and cancelling an inverse:
\begin {eqnarray} \label {2transition sstrict unit left}
&
\matrixxy
{ && && && \dstring \gg && && \uuuustring \gg \\ \xx \\
\ddstring \gg \rstring \iota && \ustring \gg &&
\dstring \gg \rrstring { \bar \gamma } && && \ddstring \gg \\\\
&& \rstring \epsilon && && \xx \\ && && && && & \yy & & \yy \\
\rrrrstring \gamma && && \dstring \gg && && \\\\ && && && && && }
& \nonumber \\ \\ & = \qquad
\matrixxy
{ & && \dstring \gg && && \uuustring \gg \\ \xx & && && & \yy \\
& \dstring \gg \rrstring { \bar \gamma } && && \dstring \gg \\
& && \xx && && & \yy \\ & \rrstring \gamma && \dstring \gg && \\\\ & && && && }
\qquad = \qquad
\matrixxy
{ & \ddstring \gg && \uustring \gg \\\\ \xx & & \yy & & \yy \\\\ & && }
& \nonumber
\end {eqnarray}
The other coherence laws may also be proved by string diagrams;
for example, to prove the associative law~(\ref {2transition ass}),
\I apply the associative law for the original \( \gg \),
then introduce \( \gammasub \xx \xx \zz \) and its inverse,
which allows me to use the original associative law again,
and then \I just remove \( \gammasub \xx \xx \yy \) and its inverse.
\end {prf}
To see that \( \gg \) and \( \gg ' \) are equivalent,
let \( \bbsub \xx \yy \)
be \( \mmsub { \ggsub \xx \yy } { \iisub { \ggsub \yy \yy } } \)
(the same as \( \ggpsub \xx \yy \)),
and let \( \bbpsub \xx \yy \) be \( \ggsub \xx \yy \).
\begin {prop}
\( \bb \) is \a \GGG \two transition morphism from \( \gg \) to \( \gg ' \),
\( \bb ' \) is \a \GGG \two transition morphism from \( \gg ' \) to \( \gg \),
there is an invertible \GGG \two transition \two morphism \( \xi \)
from \( \bb ; \bb ' \) to the identity on \( \gg \),
and there is an invertible \GGG \two transition \two morphism \( \xi ' \)
from \( \bb ' ; \bb \) to the identity on \( \gg ' \).
\end {prop}
In other words, \( \gg \) and \( \gg ' \) are equivalent.
\begin {prf}
\( \sigma \), \( \delta ' \), and \( \xi ' \) work just like \( \gamma \);
while \( \delta \), \( \sigma ' \) and \( \xi \) work like \( \gamma ' \).
(Sometimes there's an extra \( \iisub { \ggsub \zz \zz } \)
multiplied on the right,
but it doesn't get involved.)
The proofs of the coherence laws
are just like the associative law for \( \gamma \)
(laws~(\ref {2transition morphism ass left},
\ref {2transition morphism ass bi})
for \( \bb \),
(\ref {2transition morphism ass right}, \ref {2transition morphism ass bi})
for \( \bb ' \),
(\ref {2transition 2morphism left}) for \( \xi \),
and (\ref {2transition 2morphism left}, \ref {2transition 2morphism right})
for \( \xi ' \));
the left unit law for \( \gamma \)
(law~(\ref {2transition morphism unit left}) for \( \bb \));
the associative law for \( \gamma ' \)
(law~(\ref {2transition morphism ass right}) for \( \bb \),
(\ref {2transition morphism ass left}) for \( \bb ' \),
and (\ref {2transition 2morphism right}) for \( \xi \));
the right unit law for \( \gamma ' \)
(law~(\ref {2transition morphism unit right}) for \( \bb \));
the left unit law for \( \gamma ' \)
(law~(\ref {2transition morphism unit left}) for \( \bb ' \));
and the right unit law for \( \gamma \)
(law~(\ref {2transition morphism unit right}) for \( \bb ' \)).
\end {prf}
Notice that the original \( \eta \)
is never used in any of these constructions.
So the result is not only
that every \GGG \two transition may be semistrictified,
but also that it's enough to give \( \gg \) and \( \gamma \)
(and to check associativity for this \( \gamma \))
to define an (equivalent) semistrict \GGG \two transition.
\part {Examples} \label {example}
\A few simple examples will help to clarify the ideas;
\a more complicated example will link them to the previous notion of gerbe.
In these examples,
\I will not hesitate to assume that the category \( \C \)
is indeed the category of smooth manifolds and smooth functions.
\par These examples,
although they live in the \two category \( \CCC \) of \two spaces,
are built of the category \( \C \) of spaces.
Accordingly, \I will use italic Roman letters for spaces and maps;
then, as in section~\ref {disc2space},
\I will use calligraphic or Fraktur letters
for the \two spaces and \two maps built out of them.
\section {Bundles as \two bundles} \label {discexample}
Just as every space is \a categorically discrete \two space,
so every group is \a categorically discrete \two group,
and every bundle is \a categorically discrete \two bundle.
Just apply section~\ref {disc2space} to every space involved;
all of the natural transformations needed can just be identities.
This is trivial and uninteresting.
\section {Vector bundles} \label {vecexample}
Let \( E \) be \a vector bundle over the base space \( B \).
That is, the fibre \( F \) is \a (real or complex) vector space,
and the action of \( G \) on \( F \) is \a linear representation of \( G \).
To define \a \two bundle over \( B \),
let the space of points of the total \two space be \( \cat \EE 1 : = B \)
and let its space of arrows be \( \cat \EE 2 : = E \).
Let the source and target maps
\( \cat \EE 2 \longto \cat \EE 1 \)
both be the original projection map \( E \maplongto p B \).
Let the identity arrow map
\( \cat \EE 1 \maplongto { \dom d { 0 0 } } \cat \EE 2 \)
be the zero section that assigns each point to the zero vector over it.
Similarly, let the composition map \( \dom d { 0 2 } \)
be given by vector addition.
(In the case where \( E \) is the tangent bundle \( T B \),
this \two space was described in the course of \cite [Example~48] {HDA5},
as the \q {tangent groupoid} of \( B \),
so called because this \two space, viewed as \a category, is \a groupoid.
But \I would call this example the \strong {tangent \two bundle},
because as you will see, it has the structure of \a \two bundle.)
Returning to the general case, \I have \a total \two space \( \EE \);
the projection \two map \( \EE \maplongto \pp B \)
is given by the identity \( B \maplongto B B \) on points
and by \( E \maplongto p B \) on arrows.
\par Now, in order to get \a \two bundle,
\I need to identify the appropriate \two group.
This was defined in \cite [Example~50] {HDA5} in terms of Lie crossed modules.
In the terms of the present paper,
the space of points of the \two group is simply \( \cat \GG 1 : = G \),
while the space of arrows is \( \cat \GG 2 : = G \times F \).
The source and target maps
are both simply
the projection \( G \times F \maplongto { G \times \hat F } G \).
The identity arrow map
maps \( x \in G \) to \( \lr ( { 0 , x } ) \in G \times F \),
while composition is again given by vector addition.
\par The multiplication \( \mm \) is more interesting.
Since the vector space \( F \) is an abelian group
and the action of \( G \) on \( F \) is linear,
the space \( \cat \GG 2 \)
can be interpreted as \a semidirect product \( G \imestimes F \).
Thus inspired, on points,
\( \cat \GG 1 \times \cat \GG 1 \maplongto { \cat \mm 1 } \cat \GG 1 \)
is simply \( G \times G \maplongto m G \).
But on arrows,
\I let the product of \( \lr ( { x , \vec w } ) \)
and \( \lr ( { x ' , \vec w ' } ) \)
be \( \lr ( { x x ' , \vec w x ' + \vec w ' } ) \).
(This formula for multiplication can be internalised with \a diagram,
thus generalising the situation to an arbitrary category \( \C \),
but \I won't do that here.
In fact, to do it perfectly properly,
\I would need to replace the vector space \( F \)
with an internal abelian group
---although this is not \a bad thing,
since this can be an interesting generalisation
even when \( \C \) is the category of smooth manifolds.)
\par Next, I'll define the fibre \two space \( \FF \).
Let the point space \( \cat \FF 1 \) be the singleton space \( 1 \),
and let the arrow space \( \cat \FF 2 \) be \( F \).
The source and target maps
are both the unique map \( F \maplongto { \hat F } 1 \).
Again, the identity arrow map is given by the zero vector,
and the composition \two map is given by vector addition.
If \( U \maplongto j B \) is \a cover map for the original vector bundle
(and in the case of \a vector bundle,
any open cover by Euclidean neighbourhoods will suffice),
then the pulled-back \two map
\( \FF \times \UU \maplongto { \tilde \jj } \EE \)
is given by \( U \maplongto j B \) on points
and \( F \times U \maplongto { \tilde j } E \) on arrows.
Some tedious checking shows that this gives \a \two pullback:
\begin {equation} \label {vecexample pullback}
\matrixxy
{ \FF \times \UU
\ar [dd] _ { \hat \FF \times \id \UU } \ar [rr] ^ { \tilde \jj } &&
\relax \EE \ar [dd] ^ { \pp } \\\\
\UU \ar [rr] _ j && \BB }
\end {equation}
The action of \( \GG \) on \( \FF \) is simple enough.
Of course, the point map
\( \cat \FF 1 \times \cat \GG 1 \maplongto { \cat \rr 1 } \cat \FF 1 \)
is necessarily \( 1 \times G \maplongto { \hat G } 1 \).
But the arrow map
\( \cat \FF 2 \times \cat \GG 2 \maplongto { \cat \rr 2 } \cat \FF 2 \)
is more interesting,
mapping \( \lr ( { \vec w , \lr ( { x , \vec w ' } ) } ) \)
to \( \vec w x + \vec w ' \).
\par Next, the \two transition \( \nerve \UU 2 \maplongto \gg \GG \)
is simply \( \nerve U 2 \maplongto { \hat \nerve U 2 } 1 \) on points,
and \( \nerve U 2 \maplongto g G \) on arrows.
Now all of the \two maps are defined.
\par \I still need to define the natural transformations
and verify the coherence laws.
But the natural transformations can all be taken to be identities,
so the coherence laws are all automatically true.
Therefore, \I have \a \GGG \two bundle.
\section {Gerbes} \label {gerbexample}
For gerbes, \I follow Breen in \cite [chapter~2] {Asterisque}.
As noted in \cite [2.13] {Asterisque},
this discussion applies most generally to \a crossed module,
denoted by Breen as \( G \maplongto \delta \Pi \).
Now, for Breen, this is not \a single crossed module
but rather \a \emph {sheaf} of crossed modules over \a space \( X \).
\I will (again) take the base \two space \( \BB \)
to be the \two space built trivially out of the space \( X \)
as in section~\ref {disc2space};
but \I do not reach the generality of an arbitrary sheaf of crossed modules.
Instead, given \a strict \two group \( \GG \),
interpreted as \a crossed module as in section~\ref {2crosmod},
let Breen's \( G \) be the sheaf of \H valued functions on \( X \),
and let his \( \Pi \) be the sheaf of \D valued functions on \( X \).
Then acting pointwise,
these form \a sheaf of crossed modules \( G \maplongto \delta \Pi \).
(Thus gerbes are more general than \two bundles,
since only sheaves of this form can appear,
much as sheaves themselves are more general than bundles;
although \two bundles are still more general than gerbes
---at least as so far as anything has been published---
in that only strict \two groups have so far appeared with gerbes.)
\par In \cite [2.4] {Asterisque},
Breen describes \a gerbe in terms of local data as follows:
\begin {closeitemize}
\item an open cover \( \famb { U _ i } { i \in I } \) of \( X \);
\item an open cover
\( \famb { U ^ \alpha _ { i j } } { \alpha \in J _ { i j } } \)
of each double intersection \( U _ i \cap U _ j \);
\item \a section \( \lambda ^ \alpha _ { i j } \) of \( \Pi \)
over each \( U ^ \alpha _ { i j } \);
and \item \a section \( g ^ { \alpha \beta \gamma } _ { i j k } \) of \( G \)
over each triple intersection
\( U ^ \alpha _ { i j } \cap U ^ \beta _ { j k } \cap U ^ \gamma _ { i k } \);
such that \item there is only one \( \alpha \) in each \( J _ { i i } \);
\item \( \msub { \lambda ^ \alpha _ { i j } } { \lambda ^ \beta _ { j k } } \)
is equal to
the value of \( g ^ { \alpha \beta \gamma } _ { i j k } \) in \( \Pi \)
multiplied by \( \lambda ^ \gamma _ { i k } \)
(\cite [(2.4.5)] {Asterisque}, or \cite [(2.1.5)] {BM});
\item \( \msub { g ^ { \alpha \beta \gamma } _ { i j k } }
{ g ^ { \gamma \eta \delta } _ { i k l } } \)
is equal to the action of \( \lambda ^ \alpha _ { i j } \)
on \( g ^ { \beta \eta \epsilon } _ { j k l } \),
multiplied by \( g ^ { \alpha \epsilon \delta } _ { i j l } \),
on each sextuple intersection
\( U ^ \alpha _ { i j } \cap U ^ \beta _ { j k } \cap U ^ \gamma _ { i k }
\cap U ^ \eta _ { k l } \cap U ^ \epsilon _ { j l }
\cap U ^ \delta _ { i l } \)
(\cite [(2.4.8)] {Asterisque}, or \cite [(2.1.6)] {BM});
\item \( \lambda ^ \alpha _ { i i } \) is the identity in \( \Pi \)
(for \( \alpha \) the unique element of \( J _ { i i } \));
\item \( g ^ { \alpha \beta \beta } _ { i i k } \) is the identity in \( G \)
(for \( \alpha \) again the unique element of \( J _ { i i } \));
and \item \( g ^ { \alpha \beta \alpha } _ { i k k } \)
is the identity in \( G \)
(for \( \beta \) now the unique element of \( J _ { k k } \)).
\end {closeitemize}
Since Breen nowhere gives precisely this list, \a few words are in order:
Breen's definition of \a \q {labeled decomposition}
is given directly in terms of certain \( \phi ^ \alpha _ { i j } \),
out of which the \( \lambda \) and \( g \) are constructed;
but \cite [2.6] {Asterisque}
shows conversely how to recover \( \phi \) from \( \lambda \) and \( g \).
Also, Breen never mentions \( \Pi \) at this point,
only the automorphism group \( \Aut G \) of \( G \);
he is only considering the case where \( \Pi \) is this automorphism group,
but he explains how to generalise this when he discusses crossed modules,
and \I have followed that here.
Finally, Breen suppresses the unique element of each \( J _ { i i } \),
but I've restored this since \I find it confusing to suppress it.
\par Here is how these data
corresponds to the local data of \a \GGG \two bundle,
that is \a \GGG \two transition:
\begin {closeitemize}
\item The open cover \( \faml { U _ i } \)
corresponds to \a cover \( U \maplongto j X \),
which (as explained in section~\ref {222})
is equivalent to \a \two cover \( \UU \maplongto \jj \BB \).
\item The open cover \( \faml { U ^ \alpha _ { i j } } \)
corresponds to the cover \( \src \gg \longto U \)
involved in the transition \two map \( \UU \maplongto \gg \GG \).
\item The section \( \lambda ^ \alpha _ { i j } \)
is simply the map \( \src \gg \maplongto { \cat \gg 1 } D = \cat \GG 1 \).
\item The (disjoint union of the) triple intersections
\( U ^ \alpha _ { i j } \cap U ^ \beta _ { j k } \cap U ^ \gamma _ { i k } \)
may be taken as the domain of the map \( \src \gamma \)
for the natural isomorphism \( \gamma \) of diagram~(\ref {2transition}).
Then \( \src \gamma \) takes values in \( \cat \GG 2 = H \semitimes D \),
and Breen's \( g \) is the \H valued component of this map.
\item The simplicity of the open cover of \( U _ i \cap U _ i \)
should be ensured before the translation begins
by passing to \a refinement.
(This doesn't work for the \( U _ i \cap U _ j \) in general,
since you might not be able to coordinate them.
Breen discusses this in \cite [2.3] {Asterisque}.)
\item The requirement
on the value of \( g ^ { \alpha \beta \gamma } _ { i j k } \) in \( \Pi \)
fixes the \D valued component of \( \src \gamma \);
that is, it indicates an appropriate source and target for \( \gamma \).
\item The requirement
for the action of \( \lambda ^ \alpha _ { i j } \)
on \( g ^ { \beta \eta \epsilon } _ { j k l } \),
following the principles laid out in section~\ref {2crosmod},
is precisely the associativity coherence law~(\ref {2transition ass}).
\item That \( \lambda ^ \alpha _ { i i } \) is trivial
is ensured by passing to the semistrictification;
see section~\ref {2semistrict}.
\item This done,
that \( g ^ { \alpha \beta \beta } _ { i i k } \)
and \( g ^ { \alpha \beta \alpha } _ { i k k } \)
are trivial
is the unit coherence law~(\ref {2transition unit}).
\end {closeitemize}
Thus, gerbes correspond to (semistrict) \two bundles.
\par But to get \a good comparison,
this should respect equivalence of each.
(Already \I had to replace an arbitrary \two bundle
with its equivalent semistrictification.)
Given two gerbes (which Breen denotes as primed and unprimed),
here are the local data for an equivalence of gerbes:
\begin {closeitemize}
\item \a common refinement
of the open covers \( \faml { U _ i } \) and \( \faml { U ' _ i } \)
(renamed \( \faml { U _ i } \) by Breen);
\item common refinements of each \( \faml { U ^ \alpha _ { i j } } \)
on \a given \( U _ i \cap U _ j \);
\item \a further refinement to allow \( \mu _ i \) below to make sense;
\item \a section \( \mu _ i \) of \( \Pi \) on each \( U _ i \);
and \item \a section \( \delta ^ \alpha _ { i j } \) of \( G \)
on each \( U ^ \alpha _ { i j } \);
such that \item \( \delta ^ \alpha _ { i i } \) is the identity in \( \Pi \)
(for \( \alpha \) the unique element of \( J _ { i i } \));
\item \( \msub { \lambda ' ^ \alpha _ { i j } } { \mu _ j } \)
is equal to the value of \( \delta ^ \alpha _ { i j } \) in \( \Pi \)
multiplied by \( \msub { \mu _ i } { \lambda ^ \alpha _ { i j } } \)
(\cite [2.4.16] {Asterisque}, or \cite [2.1.11] {BM});
and \item \( \msub { g ' ^ { \alpha \beta \gamma } _ { i j k } }
{ \delta ^ \gamma _ { i k } } \)
is equal to the action of \( \lambda ' ^ \alpha _ { i j } \)
on \( \delta ^ \beta _ { j k } \),
multiplied by \( \delta ^ \alpha _ { i j } \),
multiplied by the action of \( \mu _ i \)
on \( g ^ { \alpha \beta \gamma } _ { i j k } \)
(\cite [(2.4.17)] {Asterisque}, or \cite [(2.1.12)] {BM}).
\end {closeitemize}
\I use primes like Breen,
but it turns out that the order here is \emph {backwards};
so as you compare the above
to the data of \a \GGG \two transition morphism,
be sure to swap the primes in the data for the \GGG \two transitions.
Here is the precise correspondence:
\begin {closeitemize}
\item The common refinement of the original open covers
is taken automatically in forming the pullback \( \UU \cap \UU ' \),
the source of \( \bb \).
\item The common refinement of the open covers' hypercovers
gives the source \( \src \bb \) of \( \cat \bb 1 \).
\item The final refinement is only taken
to ensure that the new \( J _ { i i } \) are again singletons.
\item \( \mu _ i \) corresponds to \( \bbsub \xx \xx \).
In general, \( \bbsub \xx \yy \) may be recovered
up to isomorphism as \( \mmsub { \ggsub \xx \yy } { \bbsub \yy \yy } \).
\item Breen's \( \delta ^ \alpha _ { i j } \)
corresponds to my \( \sigmasub \xx \yy \yy \)
followed by \( \ndeltasub \xx \xx \yy \).
\item That \( \sigmasub \xx \xx \xx \) followed by \( \deltasub \xx \xx \xx \)
be trivial
fixes \( \sigma \) and \( \delta \) up to isomorphism.
\item The value of \( \delta ^ \alpha _ { i j } \) in \( \Pi \) holds
because \( \mmsub { \ggsub \xx \yy } { \bbsub \yy \yy }
\maplongTo { \sigmasub \xx \yy \yy } \bbsub \xx \yy
\maplongTo { \ndeltasub \xx \xx \yy }
\mmsub { \bbsub \xx \xx } { \ggpsub \xx \yy } \).
\item The last item
is \a combination of the associativity coherence laws
(\ref {2transition morphism ass left},
\ref {2transition morphism ass right}, \ref {2transition morphism ass bi}).
\end {closeitemize}
Specifically, the formulation of the last item in \cite {BM},
translated into my notation, gives these string diagrams:
\begin {equation} \label {bundle morphism coherence gerb}
\matrixxy
{ \dstring \gg && && \dstring \gg && \ddstring \bb \\ && \yy \\
\rrstring \gamma && \dstring \gg && \\ && && & \zz \\
&& \rrstring \sigma && \dstring \bb && & \zz \\ \xx \\
&& \dstring \bb \rrstring { \bar \delta } && && \dstring { \gg ' } \\
&& && \xx \\ && && && }
\qquad = \qquad
\matrixxy
{ & \dddstring \gg && \dstring \gg && \dstring \bb \\ & && & \zz \\
& && \rstring \sigma && \dstring \bb & \\ & & \yy \\
& && \dstring \bb \rstring { \bar \delta } && \dddstring { \gg ' } \\
\xx \\ & \rstring \sigma & \dstring \bb & && & \zz \\ & && & \yy \\
& \ddstring \bb \rstring { \bar \delta } && \dstring { \gg ' } \\\\
& & \xx & \rstring \gamma & \dstring { \gg ' } & \\\\ & && & }
\end {equation}
(Note that \cite {Asterisque} uses slightly more complicated phrasing,
as \cite {BM} points out in a footnote;
this corresponds to introducing into the string diagram an extra zig-zag,
which can be removed by (\ref {first zig-zag string}),
the first zig-zag identity).
Thus, equivalent gerbes and equivalent \two bundles correspond.
\par To make this correspondence even more precise,
\I should consider composition of \two bundle morphisms
as well as \two bundle \two morphisms.
The analogues of these for gerbes are not discussed in \cite {Asterisque},
but they are discussed by Breen \& Messing in \cite {BM},
to which \I now turn exclusively.
\I should note that Breen \& Messing are less general in \cite {BM}
than Breen alone was in \cite {Asterisque},
in that they assume that \a single open cover may be chosen once and for all
that obviates the need for the upper Greek superscripts.
(As they remark in \cite [2.1] {BM},
this assumption is valid for certain choices of \( \C \),
such as affine schemes or paracompact manifolds,
but not for more general schemes or manifolds.)
\I will restore the missing indices.
Also, they change the notation \q {\( \mu \)} to \q {\( m \)}
and use superscripts on this instead of primes;
\I will continue to write \q {\( \mu \)}
and to use primes to distinguish morphisms
(although again they will be reverse from my primes,
since composition in both \cite {BM} and \cite {Asterisque}
is in the opposite order from this paper).
\par This in mind, \a \two morphism between morphisms of gerbes
is described locally by the following data:
\begin {closeitemize}
\item \a common refinement of the open covers involved in the gerbe morphisms
(suppressed in \cite {BM});
\item \a further refinement to allow \( \theta _ i \) below to make sense
(again suppressed in \cite {BM},
but analogous to the refinement for \( \mu \) last time);
and \item sections \( \theta _ i \) of \( G \) on each \( U _ i \);
such that \item \( \mu ' _ i \)
is equal to value of \( \theta _ i \) in \( \Pi \),
multiplied by \( \mu _ i \)
(\cite [(2.3.7)] {BM});
and \item \( \msub { \delta ' ^ \alpha _ { i j } } { \theta _ i } \)
is equal to the action of \( \lambda ' ^ \alpha _ { i j } \)
on \( \theta _ j \),
multiplied by \( \delta ^ \alpha _ { i j } \)
(\cite [(2.3.8)] {BM});
\end {closeitemize}
These correspond to \a \two transition \two morphism as follows:
\begin {closeitemize}
\item The common refinement is unnecessary,
since the same pullback \( \UU \cap \UU ' \) applies to both morphisms.
\item The other refinement again ensures that each new \( J _ { i i } \)
is \a singleton.
\item \( \theta _ i \) corresponds to \( \xisub \xx \xx \).
\item The value of \( \theta _ i \) in \( \Pi \)
shows that \( \xisub \xx \xx \)
is an arrow from \( \bbsub \xx \xx \) to \( \bbpsub \xx \xx \).
\item One of the coherence laws
(\ref {2transition 2morphism left}, \ref {2transition 2morphism right})
may be taken to define \( \xisub \xx \yy \) in general;
then the other
gives the desired action
of \( \lambda ' ^ \alpha _ { i j } \) on \( \theta _ j \).
\end {closeitemize}
Here are the full string diagrams for the last item:
\begin {equation} \label {bundle 2morphism coherence gerb}
\matrixxy
{ & \dstring \gg && && \dstring \bb \\ & && \yy \\
& \rrstring \sigma && \dstring \bb && \\\\
\xx & \dstring \bb \rrstring { \bar \delta } && &&
\ddstring { \gg ' } & \yy \\\\
& \dstring { \bb ' } \pstring \chi && \xx \\\\ && && && }
\qquad = \qquad
\matrixxy
{ & \ddstring \gg && && \dstring \bb \\\\
& && \yy && \dstring { \bb ' } \pstring \chi \\ \xx & && && & \yy \\
& \rrstring { \sigma ' } && \dstring { \bb ' } && \\\\
& \dstring { \bb ' } \rrstring { \bar \delta ' } && && \dstring { \gg ' } \\
& && \xx \\ & && && }
\end {equation}
\par Both here and in \cite {BM},
composition of \two morphisms is given by multiplication in \( G \),
so the correspondence between \two morphisms is functorial.
In \cite {BM}, whiskering \a \two morphism by \a morphism on the left
does nothing,
while whiskering on the right applies \( \mu \) to \( \theta \).
This is precisely what should happen,
since \q {left} and \q {right} are reversed,
in the string diagram for \a crossed module.
Thus \I have proved the following theorem:
\begin {thm} \label {principal 2bundle 2cat equals gerbe 2cat}
Given \a space \( B \) and \a strict \two group \( \GG \),
the \two category \( \BB ^ \GG \) of \GGG \two transitions over \( \BB \)
(or equivalently, the category of principal \GGG \two bundles over \( \BB \))
is equivalent to the \two category of \GGG gerbes over \( B \)
(where \( \GG \) is identified with its corresponding crossed module).
\end {thm}
\par It's worth mentioning that the usual sort of nonabelian gerbe,
where \( \cat \GG 1 \) is \a group of automorphisms,
may not actually work in \a given category \( \C \).
(See \cite [8.1] {HDA5} for precise instructions
on how to build the \q {automorphism \two group} of \a group.)
For example, if \( \C \) is the category of finite-dimensional manifolds,
then the automorphism group of \a Lie group
has (in general) infinite dimensions,
so it's not an object of \( \C \).
One can either form \a space of only those automorphisms desired,
passing to crossed modules,
or else extend \( \C \) to include the necessary automorphism groups.
This is just one reason
why Baez \& Schreiber~\cite {HGT} would replace the category of manifolds
with \a more general category of \q {smooth spaces}.
\par Abelian gerbes, in contrast, are easy;
they use \a strict \two group where \( \cat \GG 1 \) is the trivial group.
As remarked in section~\ref {disc2space}, this requires an abelian group.
\part * {Letters used}
For reference, \I include tables of letters used and their meanings.
\par \emph {Uppercase Latin letters} normally refer
to spaces (in italics) and \two spaces (in calligraphics).
But an underlined uppercase Latin letter denotes an identity \too map,
and \a doubly underlined uppercase Latin letter
denotes the identity natural transformation of an identity \two map.
Also, \( \C \) is the category of spaces (section~\ref {1C}),
and \( \CCC \) is the \two category of \two spaces (section~\ref {2C}).
\smallbreak
\par \noindent \begin {tabular} {l c r}
Letter& Meanings& Reference sections\\
\hline
\( B \)& base space of \a bundle& \ref {1cover}, \ref {1bun}\\
\( \BB \)& base \two space of \a \two bundle& \ref {2cover}, \ref {2bun}\\
\( C \)& pullback cone& \ref {1pull}\\
\( \C \)& the category of spaces& \ref {1C}\\
\( \CC \)& pullback \two cone& \ref {2pull}\\
\( \CCC \)& the \two category of \two spaces& \ref {2C}\\
\( D \)& automorphism group of \a crossed module& \ref {2crosmod}\\
\( E \)& total space of \a bundle& \ref {1bun}\\
\( \EE \)& total \two space of \a \two bundle& \ref {2bun}\\
\( F \)& fibre of \a bundle, \G space&
\ref {1trivbun}, \ref {1loctrivbun}, \ref {1act}\\
\( \FF \)& fibre of \a \two bundle, \GGG \two space&
\ref {2trivbun}, \ref {2loctrivbun}, \ref {2act}\\
\( G \)& group& \ref {1group}\\
\( \GG \)& \two group& \ref {2group}\\
\( H \)& base group of \a crossed module& \ref {2crosmod}\\
\( K \)& index set (not \a space)& \ref {1opencover}\\
\( N \)& kernel of action of \( G \) on \( F \)& \ref {1Bun}\\
\( P \)& pullback& \ref {1pull}\\
\( \PP \)& \two pullback& \ref {2pull}\\
\( R \)& spaces of \a relation& \ref {1eqrel}\\
\( \RR \)& \two spaces of \a \two relation& \ref {2eqrel}\\
\( U \)& cover, subspace& \ref {1cover}, \ref {1restrict}\\
\( \UU \)& \two cover, \two subspace& \ref {2cover}, \ref {2restrict}\\
\( W , X , Y , Z \)& generic space\\
\( \XX , \YY , \ZZ \)& generic \two spaces
\end {tabular}
\smallbreak
\par \emph {Lowercase Latin letters}
normally refer to \too maps between \too spaces
(in italics and Fraktur, respectively).
(see section~\ref {disc2space}).
But an underlined lowercase Latin letter
denotes an identity natural transformation,
and \a lowercase Latin letter inside vertical bars
denotes \a domain \two space (see section~\ref {2map}).
\smallbreak
\par \noindent \begin {tabular} {l c r}
Letter& Meanings& Reference sections\\
\hline
\( b \)& \G transition morphism& \ref {1G^B}\\
\( \bb \)& \GGG \two transition morphism& \ref {2G^B}\\
\( d \)& structure maps in \a \two space or crossed module&
\ref {2C}, \ref {2crosmod}\\
\( e \)& identity in \a group& \ref {1group}\\
\( \ee \)& identity in \a \two group& \ref {2group}\\
\( f \)& \gee bundle morphism& \ref {1C/B}, \ref {1Bun}\\
\( \ff \)& \Gee \two bundle morphism& \ref {2C/B}, \ref {2Bun}\\
\( g \)& \G transition&\ref {1Gbun}\\
\( \gg \)& \GGG \two transition& \ref {2Gbun}\\
\( i \)& inverse in \a group& \ref {1group}\\
\( \ii \)& inverse in \a \two group& \ref {2group}\\
\( j \)& cover map, subspace map, relation map&
\ref {1cover}, \ref {1restrict}, \ref {1eqrel}\\
\( \jj \)& \two cover \two map, \two subspace \two map, \two relation \two map&
\ref {2cover}, \ref {2restrict}, \ref {2eqrel}\\
\( k \)& index (not \a map)& \ref {1opencover}\\
\( l \)& action of \( D \) on \( H \)& \ref {2crosmod}\\
\( m \)& multiplication in \a group& \ref {1group}\\
\( \mm \)& multiplication in \a \two group& \ref {2group}\\
\( n \)& natural number (not \a map)& \ref {2C}\\
\( p \)& projection map& \ref {1bun}\\
\( \pp \)& projection \two map& \ref {2bun}\\
\( r \)& action on \a right \G space& \ref {1act}\\
\( \rr \)& action on \a right \GGG \two space& \ref {2act}\\
\( t \)& \G map& \ref {1C^G}\\
\( \tt \)& \GGG \two map& \ref {2C^G}\\
\( u \)& pullback cone morphism& \ref {1pull}\\
\( \uu \)& \two pullback cone morphism& \ref {2pull}\\
\( w \)& generic, often an element, often of \( F \)\\
\( \ww \)& generic, often an element, often of \( \FF \)\\
\( x , y , z \)& generic, often an element, often of \( U \)\\
\( \xx , \yy , \zz \)& generic, often an element, often of \( \UU \)
\end {tabular}
\smallbreak
\par \emph {Lowercase Greek letters}
refer to equations between maps and natural transformations between \two maps.
But \a lowercase Greek letter inside vertical bars
denotes the underlying map of \a natural transformation
(see section~\ref {nat}).
\smallbreak
\par \noindent \begin {tabular} {l c r}
Letter& Meanings& Reference digrams\\
\hline
\( \alpha \)& associative law in \a group, associator in \a \two group&
(\ref {monoid}, \ref {2monoid})\\
\( \beta \)&
\begin {tabular} c
composition law for \G transition morphisms,\\
compositor of \GGG \two transition morphisms
\end {tabular} &
(\ref {composite transition morphism universal property},
\ref {composite 2transition morphism universal property})\\
\( \gamma \)&
multiplication law of \a \G transition,
multiplicator in \a \GGG \two transition&
(\ref {transition}, \ref {2transition})\\
\( \delta \)&
\begin {tabular} c
right action law of \a \G transition morphism,\\
right actor of \a \G \two transition morphism
\end {tabular} &
(\ref {transition morphism right}, \ref {2transition morphism right})\\
\( \epsilon \)& left inverse law in \a group, left invertor in \a \two group&
(\ref {inverse}, \ref {2inverse})\\
\( \zeta \)&
\begin {tabular} c
when \a bundle morphism between \G bundles is \a \G bundle morphism,\\
when \a \two bundle morphism between \GGG \two bundles
is \a \GGG \two bundle morphism
\end {tabular} &
(\ref {Gbundle morphism}, \ref {G2bundle morphism})\\
\( \eta \)&
identity law in \a \G transition, identitor in \a \GGG \two transition&
(\ref {transition}, \ref {2transition})\\
\( \theta \)&
\begin {tabular} c
when \a locally trivial bundle is \a \G bundle,\\
when \a locally trivial \two bundle is \a \GGG \two bundle
\end {tabular} &
(\ref {Gbundle}, \ref {G2bundle})\\
\( \iota \)& right inverse law in \a group, right invertor in \a \two group&
(\ref {inverse}, \ref {2inverse})\\
\( \kappa \)& equivalence of \gee bundles, \two morphism of \two \Gee bundles&
(\ref {bundle equivalence}, \ref {2bundle 2morphism})\\
\( \lambda \)& left unit law in \a group, left unitor in \a \two group&
(\ref {monoid}, \ref {2monoid})\\
\( \mu \)& associative law in \a \G space, associator in \a \GGG \two space&
(\ref {Gspace}, \ref {G2space})\\
\( \nu \)& universality of \a \two pullback, universality of \a \two quotient&
(\ref {2pullback universality}, \ref {2quotient 2universality})\\
\( \xi \)& \two morphism of \GGG \two transitions&
(\ref {2transition 2morphism})\\
\( \pi \)& when \a \too map between \too bundles is \a \too bundle morphism&
(\ref {bundle morphism}, \ref {2bundle morphism})\\
\( \rho \)& right unit law in \a group, right unitor in \a \two group&
(\ref {monoid}, \ref {2monoid})\\
\( \sigma \)&
\begin {tabular} c
left action law of \a \G transition morphism,\\
left actor of \a \G \two transition morphism
\end {tabular} &
(\ref {transition morphism left}, \ref {2transition morphism left})\\
\( \tau \)& equivalence of \G spaces, \two morphism of \GGG \two spaces&
(\ref {Gspace equivalence}, \ref {Gnattrans})\\
\( \upsilon \)& unit law in \a \G space, unitor in \a \GGG \two space&
(\ref {Gspace}, \ref {G2space})\\
\( \phi \)&
\begin {tabular} c
when \a map between \G spaces is \a \G map,\\
when \a \two map between \GGG \two spaces is \a \GGG \two map
\end {tabular} &
(\ref {Gmap}, \ref {G2map})\\
\( \chi , \psi \)& generic\\
\( \omega \)& \too pullback, \too quotient of an equivalence \too relation&
\begin {tabular} r
(\ref {pullback cone}ff, \ref {quotient cone}ff, \ref {cover kernel pair}ff,
\ref {restriction}, \ref {summary pullback},\\ \ref {2pullback cone}ff,
\ref {2quotient cone}ff, \ref {2cover 2kernel pair}ff,
\ref {2restriction}, \ref {summary 2pullback})
\end {tabular}
\end {tabular}
\smallbreak
\section * {Acknowledgements}
\I thank James Dolan for helpful discussions about category theory,
Miguel Carri\'on \'Alvarez for helpful discussions about geometry,
Aaron Lauda for helpful discussions about \two groups,
Urs Schreiber for helpful discussions about \two spaces,
and Danny Stevenson for helpful discussions about covers.
And of course \I thank John Baez for all of the above,
as well as inspiration, guidance, and encouragement.
\par \raggedright


\begin{thebibliography}{HDA6}

\bibitem[ACJ]{nonabbungerb}
Paolo Aschieri{,}~Luigi Cantini{,} and Branislav~Jur\v co.
\newblock \textsl {Nonabelian Bundle Gerbes, their Differential Geometry and
  Gauge Theory}.
\newblock {\em Communications In Mathematical Physics}, 254:367--400, 2005.
\newblock arXiv:hep-th/0312154.

\bibitem[Att]{Attal}
Romain Attal.
\newblock \textsl {Combinatorics of Non-Abelian Gerbes with Connection and
  Curvature}.
\newblock arXiv:math-ph/0203056.

\bibitem[Baez]{gauge.tex}
John~C. Baez.
\newblock \textsl {Higher Yang--Mills Theory}.
\newblock arXiv:hep-th/0206130.

\bibitem[BD]{categorification}
John~C. Baez and James Dolan.
\newblock \textsl {Categorification}.
\newblock In Ezra Getzler and Mikhail Kapronov, editors, {\em Higher Category
  Theory}.
\newblock Number 230 in Contemporary Mathematics. American Mathematical
  Society.
\newblock arXiv:math.QA/9802029.

\bibitem[B\'en]{bicat}
Jean B\'enabou.
\newblock {\em Introduction to Bicategories}.
\newblock Number~40 in Lecture Notes in Mathematics. Springer--Verlag, 1967.

\bibitem[BM]{BM}
Lawrence Breen and William Messing.
\newblock \textsl {Differential Geometry of Gerbes}.
\newblock arXiv:math.AG/0106083.

\bibitem[Bre90]{crosmodgerb}
Lawrence Breen.
\newblock \textsl {Bitorseurs et Cohomologie Non Ab\'elienne}.
\newblock In P.~Cartier et~al., editors, {\em The Grothendieck Festschrift},
  volume~I, pages 401--476.
\newblock Number~86 in Progress in Mathematics. 1990.

\bibitem[Bre94]{Asterisque}
Lawrence Breen.
\newblock {\em On the Classification of \Two Gerbes and \Two stacks}.
\newblock Number 225 in Ast\'erisque. Soci\'et\'e Math\'ematique de France,
  1994.

\bibitem[CF]{CF}
Louis Crane and I.~Frenkel.
\newblock \textsl {Four-Dimensional Topological Quantum Field Theory, Hopf
  Categories, and the Canonical Bases}.
\newblock {\em Journal of Mathematical Physics}, 35, 1994.

\bibitem[CWM]{CWM}
Saunders~Mac Lane.
\newblock {\em Categories for the Working Mathematician}.
\newblock Number~5 in Graduate Texts in Mathematics. Springer--Verlag, 2nd
  edition, 1998.

\bibitem[Dusk]{Duskin}
J.~Duskin.
\newblock \textsl {An Outline of a Theory of Higher Dimensional Descent}.
\newblock {\em Bulletin de la Soci\'et\'e Math\'ematique de Belgique
  S\'eries~A}, 41(2):249--277, 1989.

\bibitem[Ehr]{Ehresmann}
Charles Ehresmann.
\newblock \textsl {Categories Structur\'ees}.
\newblock {\em Annales de l'\'Ecole Normale et Superieure}, 80, 1963.

\bibitem[Ele1]{Elephant}
Peter~T. Johnstone.
\newblock {\em Sketches of an Elephant: A Topos Theory Compendium}, volume~1.
\newblock Number~43 in Oxford Logic Guides. Oxford University Press, 2002.

\bibitem[Euc]{Euclid}
Euclid.
\newblock Thomas~L. Heath, editor.
\newblock {\em Elements}.
\newblock The Perseus Digital Library.
\newblock \texttt
  {http://www.perseus.tufts.edu/cgi-bin/ptext?doc=Perseus:text:1999.01.0086:id%
=elem.1.c.n.1}.

\bibitem[F-B]{crosmod}
Magnus Forrester-Barker.
\newblock \textsl {Group Objects and Internal Categories}.
\newblock arXiv:math.CT/0212065.

\bibitem[GP]{HP}
Florian Girelli and Hendryk Pfeiffer.
\newblock \textsl {Higher Gauge Theory --- Differential versus Integral
  Formulation}.
\newblock arXiv:hep-th/0309173.

\bibitem[HCA1]{Borceux}
F.~Borceux.
\newblock {\em Handbook of Categorical Algebra}, volume 1: \textsl {Basic
  Category Theory}.
\newblock Number~50 in Encyclopedia of Mathematics and its Applications.
  Cambridge University Press, 1994.

\bibitem[HDA5]{HDA5}
John~C. Baez and Aaron~D. Lauda.
\newblock \textsl {Higher-Dimensional Algebra V: \Two groups}.
\newblock arXiv:math.QA/0307200.

\bibitem[HDA6]{HDA6}
John~C. Baez and Alissa~S. Crans.
\newblock \textsl {Higher-Dimensional Algebra VI: Lie \Two algebras}.
\newblock arXiv:math.QA/0307263.

\bibitem[HGT]{HGT}
John~C. Baez and Urs Schreiber.
\newblock \textsl {Higher Gauge Theory}.
\newblock arXiv:math.DG/0511710.

\bibitem[HGT2]{HGT2}
John~C. Baez and Urs Schreiber.
\newblock \textsl {Higher Gauge Theory II: \Two connections on \Two bundles}.
\newblock arXiv:hep-th/0412325.

\bibitem[JST]{JST}
George Janelidze{,}~Manuela Sobral{,} and Walter Tholen.
\newblock \textsl {Beyond Barr Exactness: Effective Descent Morphisms}.
\newblock In Maria~Cristina Pedicchio and Walter Tholen, editors, {\em
  Categorical Foundations: Special Topics in Order, Topology, Algebra, and
  Sheaf Theory}, chapter VIII.
\newblock Number~97 in Encyclopedia of Mathematics and its Applications.
  Cambridge University Press, 2004.

\bibitem[Jur]{crosmodbungerb}
Branislav~Jur\v co.
\newblock \textsl {Crossed Module Bundle Gerbes; Classification, String Group
  and Differential Geometry}.
\newblock arXiv:math.DG/0510078.

\bibitem[Lang]{Lang}
Serge Lang.
\newblock {\em Algebra}.
\newblock Addison--Wesley, 3rd edition, 1995.

\bibitem[Mak]{Makkai}
Michael Makkai.
\newblock \textsl {Avoiding the Axiom of Choice in General Category Theory}.
\newblock {\em Journal of Pure and Applied Algebra}, 108:109--173, 1996.

\bibitem[Met]{surjsubm}
David~S. Metzler.
\newblock \textsl {Topological and Smooth Stacks}.
\newblock arXiv:math.DG/0306176.

\bibitem[ML]{Noether}
Saunders~Mac Lane.
\newblock \textsl {Topology Becomes Algebraic with Vietoris and Noether}.
\newblock {\em Journal of Pure and Applied Algebra}, 39(3):305--307, 1986.

\bibitem[Moer]{Moerdijk1228}
Ieke Moerdijk.
\newblock \textsl {On the Classification of Regular Lie groupoids}.
\newblock Preprint 1228, Department of Mathematics, Utrecht University, 2002.

\bibitem[Mur]{bungerb}
Michael~K. Murray.
\newblock \textsl {Bundle Gerbes}.
\newblock {\em Journal of the London Mathematical Society 2}, 54:403--416,
  1996.
\newblock arXiv:dg-ga/9407015.

\bibitem[Sch]{Urs}
Urs Schreiber.
\newblock {\em From Loop Space Mechanics to Nonabelian Strings}.
\newblock PhD thesis, Universit\"at Duisburg-Essen, 2005.
\newblock arXiv:hep-th/0509163.

\bibitem[Spi]{submanifold}
Michael Spivak.
\newblock {\em A Comprehensive Introduction to Differential Geometry},
  volume~1.
\newblock Publish or Perish, Inc., 2nd edition, 1979.

\bibitem[Ste]{Steenrod}
Norman Steenrod.
\newblock {\em The Topology of Fibre Bundles}.
\newblock Princeton University Press, 1951.

\bibitem[Str]{string}
Ross Street.
\newblock \textsl {Low-Dimensional Topology and Higher-Order Categories}.
\newblock Halifax, July 9--15 1995. International Category Theory Meeting.

\bibitem[Weyl]{Weyl}
Hermann Weyl.
\newblock {\em The Classical Groups: Their Invariants and Representations}.
\newblock Princeton University Press, 1939.

\end{thebibliography}
\end {document}